\documentclass[11pt,a4paper]{article}

\usepackage{amsmath,amssymb}
\newcommand{\R}{\mathbb{R}}
\newcommand{\D}{\Delta}

\newcommand{\p}{\partial}

\newtheorem{cor}{Corollary}

\newcommand{\e}{\epsilon}
\newcommand{\va}{\varphi}

\newtheorem{definition}{Definition}
\newtheorem{theorem}{Theorem}

\newtheorem{proposition}{Proposition}

\newtheorem{remarka}{Remark}

\newtheorem{lemme}{Lemma}

\title{Vortex solutions for the compressible Navier-Stokes equations with general viscosity coefficients in 1D: regularizing effects or not on the density}
\author{Boris Haspot  \thanks{Universit\'e Paris Dauphine, PSL Research University, Ceremade, Umr Cnrs 7534, Place du Mar\' echal De Lattre De Tassigny 75775 Paris cedex 16 (France), haspot@ceremade.dauphine.fr } \thanks{ANGE project-team (Inria, Cerema, UPMC, CNRS), 2 rue Simone Iff, CS 42112, 75589 Paris, France.} }
\date{}
\begin{document}
%\tableofcontents
\maketitle

%The title of your section 1
%Invariance du scaling nous donne ici:
%$$(\rho_{0},u_{0})\in B^{\N+\frac{2}{2\lambda-\gamma-1}}_{2,\infty}\times B^{\N-1+2\frac{\lambda-1}{2\lambda-\gamma-1}}_{2,\infty}.$$
% Text of your Version française abrégée here.
% Note you do not need to repeat here equations that you use in the
% main text - for example 'voir (3)' is quite acceptable.

%\selectlanguage{english}
% main text
% It is also important to point out that such 
%Mwith large initial data for compressible Navier Stokes equation with viscosity coefficients of the form $\mu(\rho)=\mu\rho^\alpha$ with $0\leq\alpha$ (it includes in particular the important physical case of the viscous shallow water system when $\alpha=1$ and of the case of constant viscosity coefficients). More precisely we deal with initial data $(\rho_0,u_0)$ which are in $BV(\R)\times{\cal M}(\R)$.
\begin{abstract}
We consider Navier-Stokes equations for compressible viscous fluids in the one-dimensional case with general viscosity coefficients. We prove the existence of global weak solution  when the initial momentum $\rho_0 u_0$ belongs to the set of the finite measure ${\cal M}(\R)$ and when the initial density $\rho_0$ is in the set of bounded variation functions $BV(\R)$.  In particular it allows to deal 
with initial momentum which are Dirac masses  and initial density which admit shocks.  We can observe in particular that this type of initial data have infinite energy. Furthermore we show that if the coupling between the density and the velocity is sufficiently strong then the initial density which admits initially shocks is instantaneously regularized and becomes continuous. This coupling is expressed via the regularity of the so called effective velocity $v=u+\frac{\mu(\rho)}{\rho^2}\p_x \rho$ with $\mu(\rho)$ the viscosity coefficient. Inversely if the coupling between the initial density and the initial velocity is too weak (typically $\rho_0 v_0\in{\cal M}(\R)$) then we prove the existence of weak energy solution in finite time but the density remains a priori discontinuous on the time interval of existence.%The key ingredient of the proof relies to a new formulation of the compressible equations involving a new effective velocity $v$ (see \cite{cras,para,CPAM,CPAM1}) such that the density verifies a parabolic equation. We estimate $v$ in $L^\infty_{t,x}$ norm which enables us to control the $L^\infty_{t,x}$ norm of $\frac{1}{\rho}$ by using the maximum principle.
\end{abstract}
\section{Introduction}
In this paper we wish to investigate the existence of global weak solutions of the following Navier-Stokes equations for compressible isentropic flow:
\begin{equation}
\begin{cases}
\begin{aligned}
&\p_t\rho+\p_x(\rho u)=0,\\
&\p_t(\rho u)+\p_x(\rho u^2)-\p_x(\mu (\rho)\p_x u)+\p_x P(\rho)=0.
\end{aligned}
\end{cases}
\label{1}
\end{equation}
with possibly degenerate viscosity coefficient $\mu(\rho)\geq 0$. Here $u=u(t,x)\in\R$ stands for the velocity field, $\rho=\rho(t,x)\in\R^{+}$ is the density and $(t,x)\in\R^+\times\R$. 
Throughout the paper, we will assume that the pressure $P(\rho)$ verifies a $\gamma$ type law
$P(\rho)=a \rho^\gamma$ with $a>0$ and $\gamma>1$ and that the viscosity coefficient can be written under the form $\mu(\rho)=\mu\rho^\alpha$ with $\mu>0$ and $\alpha\geq 0$. A large amount of literature is dedicated to the study of the compressible Navier-Stokes equations with
 constant viscosity case, however physically the viscosity of a gas depends on the temperature and  on the density (in the isentropic case). Let us mention the case of the Chapman-Enskog viscosity law (see \cite{CC70}) or the case of monoatomic gas ($\gamma=\frac{5}{3}$) where $\mu(\rho)=\rho^{\frac{1}{3}}$. More generally, the viscosity coefficient $\mu(\rho)$ is expected to vanish as a power of the density $\rho$ on the vacuum. We emphasize that the case $\alpha=1$ corresponds to the so called viscous shallow water system. This system with friction has been derived by Gerbeau and Perthame in \cite{GP} from the Navier-Stokes system with a free moving boundary in the shallow water regime at the first order (it corresponds to a small shallowness parameter). This derivation relies on the hydrostatic approximation where the authors follow the role of viscosity and friction on the bottom. 
\\[1mm]
In the case of compressible Navier-Stokes equations in one dimension, the existence of global weak solutions was first obtained by Kazhikhov and Shelukin \cite{KS77} for smooth enough data close to the equilibrium (in particular the initial density $\rho_0$ is bounded away from zero) with constant viscosity coefficient. The case of initial density admitting shocks has been treated by Hoff \cite{Hoff86} where the initial density belongs to $BV(\R)$, in addition the author needs smallness assumption on the initial data. In \cite{Hof98}, Hoff proved the existence of global weak energy solution for constant viscosity coefficients provided that $\rho_0$ is only bounded in $L^\infty$ norm and is far away from zero. In this work there is no smallness restriction on the initial data. In passing let us mention that the existence of global weak solution for constant viscosity coefficients in any dimension $N\geq 2$ has been proved for the first time by Lions in \cite{Lio98} and the result has been later refined by Feireisl et al \cite{FNP01}. The existence of global strong solution in one dimension with large initial data for initial density far away from the vacuum has been proved for the first time by Kanel \cite{Ka}  (see also \cite{Hof87}) in the case of constant viscosity coefficients.
\\
The study of the compressible Navier Stokes equations with degenerate viscosity coefficients is more recent and has been in particular motivated by the introduction of a new entropy (see \cite{BD1}) which provides new regularity estimate on the gradient of the density. In \cite{Jiu} Jiu and Xin proved the existence of global weak energy solution in one dimension when $\mu(\rho)=\mu\rho^\alpha$ with $\alpha>\frac{1}{2}$ and with large initial energy data.
In \cite{MV}, Mellet and Vasseur showed the existence of global strong solution when $0\leq\alpha<\frac{1}{2}$ for large initial data provided that the initial data $(\rho_0-\bar{\rho},u_0)$ are in $H^1(\R)\times H^1(\R)$ with $0<c\leq\rho_0\leq M<+\infty$. In \cite{Glob} we extends this result to the case $\frac{1}{2}<\alpha\leq 1$ including in particular the case of the so called viscous shallow water system ( $\alpha=1$, see \cite{GP}). Constantin et al in \cite{PC} propose a very interesting other method to prove the existence of global strong solution for $0\leq\alpha\leq 1$ and also obtain new results in the case $\alpha>1$.\\
In passing we would like to give few words on the existence of global weak solution for degenerate viscosity coefficients when $N$ the dimension verifies $N\geq 2$. Mellet and  Vasseur proved in \cite{MV06} the stability of the global weak solution for compressible Navier-Stokes equation with viscosity coefficient verifying the so called BD entropy (see \cite{BD1}) in dimension $N= 2,3$ (we refer also to \cite{BD1,BrDeZa} when we consider in addition friction terms). Let us mention in particular that the case $\mu(\rho)=\mu\rho$ with $\mu>0$ and $\lambda(\rho)=0$ verifies the algebraic relation related to the new entropy discovered in \cite{BD1}, it corresponds here to the so called viscous shallow water system.
For $N=2,3$ the important problem of the existence of global weak solutions has been recently resolved independently by Vasseur and Yu \cite{V,V1} and Li and Xin in \cite{Li}.\\
We would like to emphasize that all these results have been achieved in a framework with finite energy initial data, indeed systematically $\sqrt{\rho_0}u_0$ belongs to $L^2(\R)$. It will be not the case of our result since $\rho_0 u_0, \rho_0v_0$ are  in our case only finite measures or a vortex. The existence of global solution for initial data which are finite measure is an important problem in mathematics and has been solved for other equations. We recall in particular that for the incompressible Navier-Stokes equations it exists a theory on the existence of global strong solution in dimension $N=2$ for initial vorticity which belongs to the set of finite measure (see for example \cite{GM,GG}). It exists also a theory for the existence of global weak solution for the incompressible Euler equations in dimension $N=2$ when the initial vorticity is a measure (see \cite{Delort,Berto}).\\
Following \cite{para,Glob,BrDeZa,CPAM,CPAM1} we can now rewrite the system (\ref{1}) by using the so called effective %momentum $m^1=m+\p_x\va_1(\rho)$ with $\va_1'(\rho)=\frac{\mu(\rho)}{\rho}$:
 velocity $v=u+\p_x \va(\rho)$ with $\va'(\rho)=\frac{\mu(\rho)}{\rho^2}$:
%\begin{equation}
%\begin{cases}
%\begin{aligned}
%&\p_t\rho-\p_{x}(\frac{\mu(\rho)}{\rho}\p_x\rho)+\p_x(\rho v)=0,\\
%& \p_t m^1+\p_x(\frac{m m^1}{\rho})+a\gamma\frac{\rho^{\gamma}}{\mu(\rho)}m^1=a\gamma\frac{\rho^{\gamma}}{\mu(\rho)}m.
%\end{aligned}
%\end{cases}
%\label{11}
%\end{equation}
\begin{equation}
\begin{cases}
\begin{aligned}
&\p_t\rho-\p_{x}(\frac{\mu(\rho)}{\rho}\p_x\rho)+\p_x(\rho v)=0,\\
&\p_t (\rho v)+\p_x(\rho u v)+a\gamma\frac{\rho^{\gamma+1}}{\mu(\rho)}v=a\gamma\frac{\rho^{\gamma+1}}{\mu(\rho)}u\\
&\p_t(\rho u)+\p_x(\rho u^2)-\p_x(\mu (\rho)\p_x u)+a\gamma\frac{\rho^{\gamma+1}}{\mu(\rho)}v=a\gamma\frac{\rho^{\gamma+1}}{\mu(\rho)}u
\end{aligned}
\end{cases}
\label{11}
\end{equation}
It is important to note that this change of unknown is true for any viscosity coefficients (at the opposite in dimension $N\geq 2$ such change of unknown involves an algebraic relation between the coefficients $\mu$ and $\lambda$ which precludes in particular the case of constant viscosity coefficients (see \cite{para,BrDeZa} for more details)). When we observe the form of the system (\ref{11}), it is then natural to consider initial momentum $\rho_0u_0$ and $\rho_0 v_0$ which are in ${\cal M}(\R)$; indeed if we choose  $\rho_0u_0$ and $\rho_0 v_0$ in $L^1(\R)$ it is easy to prove that the $L^1$ norm is preserved all along the time. The problem of working with initial data describing some vortex seems then natural.\\
The questions addressed in this paper are also motivated by the fact that discontinuous solutions are fundamental both in the physical theory of non-equilibrium thermodynamics as well as in the mathematical study of inviscid models for compressible flow. It is important, therefore, to express a rigorous theory which accommodates these discontinuities or shocks in the theory of compressible Navier Stokes equations. In order to obtain such results, a key point is to work with initial data with minimal regularity assumption as much as it is possible to do. The results of this paper effectively give a first direction to the existence of global weak solution with initial density admitting shocks. More precisely we should deal with initial density $\rho_0$ in $BV(\R)$   with $0<c\leq\rho_0\leq M<+\infty$ and initial momentum $\rho_0u_0, \rho_0 v_0$ in ${\cal M}(\R)$ the set of finite real measure (including in particular Dirac masses).  If $\mu_1\in{\cal M}(\R)$, the total variation of $\mu_1$ is defined by:
$$\|\mu_1\|_{{\cal M}(\R)}=\sup\{\int_{\R}\varphi d\mu_1\;/ \;\varphi\in C_0(\R),\;\|\varphi\|_{L^\infty(\R)}\leq 1\;\}$$
where $C_0(\R)$ is the set of all real-valued functions on $\R$ vanishing at infinity. In addition we say that $\rho_0$ is in $BV(\R)$ if $\rho_0$ is measurable, locally integrable, continuous on the left and such that the derivative in the sense of the distribution is a finite measure. An other way to describe the set $BV(\R)$ is to consider the function $\rho_0$ with finite total variation:
$$TV(\rho_0)=\sup\{ \sum_{j=1}^N |\rho_0(x_j)-\rho_0(x_{j+1})|\;N\in\mathbb{N}^*,\;-\infty<x_0<\cdots<x_N<+\infty\}.$$
The set $BV(\R)$ of these functions is a Banach space when we equip it with the norm $TV(\rho_0)+|\rho_0(y)|$ with $y\in\R$. The set $BV(\R)$ is particularly adapted to our study since it is sufficiently large for admitting shocks. It is also the natural space fas we will see for studying the hyperbolic system in one dimension.\\
We will observe that if in addition we assume that $v_0$ and $\rho_0-\bar{\rho}$ belong to $L^2(\R)$ then the density is instantaneously regularizing inasmuch as the density $\rho(t,\cdot)$ becomes continuous for $t>0$.\\
It is interesting to compare this kind of surprising phenomenon with what is known on the behavior of the solutions of hyperbolic systems (typically the compressible isentropic Euler system which corresponds to $\mu(\rho)=0$).
In particular when the system is strictly hyperbolic (for a $N\times N$ system with $x\in\R$) and genuinely non linear, we know from the Glimm theorem  (see \cite{Glimm}) that it exists a global weak solution with small initial data in $(BV(\R))^N$ (see also \cite{Dafermos} for some extension in the case of a $2\times2$ system where the smallness assumption is weaken due to the existence of Riemann invariants). Concerning the uniqueness of such solutions, we refer to the works of Bressan et al (see \cite{Bressan} for a excellent review on all these results using in particular the so-called front tracking methods which generate a contractive $L^1$ semi-group). It is important to note that for such hyperbolic systems there is a priori no regularizing effects on the unknown in the sense that the solution admits shocks all along the time (however a solution with a $L^\infty$ initial data can become $BV(\R)$ instantaneously, we refer to the scalar case with strictly convex flux which was in particular studied by Oleinik \cite{Oleinik}). %This explains why in particular the $BV(\R)$ space is suitable for the study of the one dimensional case since this space contains functions with shocks. The simpler case is the solution of the Riemann problem due to Lax \cite{Lax}. 
In particular it is known that for regular general initial data, there is creation of shocks in finite time (except if we can use a characteristic methods, typically in the case of $2\times2$ systems where we have Riemann invariants). It is then remarkable to observe that in the case of compressible Navier Stokes systems the density can be instantaneously regularized even if the density seems to be governed by a transport equation in (\ref{1}).
%We emphasize in addition that the introduction of the effective velocity $v$ allows to transform the system (\ref{1}) into a parabolic equation on the density and a transport equation on the velocity (we will see in the sequel that $v$ is not so far to verify a damped transport equation). It seems surprising to observe that contrary to $u$ which has a parabolic behavior, $v$  has a hyperbolic behavior. Roughly speaking, the compressible Navier Stokes equations in the one dimensional-case can be seen as the compressible Euler equations with a viscous regularizing term on the density of the type $-\p_{x}(\frac{\mu(\rho)}{\rho}\p_x\rho)$. Note that it is less obvious in dimension $N\geq 2$ since the momentum equation of (\ref{11}) has in addition a term of the form ${\rm div}(\mu(\rho){\rm curl}v)$, see \cite{para}).\\
 \\
% \\
%We are going now to recall some results on the existence of solutions for the one-dimensional case when the viscosity coefficient is constant positive. 

\section{Main result}
We are going now to state the main results of this paper. Before we would like to recall a very interesting result due to Hoff, this one indicates that there is a priori no regularizing effects on the density (in the sense that the density does not become instantaneously continuous)  when the effective velocity $v$ is not sufficiently regular.
\begin{theorem}
[Hoff \cite{Hof98}] 
Let $\bar{\rho}>0$ and $\bar{u}\in\R$. Assume that $(\rho_0,\frac{1}{\rho_0})\in L^\infty(\R)$ and $(\rho_0-\bar{\rho},u_0-\bar{u})\in (L^2(\R))^2$ then the initial value problem (\ref{1}) has a global weak solution $(\rho,u)$ for which
\begin{equation}
\begin{cases}
\begin{aligned}
&\rho-\bar{\rho},\,\rho u \in C([0,+\infty),H^{-1}(\R)),\\
&u\in C((0,+\infty),L^2(\R)),\\
&u(t,\cdot), \mu \p_x u(t,\cdot)-P(\rho(t,\cdot))+P(\bar{\rho})\in H^1(\R),\,t>0,\\
&\p_t u(t,\cdot),\,\dot{u}(t,\cdot)\in L^2(\R),\;t>0,
\end{aligned}
\end{cases}
\end{equation}
where $\dot{u}$ is the convective derivative $\dot{u}=\p_t u+u\p_x u$. Additionally, given $T>0$, there is a positive constant $C(T)$ depending on $T$, $\mu$, $P$ and on upper bounds for $\|\rho_0-\bar{\rho}\|_{L^2(\R)}$, $\|u_0\|_{L^2(\R)}$, $\|\rho_0\|_{L^\infty}$, $\|\rho_0^{-1}\|_{L^\infty}$ such that if $\sigma(t)=\min(1,t)$, then:
\begin{equation}
C(T)^{-1}\leq\rho(t,\cdot)\leq C(T)\;\;\mbox{a.e},
\label{1.17}
\end{equation}
\begin{equation}
\begin{aligned}
&\sup_{0\leq t\leq T}[\|\rho(t,\cdot)-\bar{\rho}\|_{L^2(\R)}+\|u(t)\|_{L^2(\R)}+\sigma(t)^{\frac{1}{2}}\|\p_x u(t,\cdot)\|_{L^2(\R)}\\
&+\sigma(t)\big(\|\dot{u}(t)\|_{L^2}+\|\mu \p_x u(t,\cdot)-P(\rho)(t,\cdot)+P(\bar{\rho})\|_{L^2(\R)}\big)\leq C(T),
\end{aligned}
\label{1.18}
\end{equation}
\begin{equation}
\begin{aligned}
&\int^T_0\big(\|\p_x u(s,\cdot)\|_{L^2}^2+\sigma(s)\|\dot{u}(s,\cdot)\|_{L^2}^2+\sigma(s)\|\mu \p_x u(s,\cdot)-P(\rho)(s,\cdot)+P(\bar{\rho})\|_{L^2(\R)}\\
&\hspace{8cm}+\sigma^2(s)\|\p_x\dot{u}(s,\cdot)\|_{L^2}^2]ds\leq C(T),
\end{aligned}
\label{1.19}
\end{equation}
and for $0<\tau<T$
\begin{equation}
\sigma(\tau)^{\frac{1}{4}}\|u\|_{L^\infty(\R\times|\tau,T])}+\sigma(\tau)^{\frac{1}{2}}\langle u\rangle^{\frac{1}{2},\frac{1}{4}}_{[\tau,T]\times\R}\leq C(T).
\label{1.20}
\end{equation}
where $\langle u\rangle^{\frac{1}{2},\frac{1}{4}}_{[\tau,T]\times\R}$ is the usaul H\"older norm
$$\sup\{\frac{|u(t,x)-u(s,y)|}{|x-y|^{\frac{1}{2}}+|t-s|^{\frac{1}{4}}}, x,y\in\R,\;t,s\in[\tau,T],\,(t,x)\ne (s,y) \}.$$
\label{theoHoff}
\end{theorem}
\begin{remarka}
In the case of the constant viscosity coefficient $\mu>0$, the effective initial velocity has the form
$v_0=u_0-\mu\p_x(\frac{1}{\rho_0})$ and under the assumptions of the theorem \ref{theoHoff} we have $v_0\in L^2(\R)+W^{-1,\infty}(\R)$. From (\ref{11}), we observe that the effective velocity $v$ verifies a damped transport equation which prevents a priori any regularizing effects. It implies that
$\p_x(\frac{1}{\rho})$ is at the best in $L^\infty(W^{-1,\infty}+L^2)$, in other words it is not sufficient 
to show via Sobolev embedding that the density $\rho$ becomes continuous. In our next theorems, the situation will be different since $v_0$ will be in $L^2(\R)$.
\end{remarka}
\begin{remarka}
In the previous theorem Hoff proved that the velocity $\p_x u-P(\rho)+P(\bar{\rho})$ becomes regular, however to be able to prove that $P(\rho)$ becomes regular it is necessary to obtain additional information on $\p_x u$. It explains also why we have no regularizing effects. The coupling with the effective velocity $v$ is different since it provides an information on $\p_x( \frac{1}{\rho})$. It explains why this coupling is stronger compared with the so called effective pressure.
\end{remarka}
%\begin{remarka}
%FAIRE UN LIEN ENTRE ETRE VITESSE LIPSCHITZ ET LA METHODE DES CARACTERISTIQUES, CF DANCHIN!
%\end{remarka}
We define now the initial momentum and the initial effective momentum which formally verify $m_0=\rho_0 u_0$ and $m^1_0=m_0+\p_x \va_1(\rho_0)$ (($m_0$ has to be considered as an unknown, $\rho_0$ corresponds to the initial density and we have $\va_1(\rho_0)=\frac{\mu}{\alpha}\rho_0^\alpha$).\\
Let us give a definition of a global weak solution for initial data verifying for $\bar{\rho}>0$:
\begin{equation}
\begin{cases}
\begin{aligned}
%&\rho_0\geq 0, m_0=0\;\;\mbox{on}\;\;\{x\in\R;\,\rho_0(x)=0\},\\
&0<c\leq\rho_0\leq M<+\infty\\
&(\rho_0-\bar{\rho})\in L^2(\R),\;\va_1(\rho_0)\in BV(\R)\\
&m^1_0\in {\cal M}(\R), \frac{m_0^1}{\rho_0}\in L^2(\R).
\end{aligned}
\end{cases}
\end{equation}
In the sequel we will define $m$, $m^1$, $m^2$ as follows:
$$m=\rho u,\;m^1=\rho u+\p_x\va_1(\rho)\;\;\mbox{and}\;\;m^2=\sqrt{\rho}u+\p_x\va_2(\rho),$$
with $\va_2(\rho)=\frac{\mu}{\alpha-\frac{1}{2}}\p_x\rho^{\alpha-\frac{1}{2}}$ with $\alpha\ne\frac{1}{2}$. We refer to \cite{Lio98} for the definition of the Orlicz space $L^\gamma_2(\R)$.
\begin{definition}
\label{def}
A pair $(\rho,u)$ is said to be a weak solution to (\ref{1}) provided that:
\begin{itemize}
\item $\rho\geq 0$ a. e and for $\bar{\rho}>0$, $\e>0$ and any $T>0$, $1\leq p<+\infty$:
$$
\begin{cases}
\begin{aligned}
&\rho-\bar{\rho}\in L^{\infty}([0,T],L^\gamma_2(\R)) ,\\
&\rho\in C([0,+\infty[,W^{-\e,p}_{loc}(\R)),\,\rho\in L^\infty([0,T],L^\infty(\R))\\
&\rho u\in L^{\infty}((0,T),L^1(\R)), \sqrt{\rho}u\in L^{2+\e}_{loc}([0,T]\times\R)\\
&\p_x\rho^{\frac{1}{2}(\gamma+\alpha-1)}\in L^2([0,T],L^2(\R)), \p_x(\rho^{\alpha-\frac{1}{2}})\in L^{2}(([0,T],L^2(\R))\\
&m^2\in L^\infty([0,T],L^2(\R))\\
&\p_x\rho^{\alpha-\frac{1}{2}}\in L^\infty([1,+\infty[,L^2(\R))\\%\sqrt{\rho}v\in L^\infty_T((0,T),L^2(\R)),\\
&m^1\in L^\infty((0,T),{\cal M}(\R)).
\end{aligned}
\end{cases}
$$
%with $\e>0$ sufficiently small, $\sqrt{\rho}v=\p_x\va_3(\rho)+\sqrt{\rho}u$ with $\va_3'(\rho)=\sqrt{\rho}\va'(\rho)$.
\item $(\rho,\sqrt{\rho}u)$ satisfies in the sense of distributions:
$$
\begin{cases}
&\p_t\rho+\p_x(\sqrt{\rho}\sqrt{\rho}u)=0\\
&\rho(0,x)=\rho_0(x).
\end{cases}
$$
%For any $t_2\geq t_1\geq 0$ and any $\varphi\in C_0^{\infty}([0,T]\times\R$:
%\begin{equation}
%\int_{\R}\rho(t_2,x)\varphi(t_2,x) dx-\int_{\R}\rho(t_1,x)\varphi(t_1,x) dx=\int^{t_2}_{t_1}\int_{\R}(\rho(s,y)\p_s \varphi(s,y)+\rho u\p_y(s,y)) ds dy.
%\label{2.3}
%\end{equation}
\item  For any $\varphi\in C_0^{\infty}([0,T[\times\R)$ and any $T>0$:
\begin{equation}
\begin{aligned}
&\int_{\R} m_0(x)\varphi(0,x) dx+\int^T_0\int_{\R}[\sqrt{\rho}(\sqrt{\rho}u)\p_t\varphi+((\sqrt{\rho}u)^2+\rho^\gamma)\p_x\varphi] dt dx\\
&+\langle\mu(\rho)\p_x u,\p_x\varphi\rangle=0,
\end{aligned}
\label{2.4}
\end{equation}
where the diffusion term makes sense when written as:
\begin{equation}
\begin{aligned}
&\langle\mu(\rho)\p_x u,\p_x\varphi\rangle=\int^T_0\int_{\R}\p_x(\va_3(\rho)) \sqrt{\rho} u\p_x\varphi dx dt+\int^T_0\int_{\R}\frac{\mu(\rho)}{\sqrt{\rho}} (\sqrt{\rho} u)\p_{xx}\varphi dx dt,
%-\int^T_0\int_{\R}\mu\rho^{\alpha-\frac{1}{2}}\sqrt{\rho}u\p_x\varphi\, dx dt-\frac{2\mu\alpha}{2\alpha-1}\int^T_0\int_{\R}\p_x(\rho^{\alpha-\frac{1}{2}})\sqrt{\rho}u\varphi\, dx dt.
\end{aligned}
\label{2.5}
\end{equation}
with $\va_3'(\rho)=\frac{\mu'(\rho)}{\sqrt{\rho}}$.
\end{itemize}
\end{definition}
We obtain the main following theorems. 
\begin{theorem}
Assume that $\mu(\rho)=\mu\rho^\alpha$ with $\mu,\alpha> 0$, $\alpha\ne\frac{1}{2}$
and  $P(\rho)=a\rho^{\gamma}$ with $\gamma>1$, $a>0$, $\gamma\geq \alpha$ and $\gamma\geq 2\alpha-1$ if $\alpha>\frac{1}{2}$.
The initial data satisfy for $\bar{\rho}>0$:
%\begin{equation}
%\begin{cases}
%\begin{aligned}
%&m_0\in {\cal M}(\R),\\
%&0<c\leq\rho_0(x)\leq C<+\infty\;\;\;\mbox{for any}\;x\in\R\\
%&\rho_0-\bar{\rho}\in L^{\gamma}_{2}(\R),\;\p_x \va_1(\rho_0)\in {\cal M}(\R)\\
%&m_0^1\in L^2(\R).
%\end{aligned}
%\end{cases}
%\label{initial}
%\end{equation}
%
\begin{equation}
\begin{cases}
\begin{aligned}
&\va_1(\rho_0)\in BV(\R),\;(\rho_0-\bar{\rho})\in L^{\gamma}_{2}(\R),\\
%&\rho_0 v_0\in {\cal M}(\R),\,\rho_0 u_0=m_0\in {\cal M}(\R),\;v_0\in L^2(\R)\\
&0<c\leq\rho_0(x)\leq C<+\infty\;\;\;\mbox{for any}\;x\in\R,\\
&v_0\in L^2(\R), \;\rho_0v_0\in{\cal M}(\R).
%&\rho_0-\bar{\rho}\in L^{\gamma}_{2}(\R).
\end{aligned}
\end{cases}
\label{initial}
\end{equation}
The momentum $m_0$ is defined as follows:
\begin{equation}
m_0=\rho_0v_0-\p_x\va_1(\rho_0),\,m_0\in{\cal M}(\R).
\label{initial1}
\end{equation}
In addition theres exists $\e_0>0$ such that if:
%\begin{equation}
%\|\rho_0v_0\|_{{\cal M}(\R)}+\|\rho_0u_0\|_{{\cal M}(\R)}\leq \e_0,
%\label{petitesse}
%\end{equation}
\begin{equation}
\|\p_x\va_1(\rho_0)\|_{{\cal M}(\R)}+\|m_0\|_{{\cal M}(\R)}\leq \e_0,
\label{petitesse}
\end{equation}
then there exists a global weak solution $(\rho,\sqrt{\rho}u)$ for the system (\ref{1}) verifying the definition \ref{def}.
In addition there exists $T_\beta>0$ and $C>0$ such that we have:
\begin{equation}
\begin{cases}
\begin{aligned}
&\|(\rho,\frac{1}{\rho})\|_{L^\infty([0,T_\beta],L^\infty(\R)}\leq C,\\
&\|u\|_{ L^{p(s)-\e}([0,T_\beta], H^s(\R))}\leq C,
%&\|\sqrt{\rho} v\|_{L^\infty((0,t),L^2(\R))}\leq C\\
%&\|\rho v(t,\cdot)\|_{{\cal M}(\R)}\leq C\\
%&\|\rho-\bar{\rho}\|_{L^\infty((0,t),L^2(\R))}\leq C\\
%&\|\p_x (\rho^{\frac{1}{2}(\gamma+\alpha-1)})\|_{L^2((0,t),L^2(\R))}\leq C\\
%&\|\p_x\rho^{\alpha-\frac{1}{2}}\|_{L^2([0,T_\beta],L^2(\R))+L^\infty([T_\beta,+\infty[,L^2(\R))}\leq C.
\end{aligned}
\end{cases}
\label{sestim1vb}
\end{equation}
with $0<s<1$ and 
$p(s)=\frac{6}{1+2s}$ and $\e>0$ sufficiently small such that $p(s)-\e\geq 1$.
There exists $C>0$ sufficiently large such that for any $t>0$ and $C(t)>0$ depending on $t$ we have:
\begin{equation}
\begin{cases}
\begin{aligned}
&\|\rho(t,\cdot)\|_{L^\infty(\R)}\leq C,\\
&\|\rho u(t,\cdot)\|_{L^1(\R)}\leq C(t),\\
&\|\sqrt{\rho}u\|_{L^\infty((0,t),L^2(\R))}\leq C(1+\frac{1}{\sqrt{t}}),\\
%&\|u\|_{ L^{p(s)-\e}([0,T_\beta], H^s(\R))}\leq C\\
&\|m^{2}\|_{L^\infty((0,t),L^2(\R))}\leq C,\\
%&\|\rho v(t,\cdot)\|_{{\cal M}(\R)}\leq C\\
&\|\rho-\bar{\rho}\|_{L^\infty((0,t),L^2(\R))}\leq C,\\
&\|\p_x (\rho^{\frac{1}{2}(\gamma+\alpha-1)})\|_{L^2((0,t),L^2(\R))}\leq C,\\
&\|\p_x\rho^{\alpha-\frac{1}{2}}\|_{L^2([0,T_\beta],L^2(\R))+L^\infty([T_\beta,+\infty[,L^2(\R))}\leq C.
\end{aligned}
\end{cases}
\label{sestim1}
\end{equation}
If $0<\alpha<\frac{1}{2}$ we have in addition for $C>0$:
\begin{equation}
\begin{cases}
\begin{aligned}
&\| \frac{1}{\rho}\|_{L^\infty(\R^+,L^\infty(\R))}\leq C\\
&\|\p_x u\|_{L^2([T_\beta,+\infty[,L^2(\R))}\leq C.
\end{aligned}
\end{cases}
\label{sestim2}
\end{equation}
Furthermore we have for any $T>0$:
\begin{equation}
\begin{cases}
\begin{aligned}
&m\in C^w([0,T],{\cal M}(\R))\\
&m^1\in C^w([0,T],{\cal M}(\R)),
\end{aligned}
\end{cases}
\label{continue}
\end{equation}
with $m=\sqrt{\rho}\sqrt{\rho}u$ and $m^1=\p_x\va_1(\rho)+m$.
%When $\alpha=\gamma$ we have the same result without the assumption (\ref{petitesse}).
%Finally if $\mu(\rho)\geq\mu>0$ for all $\rho\geq 0$, if $\mu$ is uniformly Lipschitz then this solution is unique in the class of weak solutions satisfying the usual entropy inequality (\ref{3.7}).
\label{theo1}
\end{theorem}
\begin{remarka}
The main interest of this result is to observe the two following points, first we can deal with initial momentum $m_0$ which are vortex. This means that we can take $m_0$ in ${\cal M}(\R)$. The second thing is that $\va_1(\rho_0)$ is only in $BV(\R)$ and by composition theorem $\rho_0$ is also in $BV(\R)$.%to observe that since:
%$$\n\va_1(\rho_0)=\rho_0v_0-\rho_0 u_0,$$
%with $\va_1(x)+x\va(x)$, it implies that since $(\rho_0u_0,\rho_0v_0)$ are in ${\cal M}(\R)$ and $\rho_0$ is on $L^\infty(\R)$
%that $\va_1(\rho_0)$ is in $BV(\R)$. 
This implies obviously that $\rho_0$ can admit initially some shocks.\\
We can now observe that these shocks at the initial time $t=0$ are instantaneously regularizing inasmuch as we verify that $\rho(t,\cdot)\in C(\R)$ for any $0<t\leq T_\beta$. Indeed we observe that:
$$\p_x \va_1(\rho(t,\cdot))=\rho u(t,\cdot)+\rho v(t,\cdot).$$
Since $\rho u(t,\cdot)$ and $\rho v(t,\cdot)$ are in $L^2(\R)$ for $t>0$ and $\va_1(\rho(t,\cdot))-\va_1(\bar{\rho})$ belongs to $L^2(\R)$ for $t\in[0,T_\beta]$, we deduce by Sobolev embedding that 
$\va_1(\rho)(t,\cdot)-\va_1(\bar{\rho})$ is continuous and then $\rho(t,\cdot)$ is also continuous for any $t>0$ (indeed we know that $(\rho,\frac{1}{\rho})$ belongs to $L^\infty([0,T_\beta],L^\infty(\R))$).\\
We can try to explain this phenomena,  $v$ verifies only a damped transport equation and has a priori no regularizing effect. In opposite the velocity $u$ admits regularizing effects since $u$ verifies a parabolic equation. Since $\p_x\va_1(\rho)=\rho v-\rho u$, taking $v$ sufficiently regular and using  regularizing effects on $u$ allows to prove that $\p_x\va_1(\rho)$ is instantaneously bounded at least in $L^1(\R)$. We can precise what we mean by sufficiently regular for $v$.
Assume that $\rho_0v_0$ is at least in $L^1(\R)$ and that $\rho_0u_0$ belongs to ${\cal M}(\R)$ we should observe regularizing effects on the density $\rho$. Indeed  we can expect that $\rho v(t,\cdot)$ remains in $L^1(\R)$ but that $\rho u(t,\cdot)$ is in $L^1(\R)$ for $t>0$ (this regularizing effect should provide from the parabolic behavior on the momentum $m$). Heuristically it implies that $\p_x\va_1(\rho(t,\cdot))\in L^1(\R)$ for $t\in[0,T_\beta]$  and we can now use the fact that $W^{1,1}(\R)$ is embedded in $C^0(\R)$.\\
In other word the regularizing effects on the density $\rho$ depend on the regularity of the coupling between the density $\rho_0$ and the velocity $u_0$ which is expressed by the unknown $\rho_0 v_0$.
\end{remarka}
%\begin{remarka}
%The condition $0<c\leq \rho_0\leq C<+\infty$ is in fact deduced from (\ref{petitesse}) and the fact that $\rho_0-\bar{\rho}$ is in $L^2(\R)\cap BV(\R)$.
%\end{remarka}
\begin{remarka}
From (\ref{initial1}), we give a sense to the initial momentum $m_0$. It is however not clear if we can define properly $u_0$. It would be natural now to write $u_0$ as $u_0=\frac{1}{\rho_0} m_0$ with $m_0\in{\cal M}(\R)$. However since $\frac{1}{\rho_0}$ is a priori not in $C(\R)$ (the space of bounded continuous function), the product $\frac{1}{\rho_0}m_0$ is not well defined. It is the same if we assume in addition that $(\rho_0-\bar{\rho})\in \dot{B}^1_{1,\infty}(\R)$ and $m_0\in \dot{B}^0_{1,\infty}(\R)$ (we refer to \cite{Danchin} for the definition of the homogeneous Besov space and on the notion of paraproduct law). In this case $\frac{1}{\rho_0}m_0$ can not be defined using classical paraproduct law, in fact the remainder $R(m_0,\frac{1}{\rho_0}-\frac{1}{\bar{\rho}})$ is not well defined.\\
However if we consider $u_0=v_0+\p_x\va(\rho_0)$ this unknown is well defined since $v_0$ belongs to $L^2(\R)$ and $\p_x\va(\rho_0)$ is in ${\cal M}(\R)$. Indeed since $\va_1(\rho_0)$ is in $BV(\R)$ and since we have $0<c\leq\rho_0\leq C<+\infty$, we deduce using the definition of the total variation and the mean value theorem that $\va(\rho_0)$ is in $BV(\R)$. In our case $\rho_0$ and $m_0$ have to be considered as the initial data. It is also remarkable that instantaneously $u$ is well defined since on $[0,T_\beta]$, $u(t,\cdot)$ is in $L^1(\R)$.
\end{remarka}
\begin{remarka}
Numerous works dedicated to the incompressible and compressible Navier-Stokes equations are dedicated to the existence of global strong solution with small initial data  $(\rho_{0},u_{0})$ belonging to \textit{critical} space. By critical, we mean that the system (\ref{1}) is solved in functional spaces with norm
 invariant by the changes of scales which leave (\ref{1}) invariant.
In the case of barotropic fluids, we can observe that the transformations:
\begin{equation}
(\rho(t,x),u(t,x))\longrightarrow (\rho(l^{2}t,lx),lu(l^{2}t,lx)),\;\;\;l\in\R,
\label{12B}
\end{equation}
have that property, provided that the pressure term has been changed accordingly. Roughly speaking we expect that such spaces are optimal in term of regularity on the initial data in order to prove the well-posedness of the system  (\ref{1}). We can easily observe in particular that $\dot{B}^{\frac{1}{p}}_{p,r}(\R)\times B^{\frac{1}{p}-1}_{p,r}(\R)$ with $p\in[1,+\infty]$, $r\in[1,+\infty]$ is a good candidate since its norm remains invariant by the transformation (\ref{12B}).\\
In this sense the previous result is a Theorem dealing with critical initial data since $BV(\R)\times{\cal M}(\R)$ is critical for (\ref{12B}) if we consider initial data $(\rho_0,m_0)$. In our case, we need to impose the smallness assumption (\ref{petitesse}), it would be interesting to get a similar result without smallness assumption even if we must work with initial data slightly more regular.\\
It is natural then to try to adapt the results of \cite{Danchin,armb,JDE2} to the case of the dimension $N=1$,
in other words it consists to prove the existence of strong solution in finite time with initial data verifying:
 $$0<c\leq \rho_0\leq M<+\infty, \,\rho_0-\bar{\rho}\in \dot{B}^1_{1,1}(\R),\,u_0\in \dot{B}^0_{1,1}(\R).$$
 Unfortunately it seems delicate to prove this type of results, indeed if we apply the method develop in \cite{Danchin,armb,JDE2} it is necessary to give a sense to the term $(\mu(\rho)-\mu(\bar{\rho}))\p_{xx}u$
 with $u\in \widetilde{L}^1([0,T],\dot{B}^2_{1,1}(\R))$, $(\rho-\bar{\rho})\in \widetilde{L}^\infty([0,T],\dot{B}^1_{1,1}(\R))$ and $\rho\geq\frac{c}{2}$ on $[0,T]$ with $T<+\infty$. It implies that
 $\p_{xx}u$ is in $\widetilde{L}^1([0,T],\dot{B}^0_{1,1}(\R))$. Unfortunately the product $(\mu(\rho)-\mu(\bar{\rho}))\p_{xx}u$
is a priori not defined, indeed if we apply the classical paraproduct law, we can observe that
$R((\mu(\rho)-\mu(\bar{\rho})),\p_{xx}u)$ is not defined in dimension $N=1$. This term is critical in the sense that we have $s_1+s_2+N\inf(0,\frac{1}{p}-\frac{1}{p_1}-\frac{1}{p_2})=0$ with $s_1=0$, $s_2=1$, $N=1$, $p=1$, $p_1=1$ and $p_2=1$.\\
It is then not obvious to prove the existence of strong solution in critical Besov spaces in one dimension, it would be however quite easy if we assume the initial data slightly subcritical $(\rho_0-\bar{\rho},u_0)\in (\dot{B}^1_{1,1}(\R)\cap \dot{B}^{1+\e}_{1,1}(\R))\times( \dot{B}^0_{1,1}(\R)\cap \dot{B}^\e_{1,1}(\R))$ with $\e>0$.\\
In this sense our result gives an other method than  \cite{Danchin,armb,JDE2}  to deal with critical initial data in one dimension.%deals with critical initial data even if we have no result of uniqueness.
\end{remarka}
%\begin{remarka}
%Expliquer la notion de $BV$
%\end{remarka}
\begin{remarka}
This Theorem is a result of global weak solution with infinite energy space, this is due to the fact that $\sqrt{\rho_0}u_0$ does not belong to $L^2(\R)$. To the best of our knowledge, this is the first result of this type for compressible Navier Stokes equations.
\end{remarka}
%\begin{remarka}
%Behavior in long time of $\frac{1}{\rho}$ in $L^\infty$ norm\\
%A CE NIVEAU LA ON A DEUX POSSIBILITES, UTILISER LE THEOREME DE JIU POUR LES TEMPS LONGS, POUR CE FAIRE IL FAUT PROUVER UN GAIN MELLET VASSEUR UNIFORME.\\ CARREMENT MONTRER QUE L'ON DEVIENT FORT INSTANTANEMENT!
%EN FAIT ON PEUT IMAGINER QUE L'ON A UNE SOLUTION FORTE EN REALITE!
%\end{remarka}
%FAIRE LE CAS $v_0\in L^1(\R)$!!!!\\
We obtain now the following two corollaries. 
\begin{cor}
Assume that $\mu(\rho)=\mu\rho^\alpha$ with $\mu,\alpha> 0$
and  $P(\rho)=a\rho^{\gamma}$ with $\gamma>0$, $a>1$ and $\gamma\geq \alpha$.
The initial data satisfy for $\bar{\rho}>0$:
%\begin{equation}
%\begin{cases}
%\begin{aligned}
%&m_0\in {\cal M}(\R),\\
%&0<c\leq\rho_0(x)\leq C<+\infty\;\;\;\mbox{for any}\;x\in\R\\
%&\rho_0-\bar{\rho}\in L^{\gamma}_{2}(\R),\;\p_x \va_1(\rho_0)\in {\cal M}(\R)\\
%&m_0^1\in L^2(\R).
%\end{aligned}
%\end{cases}
%\label{initial}
%\end{equation}
%
\begin{equation}
\begin{cases}
\begin{aligned}
&\va_1(\rho_0)\in BV(\R),\;(\rho_0-\bar{\rho})\in L^{\gamma}_{2}(\R)\\
&0<c\leq\rho_0(x)\leq C<+\infty\;\;\;\mbox{for any}\;x\in\R\\
%&\rho_0 v_0\in {\cal M}(\R),\,\rho_0 u_0=m_0\in {\cal M}(\R),\;v_0\in L^2(\R)\\
%&0<c\leq\rho_0(x)\leq C<+\infty\;\;\;\mbox{for any}\;x\in\R\\
&u_0\in L^2(\R), \;m_0\in{\cal M}(\R).
%&\rho_0-\bar{\rho}\in L^{\gamma}_{2}(\R).
\end{aligned}
\end{cases}
\label{initial}
\end{equation}
The momentum $m^1_0$ is defined as follows:
\begin{equation}
m^1_0=\rho_0u_0+\p_x\va_1(\rho_0),\,m^1_0\in{\cal M}(\R).
\label{initial1b}
\end{equation}
In addition theres exists $\e_0>0$ such that if:
%\begin{equation}
%\|\rho_0v_0\|_{{\cal M}(\R)}+\|\rho_0u_0\|_{{\cal M}(\R)}\leq \e_0,
%\label{petitesse}
%\end{equation}
\begin{equation}
\|\p_x\va_1(\rho_0)\|_{{\cal M}(\R)}+\|m_0\|_{{\cal M}(\R)}\leq \e_0,
\label{petitesse1}
\end{equation}
then there exists $T>0$ and a weak solution $(\rho,u)$ for the system (\ref{1}) verifying the definition \ref{def} on the time interval $[0,T]$. We have in addition:
\begin{equation}
\begin{cases}
\begin{aligned}
&u\in L^\infty([0,T],L^2(\R)),\;\p_x u\in L^2([0,T],L^2(\R))\\
&(\rho,\frac{1}{\rho})\in L^\infty([0,T],L^\infty(\R))^2,\; \p_x\rho\in L^\infty([0,T],{\cal M}(\R))\\
&(\rho-\bar{\rho})\in L^\infty([0,T],L^\gamma_2(\R)),\;\rho\in C([0,T],L^p_{loc}(\R))\\
&\rho u\in C([0,T],{\cal M}(\R)^*)\\
&m^1\in C([0,T],{\cal M}(\R)^*).
%&\frac{1}{2}\int_{\R}\big(\rho(T,x)|u(T,x)|^2+(\Pi(\rho(T,x))-\Pi(\bar{\rho})) \big) dx\\
%&+\int^T_0\int_{\R}\mu_n(\rho(s,x))(\p_x u(s,x))^2 dx\leq \frac{1}{2}\int_{\R}\big(\rho_n(0,x)|u_n(0,x)|^2+(\Pi(\rho_n(0,x))-\Pi(\bar{\rho})) \big) dx.
\end{aligned}
\end{cases}
\label{prec}
\end{equation}
%COMPLETER LES ESTIMATIONS!\\
%In addition when $\gamma=\alpha=1$ and $u_0\in L^\infty(\R)$, there exists a global weak $(\rho,u)$ for the system (\ref{1}) verifying the estimates (\ref{prec}). In this case we do not need the smallness assumption (\ref{petitesse1})
%Finally if $\mu(\rho)\geq\mu>0$ for all $\rho\geq 0$, if $\mu$ is uniformly Lipschitz then this solution is unique in the class of weak solutions satisfying the usual entropy inequality (\ref{3.7}).
\label{corbis}
\end{cor}
\begin{remarka}
The previous corollary is a generalization of \cite{Hoff86} for general viscosity coefficients (in \cite {Hoff86} the author deals with constant viscosity coefficients). In the previous theorem $u_0$ is well defined and we can deal with initial density admitting shocks since $\rho_0$ is in $BV(\R)$.
\end{remarka}
%\begin{cor}
%Assume that $\mu(\rho)=\mu\rho^\alpha$ with $\mu,\alpha> 0$
%and  $P(\rho)=a\rho^{\gamma}$ with $\gamma>1$, $a>1$ and $\gamma\geq \alpha$.
%The initial data satisfy for $\bar{\rho}>0$:
%\begin{equation}
%\begin{aligned}
%&(\rho_0-\bar{\rho})\in B^1_{\infty,1}(\R),\,u_0\in B^0_{\infty,1}(\R),\;v_0\in L^2(\R)\\
%&0<c\leq\rho_0(x)\leq C<+\infty\;\;\;\mbox{for any}\;x\in\R\\
%&\rho_0-\bar{\rho}\in L^{\gamma}_{2}(\R).
%\end{aligned}
%\end{equation}
%Then there exists a global weak solution for the system (\ref{1}) verifying the definition \ref{def}.
%Furthermore we have the same estimates (\ref{sestim1}) and (\ref{sestim2}) then in the previous theorem.
%\begin{equation}
%\begin{aligned}
%&\|(\rho(t,\cdot),\frac{1}{\rho}(t,\cdot))\|_{L^\infty(\R)}\leq C,\\
%&\|\rho u(t,\cdot)\|_{L^1(\R)}\leq C,\\
%&\|u\|_{ L^{p(s)-\e}([0,T_\beta], H^s(\R))}\leq C\\
%&\|v(t,\cdot)\|_{L^2(\R)}\leq C\\
%&\|\rho v(t,\cdot)\|_{L^1(\R)}\leq C,\\
%&\|\rho(t,\cdot)-\bar{\rho}\|_{L^2(\R)}\leq C
%\end{aligned}
%\end{equation}
%with $0<s<1$ and $p(s)=$
%METTRE DES ESTIMEES
%Furthermore we have for any $T>0$:
%\begin{equation}
%\begin{aligned}
%&\rho u\in C([0,T],L^1(\R))\\
%&\rho v\in C([0,T],L^1(\R)).
%\end{aligned}
%\label{continue}
%\end{equation}
%\label{cor1}
%\end{cor}
%\begin{remarka}
%Compared with the previous theorem, we have no smallness assumption (\ref{petitesse}) on the initial data. In counterpart we need to assume that $\rho_0u_0$ and $\rho_0v_0$ are more regular since they are in $L^1(\R)$ and not in ${\cal M}(\R)$.
%\end{remarka}
Let us deal now with the case of constant viscosity coefficients, we obtain the following Theorem.
\begin{theorem}
Assume that $\mu(\rho)=\mu$ with $\mu> 0$
and  $P(\rho)=a\rho^{\gamma}$ with $\gamma>1$, $a>0$.
The initial data satisfy for $\bar{\rho}>0$:
\begin{equation}
\begin{cases}
\begin{aligned}
&m^1_0\in {\cal M}(\R),\,\va_1(\rho_0)\in BV(\R),v_0\in L^2(\R)\\
&0<c\leq\rho_0(x)\leq C<+\infty\;\;\;\mbox{for any}\;x\in\R\\
&\rho_0-\bar{\rho}\in L^{\gamma}_{2}(\R).
\end{aligned}
\end{cases}
\end{equation}
The momentum $m_0$ is defined as follows:
\begin{equation}
m_0=\rho_0v_0-\p_x\va_1(\rho_0),\,m_0\in{\cal M}(\R).
\label{initial1vn}
\end{equation}
Then there exists a global weak solution $(\rho,u)$ for the system (\ref{1}) verifying the definition \ref{def}, in addition for any $t>0$ we have:
\begin{equation}
\begin{cases}
\begin{aligned}
&\|\rho u(t,\cdot)\|_{L^1(\R)}\leq C_1(t)\\
%&\|\rho v(t,\cdot))\|_{{\cal M}(\R)}\leq C_1(t)\\
&\|(\rho(t,\cdot),\frac{1}{\rho}(t,\cdot)\|_{L^\infty}\leq C_1(t)\\
&\| \sqrt{\rho}u(t,\cdot)\|_{L^2(\R)}\leq C(\frac{1}{\sqrt{t}}+1)\\
%&\|(\sqrt{s}1_{\{s\leq 1\}}+1)\p_x u\|_{L^2_t(L^2(\R))}\leq C\\
%&\|\sqrt{\rho}v(t,\cdot)\|_{L^2}\leq C\\
&\|\rho(t,\cdot)-\bar{\rho}\|_{L^\gamma_2(\R)}\leq C\\
&\|\p_x\rho\|_{L^2_t(L^2(\R))}\leq C C_1(t)\\
&\|m^2\|_{L^\infty(\R^+,L^2(\R))}\leq C.
%&\|\p_x(\rho^{\frac{1}{2}(\gamma-1)})\|_{L^2_t(L^2(\R))}\leq C\\
%&\|\p_x\rho(t,\cdot)\|_{L^2(\R)}\leq C C_1(t) (\frac{1}{\sqrt{t}}+1).
\end{aligned}
\end{cases}
\label{resumad2di}
\end{equation}
In addition $u$ belongs to $L^{p(s)-\e}_{loc}(H^s(\R))$ with $0\leq s<1$ and $p(s)=\frac{3}{1+2s}$ with $\e>0$ sufficiently small such that
$p(s)-\e\geq 1$. We have also $\rho\in C([0,+\infty[,W^{-\e,p}_{loc}(\R))$ and $\rho u\in C([0,+\infty[,{\cal M}(\R)^*)$.
%METTRE DES ESTIMEES
\label{theo2}
\end{theorem}
\begin{remarka}
Compared with the Theorem \ref{theo1}, we have no smallness assumption on the initial data.
\end{remarka}
\begin{remarka}
As in the Theorem \ref{theo1}, we can observe regularizing effects on the density. In \cite{Hof98}, the initial density verifies only $\rho_0-\bar{\rho}\in L^2(\R)\cap L^\infty(\R)$ and $\frac{1}{\rho_0}\in L^\infty(\R)$ (in particular the assumptions on the initial data are weaken as our Theorem even if they are stronger on the initial velocity since $u_0\in L^2(\R)$). In particular the Theorem of \cite{Hof98} allows to deal with initial shocks on the density and a priori there is no reasons to observe regularizing effects on the density. One of the reasons is probably that $\rho_0v_0$ is not sufficiently regular, indeed this term is a priori only in $W^{-1,\infty}(\R)$.
\end{remarka}
We recall now some classical lemma of compactness in particular the classical lemma of Aubin-Lions.
\begin{lemme}
\label{Aubin}
Let $X_0, X, X_1$ Banach spaces. Assume that $X_0$ is compactly embedded in $X$ and $X$ is continuously embedded in $X_1$. Let $1\leq p,q\leq +\infty$. We set for $T>0$:
$$W_T=\{u\in L^p([0,T],X_0),\,\frac{d}{dt}u\in L^q([0,T],X_1)\}.$$
Then if:
\begin{itemize}
\item $p<+\infty$ then the embedding of $W$ into $L^p([0,T],X)$ is compact.
\item $p=+\infty$ and $q>1$, then the embedding of $W$ into $C([0,T],X)$ is compact.
\end{itemize}
\label{Aubin}
\end{lemme}
We refer to \cite{Serre} for the following Lemma which is a consequence of the Helly's Theorem.
\begin{lemme}
We consider a sequence $(\rho_n)_{n\in\mathbb{N}}$ of functions defined on $[0,T]\times\R$ with values in $\R$ satisfying the following hypotheses:
\begin{enumerate}
\item There exists $M>0$ such that $TV(\rho_n(t,\cdot))\leq M$ for all $n\in\mathbb{N}$.
\item $|\rho_n(t,\cdot)|\leq M$ for all $n\in\mathbb{N}$.
\item There exists a sequence $(\e_n)_{n\in\mathbb{N}}$ which converges to $0^+$, such that:
$$\int_{\R}|\rho_n(t,x)-\rho_n(s,x)| dx\leq \e_n+M|t-s|^\beta$$
for all $n\in\mathbb{N}$, all $s,t\in[0,T]$ and $\beta>0$.
\end{enumerate}
Then this sequence is relatively compact in $L^1_{loc}((0,T)\times\R)$. 
\label{compac}
\end{lemme}
\begin{lemme}
Let $K$ a compact subset of $\mathbb{R}$ and $(v_n)_{n\in\mathbb{N}}$ a sequel such that:
\begin{itemize}
\item $(v_n)_{n\in\mathbb{N}}$ is uniformly bounded in $L^{p+\alpha}(K)$ with $\alpha>0$ and $p\geq1$,
\item $(v_n)_{n\in\mathbb{N}}$ converges almost everywhere to $v$,
\end{itemize}
then $(v_n)_{n\in\mathbb{N}}$ strongly  converges to $v$ in $L^{p}(K)$ with $v\in L^{p+\alpha}(K)$.
\label{lemmeimp}
\end{lemme}
{\bf Proof}: First by the Fatou lemma $v$ is in $L^{p+\alpha}(K)$. Next we have for any $M>0$:
\begin{equation}
\int_{K}|v_n-v|^{p}dx\leq \int_{K\cap\{|v_n-v|\leq M\}}|v_n-v|^{p}dx+ \int_{K\cap\{|v_n-v|\geq M\}}|v_n-v|^{p}dx.
\label{lemme1}
\end{equation}
We are dealing with the second member of the right hand side, by H\"older inequality and Tchebychev lemma we have for a $C>0$:
\begin{equation}
\begin{aligned}
\int_{K\cap\{|v_n-v|\geq M\}}|v_n-v|^{p}dx&\leq (\int_{K}|v_n-v|^{p}\frac{|v_n-v|^{\alpha}}{M^{\alpha}}dx)%^{\frac{1}{1+\alpha}}(\{|v^{\e}-v|\geq M\})^{\frac{\alpha}{1+\alpha}},\\
&\leq \frac{C}{M^{\alpha}}.
\end{aligned}
\label{lemme2}
\end{equation}
In particular we have shown the strong convergence of $v^{\e}$ to $v$, indeed from the inequality (\ref{lemme1}) it suffices to use the Lebesgue theorem for the first term on the right hand side and the estimate (\ref{lemme2}) with $M$ going to $+\infty$. {\hfill $\Box$}\\
\\
Section \ref{section2} deals with the proof of  the theorem \ref{theo1}. In the section \ref{section3} we prove the corollaries \ref{corbis}. We show the theorem \ref{theo2} which consider constant viscosity coefficients in the section \ref{section4}.
\section{Proof of theorem \ref{theo1}}
\label{section2}
%\subsection{Solutions faibles 1D avec donnŽes discontinues}
We start by constructing a sequence $(\rho_n,u_n)_{n\in\mathbb{N}}$ of global strong solutions of the system (\ref{1}) provided that the initial data is sufficiently regular, in addition we change slightly the viscosity coefficient. More precisely $(\rho_n,u_n)_{n\in\mathbb{N}}$ are solutions of the following systems:
\begin{equation}
\begin{cases}
\begin{aligned}
&\p_t\rho_n+\p_x(\rho_n u_n)=0,\\
&\p_t(\rho_n u_n)+\p_x(\rho_n u_n^2)-\p_x(\mu_n (\rho_n)\p_x u_n)+\p_x P(\rho_n)=0,
\end{aligned}
\end{cases}
\label{1bis}
\end{equation}
with $\mu_n(\rho_n)=\frac{1}{n}\rho_n^\theta+\mu\rho_n^\alpha$ with $\theta\in (0,\frac{1}{2})$ (we will fix $\theta$ later). In addition we assume that:
$$
\begin{cases}
\begin{aligned}
&m_n(0,\cdot)=e^{\frac{1}{n}\p_{xx}}(m_0)\\
&\va_1(\rho_n)(0,\cdot) =e^{\frac{1}{n}\p_{xx}}(\va_1(\rho_0))\\
&m_n^1(0,\cdot)=m_n(0,\cdot)+\frac{\mu_n(\rho_n(0,\cdot))}{\rho_n(0,\cdot)}\p_x\rho_n(0,\cdot)=e^{\frac{1}{n}\p_{xx}}(m^1_0)+\frac{1}{\mu n}\rho_n^{\theta-\alpha}(0,\cdot) e^{\frac{1}{n}\p_{xx}}(\p_x \va_1(\rho_0)),%=e^{\frac{1}{n}\D}m^1_0,
\end{aligned}
\end{cases}
$$
with $m^1_0=m_0+\p_x \va_1(\rho_0)$, $\va_1(\rho)=\frac{\mu}{\alpha}\rho^{\alpha}$ and $\alpha>0$. We verify easily since $\rho_0-\bar{\rho}$ is in $L^\gamma_2(\R)$ and we have $0<c\leq\rho_0\leq C<+\infty$ that $\rho_0-\bar{\rho}\in L^2(\R)$. In addition by composition theorem we deduce that $\va_1(\rho_0)-\va_1(\bar{\rho})$ belongs to $L^2(\R)$. %Similarly we observe that there exist $C_1,c_1>0$ independent on $n$ such that $0<c_1\leq \va_1(\rho_n)(0,\cdot)\leq C_1$, it implies in particular that we have $0<c_2\leq\rho_n(0,\cdot)\leq C_2$ for $0<c_2<C_2<+\infty$ independent on $n$.\\
It implies from (\ref{initial}) that there exists $C_{n,s}$ depending on $n$ and $s> 0$ such that for a continuous function $C$ independent on $n$ and $C_1>0$ large enough independent on $n$ we have: %(VERIFIER LE CAS MESURE SIMILAIRE A GIGA):
\begin{equation}
\begin{cases}
\begin{aligned}
&\|m_n(0,\cdot)\|_{L^1(\R)}%\leq \|m_0\|_{{\cal M}(\R)}
\leq \|m^1_0\|_{{\cal M}(\R)}+\|\va_1(\rho_0)\|_{BV(\R)}\\
&\|\p_x \va_1(\rho_n)(0,\cdot)\|_{L^1(\R)}\leq \|\va_1(\rho_0)\|_{BV(\R)}\\%\|\p_x \va_1(\rho_0)\|_{{\cal M}(\R)}\\
&0<\min(\va_1(c),\va_1(C))\leq \va_1(\rho_n)(0,\cdot)\leq\max( \va_1(C),\va_1(c))<+\infty\\
&0<c\leq \rho_n(0,\cdot)\leq C<+\infty\\
&\|m^1_n(0,\cdot)\|_{L^1(\R)}%\leq \|m_0\|_{{\cal M}(\R)}
\leq \|m^1_0\|_{{\cal M}(\R)}+\frac{1}{n\mu}\|\rho_n^{\theta-\alpha}(0,\cdot)\|_{L^\infty}\|\va_1(\rho_0)\|_{BV(\R)}
\end{aligned}
\end{cases}
\label{initial1v}
\end{equation}
and similarly we get for $s_1>0$ large enough:
\begin{equation}
\begin{cases}
\begin{aligned}
&\|\va_1(\rho_n(0,\cdot))-\va_1(\bar{\rho})\|_{H^s (\R)}\leq C_{n,s}\|\va_1(\rho_0)-\va_1(\bar{\rho})\|_{L^2(\R)}\\
&\|\rho_n(0,\cdot)-\bar{\rho}\|_{H^s (\R)}\leq C_{n,s} C(\|\rho_0\|_{L^\infty(\R)},\|\frac{1}{\rho_0}\|_{L^\infty(\R)}) \|\rho_0-\bar{\rho}\|_{L^2(\R)}\\
&\|\rho_n(0,\cdot)-\bar{\rho}\|_{L^2(\R)}\leq  C(\|\rho_0\|_{L^\infty(\R)},\|\frac{1}{\rho_0}\|_{L^\infty(\R)})\|\rho_0-\bar{\rho}\|_{L^2(\R)}\\
&\|m^1_n(0,\cdot)\|_{L^2(\R)}\leq \|m^1_0\|_{L^2(\R)}(1+\frac{C_1}{n^{\frac{3}{4}}}\|\rho_n(0,\cdot)\|_{L^\infty(\R)}^{\theta-\alpha})\\
&\|m^1_n(0,\cdot)\|_{H^{s_1}(\R)}\leq C_{n,s_1}\|m^1_0\|_{L^2(\R)}\\
&\hspace{1cm}+\frac{1}{n}(1+C_{n,s_1} C(\|\rho_0\|_{L^\infty(\R)},\|\frac{1}{\rho_0}\|_{L^\infty(\R)})\|\rho_0-\bar{\rho}\|_{L^2(\R)}) C_{n,s_1+1}\|\rho_0-\bar{\rho}\|_{L^2(\R)} \\
&\|m_n(0,\cdot)\|_{H^{s_1}(\R)}\leq C_{n,s}\|m^1_0\|_{L^2(\R)}\\
&\hspace{1cm}+\frac{1}{n}(1+C_{n,s_1} C(\|\rho_0\|_{L^\infty(\R)},\|\frac{1}{\rho_0}\|_{L^\infty(\R)})\|\rho_0-\bar{\rho}\|_{L^2(\R)}) C_{n,s_1+1}\|\rho_0-\bar{\rho}\|_{L^2(\R)} .
\end{aligned}
\end{cases}
\label{initial1v3}
\end{equation}
We have previously used composition theorem to pass from the $H^s$ norm on $\va_1(\rho_n(0,\cdot))-\va_1(\bar{\rho})$ to the $H^s$ norm on $\rho_n(0,\cdot)-\bar{\rho}$ and to deal
with the term $\|m^1_n(0,\cdot)\|_{H^{s_1}(\R)}$. Concerning the term $\|m^1_n(0,\cdot)\|_{L^2(\R)}$, we have used the fact that for $C_1>0$ large enough we have:
$$
\begin{aligned}
&\frac{1}{n}\|\rho_n^{\theta-\alpha}(0,\cdot)e^{\frac{1}{n}\p_{xx}}\p_x \va_1(\rho_0)\|_{L^2(\R)}\leq \frac{1}{n}\|\rho_n^{\theta-\alpha}(0,\cdot)\|_{L^\infty(\R)}C_1 n^{\frac{1}{4}}\|e^{\frac{1}{2n}\p_{xx}}\p_x \va_1(\rho_0)\|_{L^1(\R)}.
\end{aligned}
$$
Let
us mention that we can now define $u_n(0,\cdot)$ and $v_n(0,\cdot)$, indeed if we write $u_n(0,\cdot)=\frac{1}{\rho_n(0,\cdot)}m_n(0,\cdot)$ and $v_n(0,\cdot)=\frac{1}{\rho_n(0,\cdot)}m^1_n(0,\cdot)$  , these terms have a sense since they belong to $L^1(\R)$.\\
In particular we deduce from (\ref{initial1v}), (\ref{initial1v3}) and using product in Sobolev space that we have for $C_{n,s}>0$ depending on $n$ and every $s$ with $s> 0$:
\begin{equation}
\begin{cases}
\begin{aligned}
&\|\frac{1}{\rho_n(0,\cdot)}-\frac{1}{\bar{\rho}}\|_{H^s(\R)}\leq C_{n,s}\\
&\|u_n(0,\cdot)\|_{H^s(\R)}\leq C_{n,s_1}\\
%&0<\min(\va_1(c),\va_1(C))\leq \va_1(\rho_n)(0,\cdot)\leq\max( \va_1(C),\va_1(c))<+\infty\\
%&0<c\leq \rho_n(0,\cdot)\leq C<+\infty\\
%&\|\va_1(\rho_n(0,\cdot))-\va_1(\bar{\rho})\|_{H^s (\R)}\leq C_{n,s}\|\va_1(\rho_0)-\va_1(\bar{\rho})\|_{L^2(\R)}\\
%&\|\rho_n(0,\cdot)-\bar{\rho}\|_{H^s (\R)}\leq C_{n,s} \|\rho_0-\bar{\rho}\|_{L^2(\R)}\\
%&\|\rho_n(0,\cdot)-\bar{\rho}\|_{L^2(\R)}\leq \|\rho_0-\bar{\rho}\|_{L^2(\R)}\\
%&\|m^1_n(0,\cdot)\|_{L^2(\R)}\leq \|m^1_0\|_{L^2(\R)}\\
%&\|m^1_n(0,\cdot)\|_{H^s(\R)}\leq C_{n,s}\|m^1_0\|_{L^2(\R)}\\
%&\|m_n(0,\cdot)\|_{H^s(\R)}\leq C_{n,s}\|m^1_0\|_{L^2(\R)}+C_{n,s+1}\|\va_1(\rho_0)-\va_1(\bar{\rho})\|_{L^2(\R)}.
\end{aligned}
\end{cases}
\label{initial2v}
\end{equation}
with $s_1\geq s$ sufficiently large. We have just used the theorem of product in Sobolev space.
From \cite{MV}, we know that there exists a global strong solution for the system (\ref{1bis}) since $(\rho_n(0,\cdot)-\bar{\rho},u_n(0,\cdot))$ belongs to $H^1(\R)\times H^1(\R)$ with $0<c\leq \rho_n(0,\cdot)\leq C<+\infty$. Indeed we observe that we have the condition  $\mu_n(\rho)\geq \frac{1}{n}\rho^\theta$ with $\theta\in [0,\frac{1}{2})$ which is the relevant condition for the existence of global strong solution in \cite{MV}.\\
We are now going to prove uniform estimates in $n$ on the sequence $(\rho^n,u^n)_{n\in\mathbb{N}}$, in a second time we will prove that $(\rho^n,u^n)_{n\in\mathbb{N}}$ converges up to a subsequence to a global weak solution solution $(\rho,u)$ of (\ref{1}) with $(\rho_0,u_0)$ verifying the condition of theorem \ref{theo1}. We can mention that the solution $(\rho_n,u_n)$ here is classical in the sense that  $(\rho_n,u_n)$ is in $C^{\infty}([0,+\infty)\times\R)$. This is due to the fact that the $H^s$ regularity is preserved all along the time.
%Obviously we assume in the sequel that $\rho_0^n u_0^n$ and $\rho_0^n v_0^n$ are uniformly bounded in $L^1(\R)$ in $n$ and verify:
%\begin{equation}
%\begin{aligned}
%&\|\rho_0^n u_0^n-\rho_0 u_0\|_{{\cal M}(\R)}\rightarrow_{n\rightarrow+\infty} 0\\
%&\|\rho_0^n v_0^n-\rho_0 v_0\|_{{\cal M}(\R)}\rightarrow_{n\rightarrow+\infty} 0\\
%\end{aligned}
%\end{equation}
From (\ref{petitesse}) and (\ref{initial1v}) we deduce that there exists $C>0$ such that:
\begin{equation}
\|\rho_n u_n(0,\cdot)\|_{L^1(\R)}+\|\rho_n v_n(0,\cdot)\|_{L^1(\R)}\leq C \e_0.
\label{ini1}
\end{equation}
\subsection{Uniform estimates on $(\rho_n,u_n,v_n)_{n\in\mathbb{N}}$}
\subsubsection*{Estimates of $\rho_n u_n$ and $\rho_n v_n$ in $L^\infty_T(L^1(\R))$ for any $T>0$}
%For simplicity of the notations, we forget in the sequel the subscript $n$ of $(\rho_n,u_n,v_n,\mu_n(\rho_n))$.
From  (\ref{11}), we recall that we have (with $v_n=u_n+\frac{\mu_n(\rho_n)}{\rho_n^2}\p_x\rho_n$):
\begin{equation}
\begin{cases}
\begin{aligned}
&\rho_n\p_t u_n+\rho_n u_n\p_x u_n-\p_x(\mu_n(\rho_n)\p_x u_n)+\frac{P'(\rho_n)\rho_n^2}{\mu_n(\rho)}v_n=\frac{P'(\rho_n)\rho_n^2}{\mu_n(\rho_n)}u_n,\\
&\rho_n \p_t v_n+\rho_n u_n\p_x v_n+\frac{P'(\rho_n)\rho_n^2}{\mu_n(\rho_n)}v_n=\frac{P'(\rho_n)\rho_n^2}{\mu_n(\rho_n)}u_n.
\end{aligned}
\end{cases}
\label{11a}
\end{equation}
We set $j_k(s)=\sqrt{s^2+\frac{1}{k}}-\sqrt{\frac{1}{k}}$ and we have $j'_k(s)=\frac{s}{\sqrt{s^2+\frac{1}{k}}}$, $j''_k(s)=\frac{1}{k(s^2+\frac{1}{k})^{\frac{3}{2}}}$.  We are going now to multiply the first and the second equation respectively by $j'_k(u_n)$ and $j'_k(v_n)$. Since the solution $(\rho_n,u_n,v_n)$ is classical and is in particular in $C^2([0,+\infty)\times\R)$ we have:
\begin{equation}
\begin{cases}
\begin{aligned}
&\rho_n \p_t j_k(v_n)+\rho_n u_n\p_x j_k (v_n)+\frac{P'(\rho_n)\rho^2}{\mu_n(\rho_n)}v_n  j_k'(v_n)=\frac{P'(\rho_n)\rho_n^2}{\mu_n(\rho_n)}u_n j'_k(v_n)\\
&\rho_n\p_t j_k(u_n)+\rho_n u_n\p_x j_k(u_n)-\p_x(\mu_n(\rho_n)\p_x j_k(u_n))+\mu_n(\rho_n) j_k''(u_n)|\p_x u_n|^2\\
&\hspace{6cm}=-\frac{P'(\rho_n)\rho_n^2}{\mu_n(\rho_n)}v_n j_k'(u_n)+\frac{P'(\rho_n)\rho_n^2}{\mu_n(\rho_n)}u_n j'_k(u_n).
\end{aligned}
\end{cases}
\label{11b}
\end{equation}
We integrate now on $(0,T)\times\R$ and we obtain using the fact that $\lim _{|x|\rightarrow+\infty} u_n(t,x)=\lim _{|x|\rightarrow+\infty} v_n(t,x)=0$ for $t\geq 0$ (this is due to the fact that $(u_n,v_n)_{n\in\mathbb{N}}$ are in Sobolev space $H^s$ with high regularity) :
$$
\begin{aligned}
&\int_{\R} \rho_n(T,x) j_k(v_n)(T,x) dx+\int^T_0\int_{\R}  \frac{P'(\rho_n)\rho_n^2}{\mu_n(\rho_n)}(s,x) \frac{v_n^2(s,x)}{\sqrt{v_n^2(s,x)+\frac{1}{k}}}dsdx\\
&=\int_{\R} \rho_n(0,x) j_k(v_n(0,x))dx+\int^T_0\int_{\R} \frac{P'(\rho_n)\rho_n^2}{\mu_n(\rho_n)} (s,x)  \frac{v_n(s,x) u_n(s,x)}{\sqrt{v_n^2(s,x)+\frac{1}{k}}}dsdx.
\end{aligned}
$$
$$
\begin{aligned}
&\int_{\R} \rho_n(T,x) j_k(u_n)(T,x) dx+\int^T_0\int_{\R} \mu_n(\rho_n) j_k''(u_n)|\p_x u_n|^2 (s,x) dsdx\\
&=\int_{\R} \rho_n(0,x) j_k(u_n(0,x)) dx-\int^T_0\int_{\R} \frac{P'(\rho_n)\rho_n^2}{\mu_n(\rho_n)} (s,x)  \frac{v_n(s,x) u_n(s,x)}{\sqrt{u_n^2(s,x)+\frac{1}{k}}}dsdx\\
&\hspace{3cm}+\int^T_0\int_{\R} \frac{P'(\rho_n)\rho_n^2}{\mu_n(\rho_n)} (s,x)  \frac{u_n^2(s,x)}{\sqrt{u_n^2(s,x)+\frac{1}{k}}}dsdx.
\end{aligned}
$$
We deduce that we have:
\begin{equation}
\begin{aligned}
&\int_{\R} \rho_n(T,x) j_k(v_n)(T,x) dx
+\int_{\R} \rho_n(T,x) j_k(u_n)(T,x) dx\leq\\
&\int_{\R} \rho_n(0,x) |u_n(0,x)| dx+\int_{\R} \rho_n(0,x) |v_n(0,x)| dx\\
&\hspace{4cm}+3\int^T_0\int_{\R} \frac{P'(\rho_n)\rho_n^2}{\mu_n(\rho_n)} (|v_n(s,x)|+|u_n(s,x)|) ds dx.
\end{aligned}
\label{ent1}
\end{equation}
We pass to the limit when $k\rightarrow+\infty$ and using the Fatou lemma we obtain:
\begin{equation}
\begin{aligned}
&\int_{\R} \rho_n(T,x)|v_n|(T,x) dx
+\int_{\R} \rho_n(T,x) |u_n|(T,x) dx\leq\\
&\int_{\R} \rho_n(0,x) |u_n(0,x)| dx+\int_{\R} \rho_n(0,x) |v_n(0,x)| dx\\
&\hspace{4cm}+3\int^T_0\int_{\R} \frac{P'(\rho_n)\rho_n^2}{\mu_n(\rho_n)} (|v_n(s,x)|+|u_n(s,x)|) ds dx.
\end{aligned}
\label{ent1}
\end{equation}
From the Gronwall lemma we deduce that for any $T>0$:
\begin{equation}
\begin{aligned}
&\|\rho_n v_n(T,\cdot)\|_{L^1(\R)}+\|\rho_n u_n(T,\cdot)\|_{L^1(\R)} \\
&\hspace{1cm}\leq(\|\rho_n v_n(0,\cdot)\|_{L^1(\R)}+\|\rho_n u_n(0,\cdot)\|_{L^1(\R)} ) e^{3\int^T_0\|\frac{P'(\rho_n)\rho_n }{\mu_n(\rho_n)}(s,\cdot)\|_{L^\infty}ds}.
\end{aligned}
\label{cru1}
\end{equation}
At this level it is important to point out that we have:
$$\frac{P'(\rho_n)\rho_n}{\mu(\rho_n)}\leq \frac{a\gamma}{\mu}\rho_n^{\gamma-\alpha}+\frac{a\gamma}{n}\rho_n^{\gamma-\theta},$$
and we are going to use the fact that $\gamma-\alpha,\gamma-\theta\geq 0$ in order to control the $L^\infty$ norm of $\frac{P'(\rho_n)\rho_n}{\mu(\rho_n)}$ in terms of the $L^\infty$ norm of the density $\rho$.
%We know that $(\rho u,\rho v)$ MONTRER belongs to $L^\infty_T(
%We are going now to multiply the first and the second equation repectively by $\sqrt{u^2+\frac{1}{k}}
%On veut faire du faible avec seulement la BD entropie. On peut alors avoir des estimations $BV$ de $ u$ et $ v$. En particulier la vitesse initiale sera une mesure. En effet on a avec $f(x)=|x|$
%$$
%\begin{aligned}
%&\rho \p_t f(v)+\rho u\p_x f(v)+C \rho^\gamma f(v)=C\rho^\gamma u f'(v)\\
%&\rho\p_t f(u)+\rho u\p_x f(u)-\mu\p_x(\rho\p_x f(u))+\mu\rho f''(u)|\p_x u|^2=-C \rho^\gamma v f'(u)+C\rho^\gamma f(u).
%\end{aligned}
%$$
%En integrant il vient donc ce que l'on veut en appliquant un Gonwall si on controle la norme $L^\infty$ de $\rho^\gamma$ (il faut probablement utiliser de la petitesse).\\
\subsubsection*{BD Entropy}
Multiplying the momentum equation (\ref{11}) by $v$ and integrating over $(0,T)\times\R$, we obtain the following  entropy:
\begin{equation}
\begin{aligned}
&\frac{1}{2}\int_{\R}(\rho |v|^2(T,x)+(\Pi(\rho)-\Pi(\bar{\rho}))(T,x) dx+\int^T_0\int_{\R}\frac{P'(\rho)\mu(\rho)}{\rho^2}(s,x)|\p_x\rho|^2 (s,x) ds dx\\
&\hspace{6cm}\leq \frac{1}{2}\int_{\R}(\rho_0 |v_0|^2(x)+(\Pi(\rho_0)-\Pi(\bar{\rho}))(x) dx,
\end{aligned}
\label{BD}
\end{equation}
with:
$$\Pi(s)=s(\int^s_{\bar{\rho}}\frac{P(z)}{z^2}dz-\frac{P(\bar{\rho})}{\bar{\rho}}).$$
It gives in our case for the the system (\ref{1bis}) when $\gamma+\alpha-3\ne 0$ with $\mu_n(\rho_n)=\mu\rho_n^\alpha+\frac{1}{n}\rho_n^\theta$ and $C>0$ large enough independent on $n$:
\begin{equation}
\begin{aligned}
&\frac{1}{2}\int_{\R}(\rho_n |v_n|^2(T,x)+(\Pi(\rho_n)-\Pi(\bar{\rho}))(T,x) dx+\frac{4a\gamma\mu}{(\gamma+\alpha-1)^2} \int^T_0\int_{\R}|\p_x\rho_n^{\frac{1}{2}(\gamma+\alpha-1)}|^2 (s,x) ds dx\\
&+\frac{4\gamma}{n(\gamma+\theta-1)^2} \int^T_0\int_{\R}|\p_x\rho_n^{\frac{1}{2}(\gamma+\theta-1)}|^2 (s,x) ds dx\\
&\hspace{4cm}\leq \frac{1}{2}\int_{\R}(\rho_n(0,x) |v_n(0,x)|^2+(\Pi(\rho_n(0,x))-\Pi(\bar{\rho})) dx\leq C,
\end{aligned}
\label{BD1}
\end{equation}
and we know that $\gamma+\alpha-1>0, \gamma+\theta-1>0$ since $\alpha,\theta\geq 0$ and $\gamma>1$. The fact that $C>0$ is independent on $n$ is a direct consequence from (\ref{initial1v}) and (\ref{initial1v3}). It is also important to note that by definition we have $v_n=u_n+\p_x\va_n(\rho_n)$ with $\va'_n(\rho_n)=\frac{\mu_n(\rho_n)}{\rho_n^2}$.
\subsubsection*{Uniform $L^\infty$ estimate on the density $(\rho_n)_{n\in\mathbb{N}}$ in finite time } % with $\va_1'(\rho)=\frac{\mu(\rho)}{\rho}$ (it gives $\va_1(\rho)=\frac{\mu}{\alpha}\rho^{\alpha}$ since $\alpha> 0$)
\label{norme}
%ESTIMATION EN $n$!\\
From (\ref{cru1}) we deduce that $\p_x\va_{1,n}(\rho_n)=\rho_n v_n-\rho_n u_n$ (with $\va_{1,n}(\rho_n)=\frac{\mu}{\alpha}\rho_n^\alpha+\frac{1}{n\theta}\rho_n^\theta$) is bounded in $L^\infty_T(L^1(\R))$ for any $T>0$ provided that $\rho_n$ belongs to $L^\infty_T(L^\infty(\R))$. Since the density $\rho_n-\bar{\rho}$ and $m_n$ are bounded in $C([0,T], H^s(\R))$ for any $s>0$, we deduce %from the mass equation that $\rho$  is uniformly continuous. Furthermore from (\ref{BD}) we deduce
 that for any $t\geq 0$ we have:
\begin{equation}
\lim_{|x|\rightarrow+\infty}\va_1(\rho_n)(t,x)=\va_1(\bar{\rho}).
\end{equation}
Now using the fact that the space $W^{1,1}(\R)$ is embedded in $L^\infty(\R)$ we deduce from (\ref{cru1}) that for any $t>0$:
\begin{equation}
\begin{aligned}
\|\va_{1,n}(\rho_n(t,\cdot))\|_{L^\infty}&\leq \va_{1,n}(\bar{\rho})+\|\p_x\va_{1,n}(\rho_n (t,\cdot))\|_{L^1(\R)}\\
&\leq \va_{1,n}(\bar{\rho})+(\|\rho_n v_n(0,\cdot)\|_{L^1(\R)}+\|\rho_n u_n(0,\cdot)\|_{L^1(\R)} ) e^{3\int^T_0\|\frac{P'(\rho_n)\rho_n }{\mu_n(\rho_n)}(s,\cdot)\|_{L^\infty}ds}.
%(\|\rho_0 v_0\|_{L^1(\R)}+\|\rho_0 u_0\|_{L^1(\R)} ) e^{3\int^t_0\|\frac{P'(\rho)\rho}{\mu(\rho)}(s,\cdot)\|_{L^\infty}ds}.
\end{aligned}
\label{Linf}
\end{equation}
%UTILISER LE DAMPING POUR AVOIR VIA LA PETITESSE UNE ESTIMATION GLOBALE!\\
The previous inequality shows that we can prove $L^\infty$ estimates on the density $\rho_n$ in finite time by using bootstrap arguments. As previously, we assume that $\gamma-\alpha\geq 0$ and we have using (\ref{initial1v}) for $C>0$, $\alpha>0$ and any $t>0$:
\begin{equation}
\begin{aligned}
&\|\rho_n\|^\alpha_{L_t^\infty(L^\infty(\R))}
\leq C(\bar{\rho}^\alpha+\frac{1}{n}\bar{\rho}^\theta)\\
&+C(\|\rho_0 v_0\|_{{\cal M}(\R)}+\|\rho_0 u_0\|_{{\cal M}(\R)})(1+\frac{1}{n}\|\rho_0\|_{L^\infty(\R)}^{\theta-\alpha}) e^{3 C t( \|\rho_n\|^{\gamma-\alpha}_{L_t^\infty(L^\infty(\R))}+\frac{1}{n} \|\rho_n\|^{\gamma-\theta}_{L_t^\infty(L^\infty(\R))})}.
\end{aligned}
\label{Linfa}
\end{equation}
And it yields that for $C'>0$ independent on $n$, $n$ large enough:
\begin{equation}
\begin{aligned}
\|\rho_n\|_{L_t^\infty(L^\infty(\R))}
&\leq C'\bar{\rho}+C'(\|\rho_0 v_0\|_{{\cal M}(\R)}+\|\rho_0 u_0\|_{{\cal M}(\R)} )^{\frac{1}{\alpha}} \,e^{\frac{3 C'}{\alpha}t( \|\rho\|^{\gamma-\alpha}_{L_t^\infty(L^\infty(\R))}+\frac{1}{n} \|\rho_n\|^{\gamma-\theta}_{L_t^\infty(L^\infty(\R)})}.
\end{aligned}
\label{Linf1a}
\end{equation}
Let us prove now that the sequence $(\rho_n)_{n\in\mathbb{N}}$ is uniformly bounded in $n$ in $L^\infty$ norm on an time interval $[0,T^*]$ with $T^*>0$ independent on $n$. More precisely we define by:
$$T_n=\sup\{ t\in(0,+\infty),\|\rho_n(t,\cdot)\|_{L^\infty}\leq \sup (2\|\rho_0\|_{L^\infty},M_1)\}% 2 C'\bar{\rho}+2C'(2\|\rho_0 v_0\|_{{\cal M}(\R)}+\|\p_x\va_1(\rho_0)\|_{{\cal M}(\R)} )^{\frac{1}{\alpha}})\},
$$
with $M_1=2 C'\bar{\rho}+2C'(2\|\rho_0 v_0\|_{{\cal M}(\R)}+\|\p_x\va_1(\rho_0)\|_{{\cal M}(\R)}$ and $M= \sup (2\|\rho_0\|_{L^\infty}, M_1 )$.\\
We observe that $T_n>0$ since $\rho_n$ belongs to $C([0,+\infty[, L^\infty(\R))$ for any $n\in\mathbb{N}$ and since $\|\rho_n(0,\cdot)\|_{L^\infty(\R)}<M$. Let us define now $T_*$ such that:
$$e^{\frac{3C'}{\alpha} T_* (M^{\gamma-\alpha}+M^{\theta-\alpha})}=\frac{3}{2}\;\;\;\mbox{and}\;\;\; T_*=\frac{\alpha \ln (\frac{3}{2})}{3 C'(M^{\gamma-\alpha}+M^{\theta-\alpha})}.$$
From the definition of $T_n$ and from (\ref{Linf1a}) we deduce that for any $n\in\mathbb{N}$ large enough we have:
$$T_n\geq T_*>0.$$
It implies that $(\rho_n)_{n\in\mathbb{N}}$ is uniformly bounded in $n$ in $L^\infty((0,T_*),L^\infty(\R))$ with $T_*$ independent on $n$. In other words we have for any $n\in\mathbb{N}^*$:
\begin{equation}
\|\rho_n\|_{L^\infty([0,T_*],L^\infty(\R))}\leq M.
\label{linf}
\end{equation}
 In the sequel we will prove in fact that $(\rho_n)_{n\in\mathbb{N}}$ is also uniformly bounded in $n$ in $L^\infty((0,\infty),L^\infty(\R))$. Now if we combine (\ref{linf}), (\ref{initial1v}) and (\ref{cru1}) we deduce that there exists $C>0$ independent on $n$ such that:
\begin{equation}
\begin{aligned}
&\|\rho_n u_n\|_{L^\infty([0,T_*],L^1(\R))}+\|\rho_n v_n\|_{L^\infty([0,T_*],L^1(\R))}\\
&\hspace{3cm}\leq (2\|\rho_0 v_0\|_{{\cal M}(\R)}+\|\p_x\va_1(\rho_0)\|_{{\cal M}(\R)})e^{C T_* M^{\gamma-\alpha}}.
\end{aligned}
\label{L1}
\end{equation}
%\\
%The proof is similar when $\alpha=0$ except that we prove in one way that $\rho$ and $\frac{1}{\rho}$ are uniformly bounded in $n$ on a finite time $[0,T_\alpha]$. This is due to the fact that we can
%estimate the $L^\infty_x$ norm of $\ln\rho$.
\subsubsection*{Control of the $L^\infty$ norm of $(\frac{1}{\rho_n})_{n\in\mathbb{N}}$ in finite time}
\label{vide}
Now we are interested in estimating the $L^\infty$ norm of $\frac{1}{\rho_n}$ on a finite time $T_\beta$ independent on $n$. As in (\ref{Linf}) we have for $C>0$ large enough, for $n$ large enough and any $t\in [0,T_*]$:
%\begin{equation}
%\begin{aligned}
%\|\va_1(\rho(t,\cdot))-\va_1(\bar{\rho})\|_{L^\infty}&\leq \|\p_x\va_1(\rho(t,\cdot))\|_{L^1(\R)}\\
%&\leq (\|\rho_0 v_0\|_{L^1(\R)}+\|\rho_0 u_0\|_{L^1(\R)} ) e^{3 t \|\frac{P'(\rho)\rho}{\mu(\rho)}\|_{L^\infty_t(L^\infty)}}.
%\end{aligned}
%\label{Linf3}
%\end{equation}
%We recall now that we have for any $t\in[0,T_*]$:
%$$\|\rho\|_{L^\infty_t(L^\infty(\R)}\leq M.$$
%It implies that we have for any $t\in [0,T_*]$ 
\begin{equation}
\begin{aligned}
\|\va_{1,n}(\rho(t,\cdot))-\va_{1,n}(\bar{\rho})\|_{L^\infty}&\leq \|\p_x\va_{1,n}(\rho_n(t,\cdot))\|_{L^1(\R)}\\
&\leq 2 (\|\rho_0 v_0\|_{{\cal M}(\R)}+\|\rho_0 u_0\|_{{\cal M}(\R)} ) e^{t C M^{\gamma-\alpha}}.
\end{aligned}
\label{Linf4}
\end{equation}
Now since we have $\|\rho_0 v_0\|_{{\cal M}(\R)}+\|\rho_0 u_0\|_{{\cal M}(\R)} \leq C_2 \e_0$ with $C_2>0$ independent on $n$ and with $\e_0$ sufficiently small, we deduce that for $T_\beta\leq T_*$ such that:
\begin{equation}
e^{T_\beta C M^{\gamma-\alpha}}\leq 2,
\end{equation}
we obtain that for $C_3>0$ large enough and any $t\in [0,T_\beta]$ and any $x\in\R$:
\begin{equation}
\va_1(\rho_n(t,x))\geq \va_1(\bar{\rho})-4C_2 \e_0-\frac{C_3}{n}\|\rho_n\|_{L^\infty_t(L^\infty(\R))}^\theta.
\label{vide2}
\end{equation}
From (\ref{vide2}) and (\ref{linf}) it implies that $\frac{1}{\rho_n}$ belongs uniformly to $L^\infty([0,T_\beta],L^\infty(\R))$ for $n$ large enough with $T_\beta>0$ independent on $n$ if $\e_0>0$ is small enough.
\subsubsection*{Gain of regularity on the velocity $u_n$}
From (\ref{BD}), we deduce that for any $T>0$ we have uniformly in $n$:
\begin{equation}
\begin{aligned}
&\p_x\rho_n^{\frac{1}{2}(\gamma+\alpha-1)}\in L^2_T(L^2(\R)),\;\sqrt{\rho_n}v_n\in L^\infty_T(L^2(\R)).
\end{aligned}
\label{gainu}
\end{equation}
Now from the definition of the effective velocity $v_n$ we have:
\begin{equation}
\begin{aligned}
\sqrt{\rho_n}u_n&=\sqrt{\rho_n}v_n-\frac{\mu_n(\rho_n)}{\rho_n^{\frac{3}{2}}}\p_x\rho_n\\
&=\sqrt{\rho_n}v_n-\frac{2\mu}{\gamma+\alpha-1}\rho_n^{\frac{1}{2}(\alpha-\gamma)}\p_x\rho_n^{\frac{1}{2}(\gamma+\alpha-1)}-\frac{2}{n(\gamma+\theta-1)}\rho_n^{\frac{1}{2}(\theta-\gamma)}\p_x\rho_n^{\frac{1}{2}(\gamma+\theta-1)}.
%=\sqrt{\rho}v-\frac{2\mu}{\gamma+\alpha-1}\frac{\rho^\alpha}{\rho^{\frac{1}{2}(\gamma+\alpha)}}\p_x\rho^{\frac{1}{2}(\gamma+\alpha-1)}\\
%&=\sqrt{\rho}v-\frac{2\mu}{\gamma+\alpha-1}\rho^{\frac{1}{2}(\alpha-\gamma)}\p_x\rho^{\frac{1}{2}(\gamma+\alpha-1)}
\end{aligned}
\label{dgainu}
\end{equation}
From (\ref{gainu}), (\ref{vide2}) and since $\alpha-\gamma,\theta-\gamma\leq 0$, we deduce that  $(\sqrt{\rho_n}u_n)_{n\in\mathbb{N}}$ is uniformly bounded in $L^2((0,T_\beta)\times\R)$, it gives that for $T_\beta>0$ independent on $n$, it exists $C>0$ independent on $n$ such that:%depending on $\e_0$ and the initial data $(\rho_0,v_0)$ that for any $n\in\mathbb{N}^*$:
\begin{equation}
\|\sqrt{\rho_n}u_n\|_{ L^2((0,T_\beta)\times\R)}\leq C.
\label{again}
\end{equation}
Let us prove now some additional regularizing effects independent on $n$ on the sequel $(u_n)_{n\in\mathbb{N}}$. It will be important in order to pass to the limit when $n$ goes to $+\infty$. ALet us multiply the momentum equation of (\ref{1bis}) by $s u_n$ and integrate over $(0,t)\times\R$ with $t\in(0,T_\beta)$, we have then:
%$$
%\begin{aligned}
%&\frac{1}{2}\p_t\int \rho t^\alpha \p_t(|u|^2) dx+\frac{1}{2}\int\rho\p_x (u^2) t^\alpha u dx+\mu\int\rho (\p_x u)^2 t^\alpha+\int t^\alpha \p_x P(\rho)\, u dx=0.
%\end{aligned}
%$$
%Et on a:
$$
\begin{aligned}
%&\frac{1}{2}\int \rho t^\alpha \p_t(|u|^2) dx+\frac{1}{2}\int u^2\p_t \rho t^\alpha u dx+\mu\int\rho (\p_x u)^2 t^\alpha+\int t^\alpha \p_x P(\rho)\, u dx=0,\\
&\frac{1}{2}t \int_{\R}  \rho_n |u_n|^2(t,x)  dx-\frac{1}{2}\int^{t}_0\int_\R  \rho_n |u_n|^2 (s,x) ds dx+\int^{t}_0 \int_{\R} \mu_n(\rho_n) (\p_x u_n)^2(s,x) s\, ds dx\\
&+\int^{t}_0 \int_{\R} s  \p_x P(\rho_n)\, u_n dx=0
\end{aligned}
$$
We recall now that we have:
$$\p_t(t(\Pi(\rho_n)-\Pi(\bar{\rho}))+t\p_x (\Pi(\rho_n) u_n)+t P(\rho_n)\p_x u_n-( \Pi(\rho_n)-\Pi(\bar{\rho}))=0.$$
We deduce then that we have:
\begin{equation}
\begin{aligned}
%&\frac{1}{2}\int \rho t^\alpha \p_t(|u|^2) dx+\frac{1}{2}\int u^2\p_t \rho t^\alpha u dx+\mu\int\rho (\p_x u)^2 t^\alpha+\int t^\alpha \p_x P(\rho)\, u dx=0,\\
&\frac{1}{2}t \int_{\R}  \rho_n |u_n|^2(t,x)  dx+\int^{t}_0 \int_{\R} \mu_n(\rho_n) s (\p_x u_n)^2(s,x)  ds dx+t \int_{\R} (\Pi(\rho_n)(t,x)-\Pi(\bar{\rho}))  dx\\
&\hspace{1cm}=\int^{t}_0\int_\R(\Pi(\rho_n)-\Pi(\bar{\rho}))(s,x) ds dx+\frac{1}{2}\int^{t}_0\int_\R  \rho_n |u_n|^2 (s,x) ds dx.
\end{aligned}
\label{energu}
\end{equation}
We deduce from %(\ref{energu}), (\ref{vide2}), 
(\ref{BD}), (\ref{vide2}), (\ref{again}) and (\ref{energu}) that the sequence $(u_n)_{n\in\mathbb{N}}$ verify uniformly in $n$ on $[0,T_\beta]$:
\begin{equation}
\begin{aligned}
&\sqrt{t}u_n\in L^\infty([0,T_\beta],L^2(\R))\;\;\;\mbox{and}\;\;\;\sqrt{t}\p_x u_n\in L^2([0,T_\beta],L^2(\R)).
\end{aligned}
\label{cru2}
\end{equation}
We can obtain additional uniform informations on the sequence $(u_n)_{n\in\mathbb{N}}$ by using interpolation estimates. Indeed from (\ref{L1}) we know that $(\rho_n u_n)_{n\in\mathbb{N}}$ is uniformly bounded in $L^\infty([0,T_\beta],L^1(\R))$. We recall that $(\frac{1}{\rho_n})_{n\in\mathbb{N}}$ is uniformly bounded on $[0,T_\beta]$ then we deduce that $(u_n)_{n\in\mathbb{N}}$ is uniformly bounded in 
$L^\infty([0,T_\beta],L^1(\R))$  and then is uniformly bounded in $L^\infty([0,T_\beta],H^{-\frac{1}{2}}(\R))$  by Sobolev embedding. It is important to mention that the bound depend only on
$\e_0$ and $\|\rho_0\|_{L^\infty(\R)}$.\\
By interpolation, we have for $0<s<1$ and any $t\in (0,T_\beta]$:
\begin{equation}
\begin{aligned}
&t^{\frac{1}{6}(1+2s)} \|u_n(t,\cdot)\|_{H^s}\leq \|u_n(t,\cdot)\|^{\frac{2}{3}(1-s)}_{H^{-\frac{1}{2}}(\R)}(\sqrt{t}\|u_n(t,\cdot)\|_{H^{1}(\R)})^{\frac{1}{3}(1+2s)}.
\end{aligned}
\label{cru3}
\end{equation}
From (\ref{vide2}), (\ref{cru2}) and (\ref{cru3}) we deduce that $(u_n)_{n\in\mathbb{N}}$ is uniformly bounded in $L^{p(s)-\e}([0,T_\alpha], H^s(\R))$ with $s\in(0,1)$, $p(s)=\frac{3}{1+2s}$ and $\e>0$ sufficiently small such that $p(s)-\e\geq 1$ ($\e>0$ depends here on $s\in(0,1)$). We have finally obtained that for any $t\in[0,T_\beta]$, we have for $C>0$ independent on $n$:
\begin{equation}
\begin{cases}
\begin{aligned}
&\|\rho_n(t,\cdot)-\bar{\rho}\|_{L^2(\R)}\leq C,\;\|(\frac{1}{\rho_n},\rho_n)(t,\cdot)\|_{L^\infty(\R)}\leq C\\
&\|v_n(t,\cdot)\|_{L^2(\R)}\leq C,\;\|v_n(t,\cdot)\|_{L^1(\R)}\leq C\\
&\|\sqrt{s}\p_x u_n\|_{L^2([0,T_\beta],L^2(\R))}\leq C\\
&\|u_n(t,\cdot)\|_{L^2(\R)}\leq C(1+\frac{1}{\sqrt{t}}),\;\|u_n(t,\cdot)\|_{L^1(\R)}\leq C\\
&\|\p_x \va(\rho_n)\|_{L^2(\R)}\leq C(1+\frac{1}{\sqrt{t}})\\
&\|u_n\|_{L^{p(s)-\e}([0,t], H^s(\R))}\leq C.
%&\|(\frac{1}{\rho_n},\rho_n)(t,\cdot)\|_{L^\infty(\R)}\leq C.%\\
%&\|\p_x\rho^{\frac{1}{2}(\gamma+\alpha-1)}(t,\cdot)\|_{L^2(\R)}\leq C.
\end{aligned}
\end{cases}
\label{impbdu}
\end{equation}
Indeed we use the fact that $\|\p_x\va(\rho_n)\|_{L^2(\R)}\leq \|\p_x\va_n(\rho_n)\|_{L^2(\R)}$ with $\va_n'(\rho)=\frac{\mu_n(\rho_n)}{\rho_n^2}$.
We deduce in particular that choosing $s_1(\e)\in(0,1)$ such that $\frac{1}{2}-s_1(\e)=\frac{1}{p(s_1(\e))-\e}$, we have  from (\ref{impbdu}) and by Sobolev embedding that:
\begin{equation}
\|u_n\|_{L^{p(s_1(\e))-\e}([0,T_\beta]\times\R)}\leq C,
\label{impbdu1}
\end{equation}
with $C>0$ independent on $n$. We observe then that $4\e s_1(\e)^2-10s_1(\e)+(1-\e)=0$, and since $s_1(\e)\in(0,1)$ we get $s_1(\e)=\frac{5}{4\e}(1-\sqrt{1-\frac{4}{25}\e(1-\e)})$ with $s_1(\e)\rightarrow_{\e\rightarrow 0}\frac{1}{10}$. We deduce in particular that $2<p(s_1(\e))-\e<3$ for $\e>0$ small enough.%and we have from (\ref{impbdu}) the following uniform bound:
%\begin{equation}
%\|u_n\|_{L^{p(s_1(\e))-\e}([0,T_\beta]\times\R)}\leq C,
%\label{impbdu1}
%\end{equation}
%with $C>0$ independent on $n$.
\subsubsection*{Estimate in long time of the $L^\infty$ norm of the density $\rho_n$ when $\alpha>\frac{1}{2}$}
We know that $(\rho_n,u_n,v_n)_{n\in\mathbb{N}}$ verify the classical energy estimate on $(0,+\infty)\times \R$ (indeed the initial data belongs to the energy space, in particular $u_n(0,\cdot)\in L^2(\R)$), however the energy estimate depends on the initial data and is not uniform in $n$ (indeed the initial velocity $u_0$ of the Theorem \ref{theo1} is not in $L^2(\R)$). However for $t\geq t_1>0$ the energy estimate will be uniform in $n$. Indeed let us multiply the momentum equation of (\ref{1bis}) by $u_n$ and integrating over $[t_1,t_2]$ with $t_1\in (0,T_\beta]$ and $t_2>t_1$, it gives then:
\begin{equation}
\begin{aligned}
&\frac{1}{2}\int_{\R}(\rho_n |u_n|^2(t_2,x)+(\Pi(\rho_n)-\Pi(\bar{\rho}))(t_2,x) dx+\int^{t_2}_{t_1}\int_{\R}\mu_n(\rho_n(s,x))(\p_x u_n(s,x))^2 ds dx\\
&\hspace{4cm}\leq \frac{1}{2}\int_{\R}(\rho_n(t_1,x) |u_n(t_1,x)|^2(x)+(\Pi(\rho_n(t_1,x))-\Pi(\bar{\rho})) dx.
\end{aligned}
\label{E1}
\end{equation}
From (\ref{energu}), (\ref{impbdu}) and (\ref{BD}) we deduce that there exists $C>0$ independent on $n$:% depending only on the initial data such that:
\begin{equation}
\begin{aligned}
&\frac{1}{2}\int_{\R}(\rho_n |u_n|^2(t_2,x)+(\Pi(\rho_n)-\Pi(\bar{\rho}))(t_2,x) dx+\int^{t_2}_{t_1}\int_{\R}\mu_n(\rho_n(s,x))(\p_x u_n(s,x))^2 ds dx\\
&\hspace{10cm}\leq C(1+\frac{1}{t_1}).
\end{aligned}
\label{E2}
\end{equation}
Combining (\ref{BD1}) and (\ref{E2}) we deduce that for any $t>0$ we have for $C>0$ independent on $n$:% (we use the fact that $\|\p_x\rho_n^{\alpha-\frac{1}{2}}(t,\cdot)\|_{L^2(\R)}\leq\|\sqrt{\rho_n}(t,\cdot)\p_x\va_n(\rho_n(t,\cdot))\|_{L^2(\R)}$):
\begin{equation}
\begin{cases}
\begin{aligned}
%&\|\rho(t,\cdot)-\bar{\rho}\|_{L^2(\R)}\leq C\\
&\|\sqrt{\rho_n}v_n(t,\cdot)\|_{L^2(\R)}\leq C\\
&\|\sqrt{\rho_n} u_n(t,\cdot)\|_{L^2(\R)}\leq C(1+\frac{1}{\sqrt{t}})\\
&\|(\sqrt{s}1_{\{s\leq T_\beta\}}+\sqrt{\mu(\rho_n)}1_{\{s\geq T_\beta\}})\p_x u_n\|_{L^2_t(L^2(\R))}\leq C\\
&\|\p_x\rho_n^{\alpha-\frac{1}{2}}(t,\cdot)\|_{L^2(\R)}+\frac{1}{n}\|\p_x\rho_n^{\theta-\frac{1}{2}}(t,\cdot)\|_{L^2(\R)} \leq C(1+\frac{1}{\sqrt{t}})\\
&\|\p_x\rho_n^{\frac{1}{2}(\gamma+\alpha-1)}\|_{L^2_t(L^2(\R))}+\frac{1}{\sqrt{n}}\|\p_x\rho_n^{\frac{1}{2}(\gamma+\theta-1)}\|_{L^2_t(L^2(\R))}\leq C\\
%&u\in L^{p(s)-\e}([0,T_\beta], H^s(\R))\\
%&\|(\frac{1}{\rho},\rho)(t,\cdot)\|_{L^\infty(\R)}\leq C.
%&\|\sqrt{\rho}\p_x\va(\rho(t,\cdot))\|_{L^2(\R)}\leq C\\
&\|\rho_n(t,\cdot)-\bar{\rho}\|_{L^\gamma_2(\R)}\leq C.
\end{aligned}
\end{cases}
\label{superimpod1}
\end{equation} 
%We deduce then in particular that for any $t\geq t_1$ there exists $c>0$ independent on $n$ such that:
%\begin{equation}
%\begin{aligned}
%&\|\p_x \rho^{\alpha-\frac{1}{2}}(t,\cdot)\|_{L^2(\R)}\leq C.
%\end{aligned}
%\label{superimpod}
%\end{equation}
\begin{lemme}
\label{tech1}
When $\alpha>\frac{1}{2}$ and $\gamma\geq 2\alpha-1$, for any $t\geq 0$ there exists $C>0$ independent on $n$ such that:
\begin{equation}
\|\rho_n(t,\cdot)\|_{L^\infty}\leq C.
\label{normeLin}
\end{equation}
\end{lemme}
{\bf Proof:} %We omit again in the sequel the subscript $n$. 
We know that for $t_1\in (0,T_\beta]$, we have $C>0$ independent on $n$ such that:
\begin{equation}
\|\rho_n\|_{L^\infty([0,t_1],L^\infty(\R))}\leq C.
\label{tempspetit}
\end{equation}
Following the lemma 3.7 in \cite{Jiu}, it remains now only to prove that there exists $C>0$ independent on $n$ such that for all $t\geq t_1>0$, we have:
\begin{equation}
\|\rho_n(t,\cdot)\|_{L^\infty}\leq C.
\end{equation}
Since $\rho_n$ is regular, %(we recall that we forget the subscript $n$) 
it yields for any $t\geq t_1$ and for $\beta\geq 0$ to determine later:
\begin{equation}
\begin{aligned}
&(\rho_n^{\alpha-\frac{1}{2}}(t,x)-\bar{\rho}^{\alpha-\frac{1}{2}})^{2\beta}(t,x)=\int_{-\infty}^x \p_y ((\rho_n^{\alpha-\frac{1}{2}}(t,y)-\bar{\rho}^{\alpha-\frac{1}{2}})^{2\beta}) (t,y)dy\\
&\leq 2\beta\big(\int_{-\infty}^x |\rho_n^{\alpha-\frac{1}{2}}(t,y)-\bar{\rho}^{\alpha-\frac{1}{2}}|^{2(2\beta-1)} dy\big)^{\frac{1}{2}}\big(\int_{-\infty}^{x}|\p_y \rho_n^{\alpha-\frac{1}{2}}(t,y)|^2 dy\big)^{\frac{1}{2}}.
\end{aligned}
\label{etape1}
\end{equation}
From (\ref{superimpod1}) we know that $\int_{-\infty}^{x}|\p_y \rho_n^{\alpha-\frac{1}{2}}(t,y)|^2 dy$ is uniformly bounded in $n$ for $t\geq t_1>0$.\\
We have now to estimate the following term:
$$
\begin{aligned}
&\int_{-\infty}^x |\rho_n^{\alpha-\frac{1}{2}}(t,y)-\bar{\rho}^{\alpha-\frac{1}{2}}|^{2(2\beta-1)} dy=\int_{-\infty}^x |\rho_n^{\alpha-\frac{1}{2}}(t,y)-\bar{\rho}^{\alpha-\frac{1}{2}}|^{2(2\beta-1)}[1_{\{|\rho_n-\bar{\rho}|\leq\frac{\bar{\rho}}{2}\}}\\
&\hspace{7cm}+1_{\{|\rho_n-\bar{\rho}|>\frac{\bar{\rho}}{2}\}}]  dy=I_1(t)+I_2(t).
\end{aligned}
$$
We have now for $C,C_1>0$ sufficiently large and taking $\beta=1$ for any $t\geq t_1$ and using
(\ref{superimpod1}):
\begin{equation}
\begin{aligned}
I_1(t)&=\int_{-\infty}^x |\rho_n^{\alpha-\frac{1}{2}}(t,y)-\bar{\rho}^{\alpha-\frac{1}{2}}|^{2(2\beta-1)}1_{\{|\rho_n-\bar{\rho}|\leq\frac{\bar{\rho}}{2}\}} dy
\\
&\leq C\int_{-\infty}^x |\rho_n(t,y)-\bar{\rho}|^{2}1_{\{|\rho_n-\bar{\rho}|\leq\frac{\bar{\rho}}{2}\}} dy\\
&\leq C_1.
\end{aligned}
\label{etape2}
\end{equation}
For $\beta=1$ we have now when $t\geq t_1$ and $C,C_1>0$ sufficiently large using (\ref{superimpod1}):% and the fact that $\gamma\geq 2\alpha-1$:
\begin{equation}
\begin{aligned}
I_2(t)&=\int_{-\infty}^x |\rho_n^{\alpha-\frac{1}{2}}(t,y)-\bar{\rho}^{\alpha-\frac{1}{2}}|^{2(2\beta-1)}1_{\{|\rho_n-\bar{\rho}|>\frac{\bar{\rho}}{2}\}} dy\\
&\leq C\int_{-\infty}^x |\rho_n(t,y)-\bar{\rho}|^{2(2\beta-1)(\alpha-\frac{1}{2})}1_{\{|\rho_n-\bar{\rho}|>\frac{\bar{\rho}}{2}\}} dy\\
&\leq C\int_{-\infty}^x |\rho_n(t,y)-\bar{\rho}|^{2(\alpha-\frac{1}{2})}1_{\{|\rho_n-\bar{\rho}|>\frac{\bar{\rho}}{2}\}} dy.
%&\leq C \|(\rho(t,\cdot)-\bar{\rho})1_{\{|\rho-\bar{\rho}|>\frac{\bar{\rho}}{2}\}}\|_{L^\gamma}.
\end{aligned}
\label{etape3}
\end{equation}
Since $2\alpha-1\leq \gamma$ it exists $C_1>0$ large enough such that:
\begin{equation}
\begin{aligned}
I_2(t)&\leq C_1(\int_{-\infty}^x |\rho_n(t,y)-\bar{\rho}|^{\gamma}1_{\{|\rho_n-\bar{\rho}|>\frac{\bar{\rho}}{2}\}} dy)^{\frac{2\alpha-1}{\gamma}} |\{|\rho_n-\bar{\rho}|>\frac{\bar{\rho}}{2}\}|^{\frac{\gamma-2\alpha+1}{\gamma}}\\
&\leq C_1\|\rho_n(t,\cdot)-\bar{\rho}\|_{L^\gamma_2(\R)}^\gamma.
%&\leq C \|(\rho(t,\cdot)-\bar{\rho})1_{\{|\rho-\bar{\rho}|>\frac{\bar{\rho}}{2}\}}\|_{L^\gamma}.
\end{aligned}
\label{etape4}
\end{equation}
From (\ref{tempspetit}), (\ref{etape1}), (\ref{etape2}) and (\ref{etape4}) we conclude the proof of (\ref{normeLin}). $\blacksquare$\\
\\
%\begin{remarka}
%ON PEUT UTILISER UN ARGUMENT DE PETITESSE POUR CONTROLER LE VIDE.
%\end{remarka}
We deduce using the estimate (\ref{cru1}) that for any $t>0$ and any $n\in\mathbb{N}^*$ we have for $C(t)>0$ depending only on $t$:
\begin{equation}
\|\rho_n u_n(t,\cdot)\|_{L^1(\R)}+\|\rho_n v_n(t,\cdot)\|_{L^1(\R)}\leq C(t).
\label{gainL1infini}
\end{equation}
\subsubsection*{Estimate in long time of the $L^\infty$ norm of $\rho_n$ and $\frac{1}{\rho_n}$ when $\alpha<\frac{1}{2}$}
We recall that we have now:
$$
\begin{aligned}
&\p_t u_n+u_n\p_x u_n-\p_x (\frac{\mu_n(\rho_n)}{\rho_n}%\rho^{\alpha-1}
\p_x u_n)-\p_x\va_n(\rho_n)\,\p_x u_n%\frac{\mu}{\alpha-1}\p_x(\rho^{\alpha-1})\p_x u
+\frac{a\gamma}{\gamma-1}\p_x(\rho_n^{\gamma-1})=0.
\end{aligned}
$$
We deduce that we have for $t\in]t_1,T_\beta]$ with $0<t_1<T_\beta$:
$$
\begin{aligned}
&\p_t((t-t_1)\p_x u_n)-\p_x u_n+\p_x(u_n\,(t-t_1)\p_x u_n)-\p_x \big(\frac{\mu_n(\rho_n)}{\rho_n}
\p_x ((t-t_1)\p_x u_n)\big)\\
&-\p_x\big( (\p_x\va(\rho_n)+\p_x(\frac{\mu_n(\rho_n)}{\rho_n}))\,(t-t_1)\p_x u_n\big)%\frac{\mu}{\alpha-1}\p_x(\rho^{\alpha-1})\p_x u
+\frac{a\gamma}{\gamma-1}\p_{xx}((t-t_1)\rho_n^{\gamma-1})=0.
%\mu\p_x (\rho^{\alpha-1}\p_x ((t-t_1)\p_x u))\\
%&-\mu\p_x(\p_x(\rho^{\alpha-1})\,(t-t_1)\p_x u)-\frac{\mu}{\alpha-1}\p_x(\p_x(\rho^{\alpha-1})(t-t_1) \p_x u)+\frac{a\gamma}{\gamma-1}\p_{xx}((t-t_1)\rho^{\gamma-1})=0
\end{aligned}
$$
%$$
%\begin{aligned}
%&\p_t((t-t_1)\p_x u)-\p_x u+\p_x(u\,(t-t_1)\p_x u)-\mu\p_x (\rho^{\alpha-1}\p_x ((t-t_1)\p_x u))\\
%&-\mu\p_x(\p_x(\rho^{\alpha-1})\,(t-t_1)\p_x u)-\frac{\mu}{\alpha-1}\p_x(\p_x(\rho^{\alpha-1})(t-t_1) \p_x u)+\frac{a\gamma}{\gamma-1}\p_{xx}((t-t_1)\rho^{\gamma-1})=0
%\end{aligned}
%$$
We set now $w_n=(t-t_1)\p_x u_n$ and we have for $t\in]t_1,T_\beta]$ :
$$
\begin{aligned}
&\p_t w_n-\p_x u_n+\p_x(u_n\,w_n)-\p_x (\frac{\mu_n(\rho_n)}{\rho_n}\p_x w_n)-\p_x(\frac{\mu_n'(\rho_n)}{\rho} \p_x\rho_n\,w_n)\\
&\hspace{7cm}+\frac{a\gamma}{\gamma-1}\p_{xx}((t-t_1)\rho_n^{\gamma-1})=0
\end{aligned}
$$
%$$
%\begin{aligned}
%&\p_t v-\p_x u+\p_x(u\,v)-\mu\p_x (\rho^{\alpha-1}\p_x v)-\frac{\mu \alpha}{\alpha-1}\p_x(\p_x(\rho^{\alpha-1})v)+\frac{a\gamma}{\gamma-1}\p_{xx}((t-t_1)\rho^{\gamma-1})=0
%\end{aligned}
%$$
Multiplying by $w_n$ and integrating over $(t_1,t)\times\R$, we obtain:
$$
\begin{aligned}
&\frac{1}{2}\int_{\R} |w_n|^2(t,x) dx+ \int^t_{t_1}\int_{\R}\frac{\mu_n(\rho_n)}{\rho_n}(s,x)|\p_x w_n(s,x)|^2 ds dx- \int^t_{t_1}\int_{\R}\p_x u_n(s,x) w_n(s,x) ds dx\\
&-\int^t_{t_1}\int_{\R}u_n(s,x) w_n(s,x) \p_x w_n(s,x) ds dx + \int^t_{t_1}\int_{\R}\frac{\mu_n'(\rho_n)}{\rho_n} \p_x\rho_n (s,x) w_n(s,x)\p_x w_n(s,x) ds dx\\
&\hspace{3cm}-
\frac{a\gamma}{\gamma-1}\int^t_{t_1}\int_{\R} (s-t_1) \p_{x}\rho_n^{\gamma-1}(s,x)\p_x w_n(s,x) ds dx=0.
\end{aligned}
$$
Using (\ref{impbdu}) we obtain then on $(t_1,t)$ with $t\in]t_1,T_\beta]$, $c>0$ and $C>0$ independent on $n$:
$$
\begin{aligned}
&\frac{1}{2}\int_{\R} |w_n|^2(t,x) dx+\mu c \int^t_{t_1}\int_{\R}|\p_x w_n(s,x)|^2 ds dx\leq \|\p_x u_n\|_{L^2((t_1,t),L^2(\R))}\sqrt{t} \|w_n\|_{L^\infty((t_1,t),L^2(\R))}\\
&+\|\frac{\mu_n'(\rho_n)}{\rho}\p_x\rho_n\|_{L^\infty((t_1,t),L^2(\R))}\|\p_x w_n\|_{L^2((t_1,t),L^2(\R))}\|w_n\|_{L^2((t_1,t),L^\infty (\R))}\\
&\hspace{4cm}+\|u_n\|_{L^\infty((t_1,t),L^2(\R))} \|\p_x w_n\|_{L^2((t_1,t),L^2(\R))}\|w_n\|_{L^2((t_1,t),L^\infty (\R))}\\
&\hspace{6cm}+Ct \|\p_x\rho_n^{\gamma-1}\|_{L^2((t_1,t),L^2(\R))}\|\p_x w_n\|_{L^2((t_1,t),L^2(\R))}.
\end{aligned}
$$
Next by Gagliardo-Niremberg inequality:
$$
\begin{aligned}
&\frac{1}{2}\int_{\R} |w_n|^2(t,x) dx+\mu c \int^t_{t_1}\int_{\R}|\p_x w_n(s,x)|^2 ds dx\leq\\
&\|\p_x u_n\|_{L^2((t_1,t),L^2(\R))}\sqrt{t} \|w_n\|_{L^\infty((t_1,t),L^2(\R))}+Ct \|\p_x\rho_n^{\gamma-1}\|_{L^2((t_1,t),L^2(\R))}\|\p_x w_n\|_{L^2((t_1,t),L^2(\R))}\\
&+C(\|\frac{\mu_n'(\rho_n)}{\rho_n}\p_x\rho_n\|_{L^\infty((t_1,t),L^2(\R))}+\|u_n\|_{L^\infty((t_1,t),L^2(\R))})  \|\p_x w_n\|^{\frac{3}{2}}_{L^2((t_1,t),L^2(\R))}\|w_n\|^{\frac{1}{2}}_{L^2((t_1,t),L^2(\R))}.
\end{aligned}
$$
Now using Young inequality we deduce that for $C_1>0$ large enough:
$$
\begin{aligned}
&\frac{1}{4}\|w_n\|_{L^\infty((t_1,t),L^2(\R))}^2+\frac{\mu c }{2}\|\p_x w_n\|^2_{L^2((t_1,t),L^2(\R))}\\
&\leq C_1(\|\p_x u_n\|^2_{L^2((t_1,t),L^2(\R))}t+t^2 \|\p_x\rho_n^{\gamma-1}\|_{L^2((t_1,t),L^2(\R))}^2\\
&\hspace{2cm}+ (\|\frac{\mu_n'(\rho_n)}{\rho_n}\p_x\rho_n\|_{L^\infty((t_1,t),L^2(\R))}+\|u_n\|_{L^\infty((t_1,t),L^2(\R))})^4 \|w_n\|^{2}_{L^2((t_1,t),L^2(\R))}\big).
\end{aligned}
$$
From (\ref{impbdu}) and (\ref{superimpod1}) we deduce that for $t\in (t_1,T_\beta]$, we have for $C(t,t_1)>0$ independent on $n$:
\begin{equation}
\begin{aligned}
&\frac{1}{4}\|w_n\|_{L^\infty((t_1,t),L^2(\R))}^2+\frac{\mu c }{2}\|\p_x w_n\|_{L^2((t_1,t),L^2(\R))}\leq C(t,t_1).%\\
%&\leq C_1(t(1+\frac{1}{t_1})+t^2(t-t_1)(1+\frac{1}{t_1}) +(1+\frac{1}{t_1})^2(t-t_1)^2(1+\frac{t}{t_1})\big).
\end{aligned}
\end{equation}
We have showed that for $t_1\in]0,T_\beta]$ %there exists $C(t,t_1)>0$ independent on $n$ and depending only on the initial data $(\rho_0,v_0,u_0)$ such that for 
and any $t\in]t_1,T_\beta]$ and any $n\in\mathbb{N}^*$ we have:
%$$
%\begin{aligned}
%&\frac{1}{4}\|v\|_{L^\infty((t_1,t),L^2(\R))}^2+\frac{\mu c }{2}\|\p_x v\|_{L^2((t_1,t),L^2(\R))}\\
%&\leq C_1(\frac{t}{t_1}+t^2 \|\p_x\rho^{\gamma-1}\|_{L^2((t_1,t),L^2(\R))}^2\\
%&\hspace{4cm}+\|\p_x\rho^{\alpha-1}\|^4_{L^\infty((t_1,t),L^2(\R))}\|v\|^{2}_{L^2((t_1,t),L^2(\R))}\big).
%\end{aligned}
%$$
\begin{equation}
\|(t-t_1) \p_x u_n(t,\cdot)\|_{L^2}\leq C (t,t_1)<+\infty.
\label{last}
\end{equation}
\begin{remarka}
It is important to mention in fact that the previous estimate is true for any $\alpha>0$.
\end{remarka}
From (\ref{last}) and (\ref{impbdu}) we have obtained that for $t_2\in]t_1,T_\beta]$ with $0<t_1<T_\beta$, we have for $C>0$ independent on $n$ and depending only on the initial data $(\rho_0,u_0,v_0)$:
\begin{equation}
\begin{cases}
\begin{aligned}
&\|\p_x u_n(t_2,\cdot)\|_{L^2(\R)}\leq C\\
&\| u_n(t_2,\cdot)\|_{L^2(\R)}\leq C\\
&\|\rho_n(t_2,\cdot)-\bar{\rho}\|_{H^1(\R)}\leq C\\
&\|(\rho_n,\frac{1}{\rho_n})(t_2,\cdot)\|_{L^\infty(\R)}\leq C.
\end{aligned}
\end{cases}
\label{Mellet}
\end{equation}
From the theorem of Mellet and Vasseur in \cite{MV}, we know that there exists a global strong solution
$(\rho_{\mu,n},u_{\mu,n})$  with initial data $(\rho_n(t_2,\cdot),u_n(t_2,\cdot))$. In addition it is proved in \cite{MV} that for $t\geq 0$ there exists $C(t)>0$ depending only on $t$, $ \| u_n(t_2,\cdot)\|_{H^1(\R)}$,
$\|\rho_n(t_2,\cdot)-\bar{\rho}\|_{H^1(\R)}$ and $\|(\rho_n,\frac{1}{\rho_n})(t_2,\cdot)\|_{L^\infty(\R)}$ such that for any $t\geq 0$ we have:
\begin{equation}
\|(\frac{1}{\rho_{\mu,n}},\rho_{\mu,n})(t,\cdot)\|_{L^\infty(\R)}\leq C(t). 
\label{solMell}
\end{equation}
We deduce that  when $0<\alpha<\frac{1}{2}$ the regular solution $(\rho_n,u_n)$ verify for any $t\geq t_2$:
\begin{equation}
(\rho_n(t,\cdot),u_n(t,\cdot))=(\rho_{\mu,n}(t-t_2,\cdot),u_{\mu,n}(t-t_2,\cdot)),
\end{equation}
This is a direct consequence of the uniqueness of the solution $(\rho_{\mu,n},u_{\mu,n})$, indeed $(\rho_n,u_n)$ is sufficiently regular such that the solution is in the class of uniqueness of $(\rho_{\mu,n},u_{\mu,n})$. It implies from (\ref{solMell}) that there exists $C(t)$ depending only on $ \| u_n(t_2,\cdot)\|_{H^1(\R)}$,
$\|\rho_n(t_2,\cdot)-\bar{\rho}\|_{H^1(\R)}$ and $\|(\rho_n,\frac{1}{\rho_n})(t_2,\cdot)\|_{L^\infty(\R)}$ (it is important to mention that we have seen that these norms are independent on $n$, in particular $C(t)$ depends only on $t$ and $(\rho_0,u_0,v_0)$) such that for any $t\geq 0$ we have:% DIRE QUE LA FORME DU $\mu_n$ NE PERTURBE PAS L'ARGUMENT
\begin{equation}
\|(\frac{1}{\rho_n},\rho_n)(t,\cdot)\|_{L^\infty(\R)}\leq C(t). 
\label{asolMell}
\end{equation}
From (\ref{cru1}) it yields that for any $t>0$ there exists $C(t)>0$ independent on $n$ such that:
\begin{equation}
\|\rho_n u_n(t,\cdot)\|_{L^1(\R)}+ \|\rho_n v_n(t,\cdot)\|_{L^1(\R)}\leq C(t). 
\label{solMell1}
\end{equation}
%with:
%$$\Pi(s)=s(\int^s_{\bar{\rho}}\frac{P(z)}{z^2}dz-\frac{P(\bar{\rho})}{\bar{\rho}}).$$
\subsubsection*{Gain of integrability on the velocity $u_n$ in long time when $\alpha>\frac{1}{2}$}
Multiplying the momentum equation of (\ref{1bis}) by $u_n|u_n|^p$ (with $p=2$) and integrating over $(t_1,t)\times\R$ with
$t_1\in (0,T_\beta]$ and $t>t_1$, we have then: 
\begin{equation}
\begin{aligned}
&\frac{1}{p+2}\int_{\R}\rho_n(t,x) |u_n|^{p+2}(t,x)dx+(p+1)\int^{t}_{t_1}\int_\R \mu_n(\rho_n(s,x))|\p_x u_n(s,x)|^2 |u_n(s,x)|^p  ds dx\\
%&+\int^{T_\beta}_0\int_\R \mu(\rho(s,x))\p_x u(s,x) |u(s,x)|^pu(s,x) \p_x\varphi(x) ds dx\\
&+\int^ {t}_{t_1}\int_\R   \p_x P(\rho_n)(s,x)\, u_n |u_n|^p(s,x) ds dx=\frac{1}{p+2}\int_{\R}\rho(t_1,x) |u_n|^{p+2}(t_1,x)dx.
\end{aligned}
\label{3.55bis}
\end{equation}
From (\ref{last}) and (\ref{impbdu}) (indeed (\ref{last}) is true for any $\alpha>0$), and using interpolation we know that there exists $C(t_1)>0$ independent on $n$ such that:
\begin{equation}
|\int_{\R}\rho_n(t_1,x) |u_n|^{p+2}(t_1,x)dx|\leq C(t_1).
\label{3.83}
\end{equation}
%\begin{equation}
%\begin{aligned}
%&\frac{1}{p+2}\int_{\R} |u|^{4}(t,x) dx
%+(p+1)\int^{t}_0\int_\R \mu(\rho(s,x))|\p_x u(s,x)|^2 |u(s,x)|^2 ds dx\\
%&+\int^ {t}_{0}\int_\R   \p_x P(\rho)(s,x)\, u |u|^2(s,x) ds dx\leq C(1+\frac{1}{t_1^{\frac{1}{4}}})^4.
%\end{aligned}
%\end{equation}
It remains now only to estimate the term coming from the pressure. We have then using Young inequality for $\e>0$ and $C_\e>0$ large enough with $p=2$:
\begin{equation}
\begin{aligned}
&|\int^ {t}_{t_1}\int_\R   \p_x P(\rho_n)(s,x)\, u_n |u_n|^p(s,x) ds dx|\\
&\leq a(p+1) |\int^ {t}_{t_1}\int_\R   \rho_n^{\gamma}(s,x)\,\p_x u_n(s,x) |u_n|^p(s,x)ds dx|\\
&\leq \e\int^{t}_{t_1}\int_\R \mu_n(\rho_n(s,x))|\p_x u_n(s,x)|^2 |u_n(s,x)|^2  ds dx\\
&\hspace{3cm}+C_\e\int^{t}_{t_1}\int_\R \frac{\rho_n^{2\gamma-1}(s,x)}{\mu_n(\rho_n(s,x))}\rho_n(s,x)|u_n|^2(s,x) ds dx.
\end{aligned}
\end{equation}
Since we have $2\gamma-1\geq \alpha$, we deduce from (\ref{impbdu}) and (\ref{normeLin}) that we have:
\begin{equation}
\begin{aligned}
&|\int^ {t}_{t_1}\int_\R   \p_x P(\rho_n)(s,x)\, u_n |u_n|^p(s,x) ds dx|\\
%&\leq a(p+1) |\int^ {t}_{t_1}\int_\R   \rho^{\gamma}(s,x)\,\p_x u(s,x) |u|^p(s,x)ds dx|\\
&\leq \e\int^{t}_{t_1}\int_\R \mu_n(\rho(s,x))|\p_x u_n(s,x)|^2 |u_n(s,x)|^2  ds dx\\
&\hspace{3cm}+C_\e\|\rho_n\|^{2\gamma-1-\alpha}_{L^\infty(\R^+,L^\infty(\R))}(t-t_1) \|\sqrt{\rho_n}u_n\|_{L^\infty([t_1,+\infty),L^2(\R))}^2.
\end{aligned}
\label{gainppp}
\end{equation}
From (\ref{superimpod1}), (\ref{gainppp}), (\ref{3.83}) and (\ref{3.55bis}) we deduce that for any $n\in\mathbb{N}^*$ and for any $t_1\in]0,T_\beta]$ and $t>t_1$ we have:
\begin{equation}
\|\rho_n^{\frac{1}{4}}u_n\|_{L^\infty([t_1,t],L^4(\R))}\leq C(t_1,t),
\label{gainfin1}
\end{equation}
with $C(t_1,t)$ depending only on $t$, $t_1$ but independent on $n$.
\subsubsection*{Compactness when $0<\alpha<\frac{1}{2}$}
We wish now to prove that the sequence $(\rho_n,u_n)_{n\in\mathbb{N}}$ converges up to a subsequence to a global weak solution $(\rho,u)$ of the system (\ref{1}). Using (\ref{superimpod1}), (\ref{impbdu}), (\ref{solMell}) and (\ref{solMell1})) we have seen 
that there exists $C>0, C(t)>0$ independent on $n$ (with $C(t)$ depending only on $t$) such that we have for any $t>0$:
\begin{equation}
\begin{cases}
\begin{aligned}
%&\|\rho(t,\cdot)-\bar{\rho}\|_{L^2(\R)}\leq C\\
&\|\sqrt{\rho_n}v_n(t,\cdot)\|_{L^2(\R)}\leq C\\
&\|\sqrt{\rho_n} u_n(t,\cdot)\|_{L^2(\R)}\leq C(1+\frac{1}{\sqrt{t}})\\
&\|(\sqrt{s}1_{\{s\leq T_\beta\}}+\sqrt{\mu(\rho_n)}1_{\{s\geq T_\beta\}})\p_x u_n\|_{L^2_t(L^2(\R))}\leq C\\
&\|\p_x\rho_n^{\alpha-\frac{1}{2}}(t,\cdot)\|_{L^2(\R)}\leq C(1+\frac{1}{\sqrt{t}})\\
&\|\p_x\rho_n^{\frac{1}{2}(\gamma+\alpha-1)}\|_{L^2_t(L^2(\R))}\leq C\\
&\|u_n\|_{L^{p(s)-\e}([0,T_\beta], H^s(\R))}\leq C\\
%&u\in L^{p(s)-\e}([0,T_\beta], H^s(\R))\\
%&\|(\frac{1}{\rho},\rho)(t,\cdot)\|_{L^\infty(\R)}\leq C.
%&\|\sqrt{\rho}\p_x\va(\rho(t,\cdot))\|_{L^2(\R)}\leq C\\
&\|\rho_n(t,\cdot)-\bar{\rho}\|_{L^\gamma_2(\R)}\leq C\\
&\|(\frac{1}{\rho_n},\rho_n)(t,\cdot)\|_{L^\infty(\R)}\leq C(t)\\
&\|\rho_n u_n(t,\cdot)\|_{L^1(\R)}+ \|\rho_n v_n(t,\cdot)\|_{L^1(\R)}\leq C(t). 
\end{aligned}
\end{cases}
\label{superimpod12}
\end{equation} 
with $s\in(0,1)$, $p(s)=\frac{3}{1+2s}$ and $\e>0$ sufficiently small such that $p(s)-\e\geq 1$. From 
(\ref{superimpod12}), we deduce that there exists $C(T)>0$ independent on $n$ such that for any $T>0$ we have:
\begin{equation}
\|\rho_n-\bar{\rho}\|_{L^2((0,T),H^1(\R))}\leq C(T).
\label{compac1}
\end{equation}
In addition we have $\p_t(\rho_n-\bar{\rho})=-\p_x(\rho_n u_n)$, then from (\ref{superimpod12}) $\p_t(\rho_n-\bar{\rho})$ is uniformly bounded in $L^1((0,T),H^{-1}(\R))$. From the Lemma \ref{Aubin}, the diagonal process and compact Sobolev embedding, it implies that up to a subsequence $(\rho_n)_{n\in\mathbb{N}}$ converges strongly to $\rho$ in $L^2_{loc}(\R^+\times\R)$. We deduce that up to a subsequence $(\rho_n)_{n\in\mathbb{N}}$ converges almost everywhere to $\rho$ on $]0,+\infty[\times\R$. From Fatou lemma and (\ref{superimpod12}), we obtain that there exists $C(t)>$ depending only on $t>0$ such that:\begin{equation}
\begin{cases}
\begin{aligned}
%&\|\rho(t,\cdot)-\bar{\rho}\|_{L^2(\R)}\leq C\\
%&\|\sqrt{\rho_n}v_n(t,\cdot)\|_{L^2(\R)}\leq C\\
%&\|\sqrt{\rho_n} u_n(t,\cdot)\|_{L^2(\R)}\leq C(1+\frac{1}{\sqrt{t}})\\
%&\|(\sqrt{s}1_{\{s\leq T_\beta\}}+\sqrt{\mu(\rho_n)}1_{\{s\geq T_\beta\}})\p_x u_n\|_{L^2_t(L^2(\R))}\leq C\\
%&\|\p_x\rho_n^{\alpha-\frac{1}{2}}(t,\cdot)\|_{L^2(\R)}\leq C(1+\frac{1}{\sqrt{t}})\\
%&\|\p_x\rho_n^{\frac{1}{2}(\gamma+\alpha-1)}\|_{L^2_t(L^2(\R))}\leq C\\
%&\|u_n\|_{L^{p(s)-\e}([0,T_\beta], H^s(\R))}\leq C\\
%&u\in L^{p(s)-\e}([0,T_\beta], H^s(\R))\\
%&\|(\frac{1}{\rho},\rho)(t,\cdot)\|_{L^\infty(\R)}\leq C.
%&\|\sqrt{\rho}\p_x\va(\rho(t,\cdot))\|_{L^2(\R)}\leq C\\
&\|\rho(t,\cdot)-\bar{\rho}\|_{L^\gamma_2(\R)}\leq C\\
&\|(\frac{1}{\rho},\rho)(t,\cdot)\|_{L^\infty(\R)}\leq C(t).
%&\|\rho_n u_n(t,\cdot)\|_{L^1(\R)}+ \|\rho_n v_n(t,\cdot)\|_{L^1(\R)}\leq C(t). 
\end{aligned}
\end{cases}
\label{superimpod12a}
\end{equation} 
Next from (\ref{superimpod12}), we deduce that $(u_n)_{n\in\mathbb{N}}$ is uniformly bounded in $L^{p(s)-\e}([0,T_\beta], H^s(\R))$ and $u_n$ is uniformly bounded in $L^2_{loc}([T_\beta,+\infty),H^1(\R))$. Using classical paraproduct laws (see\cite{Danchin}) we obtain that for $t\in ]0,T_\beta]$ and $0<s<\frac{1}{2}$:
$$
\begin{aligned}
&\|T_{(\rho_n(t,\cdot)-\bar{\rho})}u_n(t,\cdot)\|_{H^s(\R)}\leq \|\rho_n(t,\cdot)-\bar{\rho}\|_{L^\infty(\R)}\|u_n(t,\cdot)\|_{H^s(\R)}\\
&\|T_{u_n(t,\cdot)}(\rho_n(t,\cdot)-\bar{\rho})\|_{H^{\frac{1}{2}+s}(\R)}\leq \|\rho_n(t,\cdot)-\bar{\rho}\|_{H^1(\R)}\|u_n(t,\cdot)\|_{H^s(\R)}\\
&\|R(u_n(t,\cdot),(\rho_n(t,\cdot)-\bar{\rho}))\|_{H^{\frac{1}{2}+s}(\R)}\leq \|\rho_n(t,\cdot)-\bar{\rho}\|_{H^1(\R)}\|u_n(t,\cdot)\|_{H^s(\R)}\\
\end{aligned}
$$
We deduce from (\ref{superimpod12}) and (\ref{compac1}) that $\rho_n u_n$ is uniformly bounded in $L^{q(s)}([0,T_\beta],H^{s}(\R))+L^{p(s)-\e}([0,T_\beta],H^s(\R))$ with $0<s<\frac{1-\e}{2(2+\e)}$ and $\frac{1}{q(s)}=\frac{1}{2}+\frac{1}{p(s)-\e}$ (  $p(s)=\frac{3}{1+2s}$ and $\e>0$ sufficiently small and $p(s)-\e\geq 2$). In particular since $p(s)-\e\geq q(s)\geq 1$, we have proved that $\rho_n u_n$ is uniformly bounded in $L^{q(s)}([0,T_\beta],H^{s}(\R))$.\\
From (\ref{superimpod12}), we observe that $u_n$ and $(\rho_n-\bar{\rho})$ are uniformly bounded in $L^2_{loc}([T_\beta,+\infty[,H^1(\R))$, by product of Sobolev space we deduce that $\rho_n u_n$ is bounded in $L^1_{loc}([T_\beta,+\infty[,H^1(\R))$. We have finally proved that $\rho_n u_n$ is bounded in
$L^1_{loc}(\R^+,H^s(\R))$ for $0<s<\frac{1-\e}{2(2+\e)}$ with $\e>0$ sufficiently small.
We recall now that we have:
$$\p_t(\rho_n u_n)=-\p_x(\rho_n u_n^2)+\p_{xx}(\mu_n(\rho_n)u_n)-\p_x (u_n\p_x(\mu_n(\rho_n)))%((\mu\rho_n^\alpha+\frac{1}{n}\rho_n^\theta)\p_x u_n)
-\p_x(a(\rho_n^\gamma -\bar{\rho}^\gamma)).$$
Let us deal with the term $\p_x (u_n\p_x(\mu_n(\rho_n)))$, from (\ref{superimpod12}) we know that 
$\p_x(\mu_n(\rho_n))$ is uniformly bounded in $L^2(\R^+,L^2(\R))$ and $u_n$ is uniformly bounded in $L^2_{loc}(\R^+,L^2(\R))$, we deduce that  $\p_x (u_n\p_x(\mu_n(\rho_n)))$ is uniformly bounded in $L_{loc}^1(\R^+,W^{-1,1}(\R))$. We proceed similarly for the other terms and we deduce by Sobolev embedding that $\p_t(\rho_n u_n)$ is uniformly bounded in $L^1_{loc}(\R^+,H^{-s_1}(\R))$ with $s_1$ sufficiently large.\\
Since we have seen that  $(\rho_n u_n)_{n\in\mathbb{N}}$ is bounded in
$L^1_{loc}(\R^+,H^s(\R))$ for $0<s<\frac{1-\e}{2(2+\e)}$ with $\e>0$ sufficiently small, using the lemma \ref{Aubin} we show that $(\rho_nu_n)_{n\in\mathbb{N}}$ converges strongly up to a subsequence to $m$ in 
$L^1_{loc}(\R^+,L^2_{loc}(\R))$. It implies in particular that up to a subsequence $(\rho_n u_n)_{n\in\mathbb{N}}$ converges almost everywhere to $m$ in $]0,+\infty[\times\R$. We define now $u$ as $u=\frac{m}{\rho}$ which is well defined since $\frac{1}{\rho}$ belongs to $L^\infty(\R^+\times\R)$ (see (\ref{superimpod12a}) and we have that $(u_n)_{n\in\mathbb{N}}$ converges up to a subsequence to $u$.\\
From Fatou lemma, we have for any $t>0$:
\begin{equation}
\begin{cases}
\begin{aligned}
%&\|\rho(t,\cdot)-\bar{\rho}\|_{L^2(\R)}\leq C\\
%&\|\sqrt{\rho_n}v_n(t,\cdot)\|_{L^2(\R)}\leq C\\
&\|\sqrt{\rho} u(t,\cdot)\|_{L^2(\R)}\leq C(1+\frac{1}{\sqrt{t}})\\
%&\|(\sqrt{s}1_{\{s\leq T_\beta\}}+\sqrt{\mu(\rho_n)}1_{\{s\geq T_\beta\}})\p_x u_n\|_{L^2_t(L^2(\R))}\leq C\\
%&\|\p_x\rho_n^{\alpha-\frac{1}{2}}(t,\cdot)\|_{L^2(\R)}\leq C(1+\frac{1}{\sqrt{t}})\\
%&\|\p_x\rho_n^{\frac{1}{2}(\gamma+\alpha-1)}\|_{L^2_t(L^2(\R))}\leq C\\
%&\|u_n\|_{L^{p(s)-\e}([0,T_\beta], H^s(\R))}\leq C\\
%&u\in L^{p(s)-\e}([0,T_\beta], H^s(\R))\\
%&\|(\frac{1}{\rho},\rho)(t,\cdot)\|_{L^\infty(\R)}\leq C.
%&\|\sqrt{\rho}\p_x\va(\rho(t,\cdot))\|_{L^2(\R)}\leq C\\
%&\|\rho(t,\cdot)-\bar{\rho}\|_{L^\gamma_2(\R)}\leq C\\
%&\|(\frac{1}{\rho},\rho)(t,\cdot)\|_{L^\infty(\R)}\leq C(t).
&\|m(t,\cdot)\|_{L^1(\R)}\leq C(t). 
\end{aligned}
\end{cases}
\label{superimpod12b}
\end{equation} 
From (\ref{superimpod12}) we observe that $(u_n)_{n\in\mathbb{N}}$ is uniformly bounded in $L^{p(s)-\e}([0,T_\beta], H^s(\R))$ and is also uniformly bounded in $L^\infty([T_\beta,+\infty[,L^2(\R))\cap L^2([T_\beta,+\infty[,\dot{H}^1(\R))$. By Gagliardo-Niremberg we deduce that $(u_n)_{n\in\mathbb{N}}$ is uniformly bounded in $L^4([T_\beta,+\infty[,L^\infty(\R))$. By interpolation it implies that $u_n$ is uniformly bounded in $L^6([T_\beta,+\infty[,L^6(\R))$. Now taking $s=\e_1>0$ sufficiently small and using the fact that $(u_n)_{n\in\mathbb{N}}$ is uniformly bounded in $L^{p(s)-\e}([0,T_\beta], H^s(\R))$, we deduce by Sobolev embedding that $(u_n)_{n\in\mathbb{N}}$ is uniformly bounded in $L^{2+\e_2}([0,T_\beta], L^{p_2}(\R))$ with $\e_2>0$ sufficiently small and $p_2>2$. It implies in particular that there exists $C>0$ independent on $n$ and $\e_2>0$ sufficiently small such that:
\begin{equation}
\|u_n\|_{L^{2+\e_2}_{loc}(\R^+\times\R)}\leq C.
\label{gainint}
\end{equation}
Using the Lemma \ref{lemmeimp}, Lemma \ref{Aubin}, (\ref{superimpod12}), (\ref{gainint}) and the fact that up to a subsequence $(\rho_n,u_n,m_n)_{n\in\mathbb{N}}$ converges almost everywhere to $(\rho,u,m)$, we deduce that up to a subsequence:
\begin{itemize} 
\item $(\rho_n,\frac{1}{\rho_n})_{n\in\mathbb{N}}$ converges strongly to $(\rho,\frac{1}{\rho})$ in $L^p_{loc}(\R^+\times\R)$ for any $1\leq p<+\infty$.
\item For any $T>0$, $(\rho_n)_{n\in\mathbb{N}}$ converges strongly to $\rho$ in $C([0,T],W^{-\e,p}_{loc}(\R))$ for any $1\leq p<+\infty$ and $\e>0$.
\item $(u_n)_{n\in\mathbb{N}}$, $(\sqrt{\rho_n}u_n)_{n\in\mathbb{N}}$, $(m_n)_{n\in\mathbb{N}}$ converge strongly in $L^{2+\e_2}_{loc}(\R^+\times\R)$ for $\e_2>0$ sufficiently small.
\end{itemize}
In particular since $(\sqrt{\rho_n}v_n)_{n\in\mathbb{N}}$ is uniformly bounded in $L^\infty(\R^+,L^2(\R))$, we deduce that up to a subsequence $(\sqrt{\rho_n}v_n)_{n\in\mathbb{N}}$ converges weakly in $L^\infty(\R^+,L^2(\R))$ to $m_3$. Since $\sqrt{\rho_n}v_n=\sqrt{\rho_n}u_n+\p_x\va_{2,n}(\rho_n)$ with $\va_{2,n}'(\rho_n)=\frac{\mu_n(\rho_n)}{\rho_n^{\frac{3}{2}}}$, we have seen using the previous informations that $(\sqrt{\rho_n}v_n)_{n\in\mathbb{N}}$ converges in the sense of the distributions to $\sqrt{\rho}u+\p_x\va_2(\rho)=m^2$. It implies that $m_3=m^2$. Indeed it suffices to prove that $(\rho_n^{\alpha-\frac{1}{2}})_{n\in\mathbb{N}}$ and $(\frac{1}{n} \rho_n^{\theta-\frac{1}{2}})_{n\in\mathbb{N}}$ converge respectively to $\rho^{\alpha-\frac{1}{2}}$ and $0$. It is again a direct consequence of (\ref{superimpod12}), dominated convergence theorem and of the fact that up to a subsequence $(\rho_n)_{n\in\mathbb{N}}$ converges almost everywhere to $\rho$. Now using the Fatou lemma for weak convergence we deduce that there exists $C>0$ such that:
\begin{equation}
\|m^{2}\|_{L^\infty(\R^+,L^2(\R))}\leq C.
\label{estimm2}
\end{equation}
Since $(\rho_n,u_n)_{n\in\mathbb{N}}$ is a regular solution, we have for any test function $\varphi\in C^\infty_0([0,+\infty[\times\R)$ after integration by parts:
\begin{equation}
\begin{aligned}
&\int^{+\infty}_0\int_{\R}\big(\rho_n u_n(s,x)\p_s\varphi(s,x)+\rho_nu_n^2(s,x)\p_x\varphi(s,x)+\mu_n(\rho_n)\,u_n(s,x)\p_{xx}\varphi(s,x)\\
&+\p_x(\mu_n(\rho_n))\, u_n(s,x)\p_x\varphi(s,x)+a \rho_n^\gamma(s,x)\p_x\varphi(s,x)\big) ds\,dx+\int_{\R}\rho_nu_n(0,x)\varphi(0,x) dx=0.
\end{aligned}
\label{faible}
\end{equation}
From (\ref{superimpod12}), we know that $(\p_x(\mu_n(\rho_n)))_{n\in\mathbb{N}}$ is uniformly bounded in $L^2(\R^+,L^2(\R))$, then up to a subsequence  $\p_x(\mu_n(\rho_n))$ converges weakly in $L^2(\R^+,L^2(\R))$ to $w$. Since $\mu_n(\rho_n)$ converges to $\mu(\rho)$ in $L^p_{loc}(\R^+\times\R)$ for any $1\leq p<+\infty$ using dominated convergence, we deduce that $w=\p_x\mu(\rho)$ in the sense of the distributions. Now since we have seen that $u_n\p_x\varphi$ converges strongly in $L^2(\R^+,L^2(\R))$ up to a subsequence to $u\p_x\varphi$ (we use the fact here that $\p_x\va$ has compact support), it yields up to a subsequence:
\begin{equation}
\lim_{n\rightarrow+\infty}\int^{+\infty}_0\int_{\R} \p_x(\mu_n(\rho_n))\, u_n(s,x)\p_x\varphi(s,x) dx ds=
\int^{+\infty}_0\int_{\R} \p_x(\mu(\rho))\, u(s,x)\p_x\varphi(s,x) dx ds
\end{equation}
Proceeding similarly for the other terms we have proved that $(\rho,u)$ is a global weak solution in the sense of the definition \ref{def}. It is important in particular to observe that $m_n(0,\cdot)$ converges to $m_0$ in the sense of the measure and we have then ATTENTION ICI $m_0(x) dx$ EST UNE MESURE:
\begin{equation}
\lim_{n\rightarrow+\infty}\int_{\R}\rho_nu_n(0,x)\varphi(0,x) dx=
\int_{\R}m_0(x)\varphi(0,x) dx.
\end{equation}
We proceed also similarly for the mass equation. From (\ref{superimpod12}) $(v_n)_{n\in\mathbb{N}}$ is uniformly bounded in $L^\infty(\R^+,L^2(\R))$ and converges weakly * to $v$ in  $L^\infty(\R^+,L^2(\R))$. We verify in particular easily that $m^2=\sqrt{\rho}v$. From the Fatou properties for
the  weak convergence we deduce that for any $t>0$:
\begin{equation}
\begin{cases}
\begin{aligned}
%&\|\rho(t,\cdot)-\bar{\rho}\|_{L^2(\R)}\leq C\\
&\|v\|_{L^\infty([0,t],L^2(\R))}\leq C\\
%&\|\sqrt{\rho_n} u_n(t,\cdot)\|_{L^2(\R)}\leq C(1+\frac{1}{\sqrt{t}})\\
&\|(\sqrt{s}1_{\{s\leq T_\beta\}}+1_{\{s\geq T_\beta\}})\p_x u\|_{L^2_t(L^2(\R))}\leq C\\
&\|\p_x\rho^{\alpha-\frac{1}{2}}\|_{L^2([0,T_\beta],L^2(\R))+L^\infty([T_\beta,+\infty[,L^2(\R))}\leq C\\
&\|\p_x\rho^{\frac{1}{2}(\gamma+\alpha-1)}\|_{L^2_t(L^2(\R))}\leq C\\
&\|u\|_{L^{p(s)-\e}([0,T_\beta], H^s(\R))}\leq C.
%&u\in L^{p(s)-\e}([0,T_\beta], H^s(\R))\\
%&\|(\frac{1}{\rho},\rho)(t,\cdot)\|_{L^\infty(\R)}\leq C.
%&\|\sqrt{\rho}\p_x\va(\rho(t,\cdot))\|_{L^2(\R)}\leq C\\
%&\|\rho_n(t,\cdot)-\bar{\rho}\|_{L^\gamma_2(\R)}\leq C\\
%&\|(\frac{1}{\rho_n},\rho_n)(t,\cdot)\|_{L^\infty(\R)}\leq C(t)\\
%&\|\rho_n u_n(t,\cdot)\|_{L^1(\R)}+ \|\rho_n v_n(t,\cdot)\|_{L^1(\R)}\leq C(t). 
\end{aligned}
\end{cases}
\label{superimpod12abis}
\end{equation} 
%ESTIMATIONS, FATOU!!
\subsubsection*{Compactness when $\alpha\geq 1$}
This case is more delicate essentially because $(\frac{1}{\rho_n})_{n\in\mathbb{N}}$ is a priori not uniformly bounded in $L^\infty(\R^+,L^\infty(\R))$. We only know that there exists $T_\beta>0$ independent on $n$ such that $(\frac{1}{\rho_n})_{n\in\mathbb{N}}$ is uniformly bounded in $L^\infty([0,T_\beta],L^\infty(\R))$. In particular we have the same estimates on $(\rho_n, u_n,v_n)$ and $(\rho,u,v)$ on the time interval $[0,T_\beta]$ as previously ($(\rho,u,v)$ are defined in a similar way). The main difficulty now is to prove that $\rho_n$ and a certain momentum $m_n=\rho_n^{\beta_1}u_n$ converge almost everywhere up to a subsequence to $\rho$ and $m$.\\
\\
\textbf{Convergence almost everywhere of $(\rho_n)_{n\in\mathbb{N}}$}
\\
\\
From (\ref{superimpod1}) and (\ref{normeLin}), we have for any $\varphi\in C^\infty_0(\R)$ with $\mbox{supp}\varphi\subset]-M,M[$, there exists $C(t_1)>0$, $C>0$ independent on $n$ such that for $t_1\in]0,T_\beta]$:
\begin{equation}
\begin{aligned}
\|\p_x(\varphi \rho_n^\alpha)\|_{L^\infty([t_1,+\infty[,L^1(\R))}&\leq C \|\varphi\sqrt{\rho_n}\| _{L^\infty([t_1,+\infty[, L^2(\R))}\|\p_x\rho_n^{\alpha-\frac{1}{2}}\|_{L^\infty([t_1,+\infty[,L^2(\R))}\\
&\hspace{4cm}+\|\rho_n^\alpha \p_x\varphi\|_{L^\infty([t_1,+\infty[,L^1(\R))}\\
&\leq( \|\varphi\|_{L^2(\R)}+\|\p_x \varphi\|_{L^1(\R)}) C(M,t_1).
\end{aligned}
\end{equation}
with $C(M,t_1)$ depending only on $t_1,M$. Next we have:
$$\p_t(\varphi\rho_n^\alpha) =(1-\alpha)\varphi \rho_n^\alpha\p_x u_n- \p_x(\varphi \rho_n^\alpha u_n)+\p_x\varphi \,\rho_n^\alpha u_n.$$
From (\ref{normeLin}) and (\ref{superimpod1}), we deduce  that $\varphi \rho_n^\alpha\p_x u_n$ is uniformly bounded in $L^2([t_1,+\infty[,L^2(\R))$, $ \p_x(\varphi \rho_n^\alpha u_n)$ is uniformly bounded (since $\alpha>\frac{1}{2}$) in $L^\infty([t_1,+\infty[,W^{-1,2}(\R)))$, $\p_x\varphi \,\rho_n^\alpha u_n$ is uniformly bounded in $L^\infty([t_1,+\infty[,L^2(\R)))$. It implies that $(\p_t( \varphi\rho_n^\alpha))_{n\in\mathbb{N}}$ is uniformly bounded in $L^2_{loc}([t_1,+\infty[,W^{-1,2}(\R)))$.
From the Lemma \ref{Aubin}, we deduce that $(\rho_n^\alpha)_{n\in\mathbb{N}}$ converges strongly up to a subsequence in $C([t_1,T],L^{\beta_1}_{loc}(\R)))$ for any $1\leq\beta_1<+\infty$ and $T>t_1$. Up to a subsequence it implies that $(\rho_n^\alpha)_{n\in\mathbb{N}}$ converges to $\rho_\alpha$ almost everywhere on $[t_1,+\infty[\times\R$. In particular up to a subsequence $(\rho_n)_{n\in\mathbb{N}}$ converges almost everywhere to $\rho$ on $[t_1,+\infty[\times\R$ with $\rho=\rho_\alpha^{\frac{1}{\alpha}}$ if $\rho_\alpha>0$ and $\rho=0$ if $\rho_\alpha=0$.\\
Now using the same arguments as the proof for $\alpha<\frac{1}{2}$ on $]0,t_1]$ with $t_1\in ]0,T_\beta]$, we obtain that $(\rho_n)_{n\in\mathbb{N}}$ converges almost everywhere on $]0,t_1]\times\R$ to $\rho$ up to a subsequence. Indeed we can use similar ideas since we know from (\ref{impbdu}) that $(\rho_n,\frac{1}{\rho_n})_{n\in\mathbb{N}}$ is uniformly bounded in $L^\infty([0,T_\beta],L^\infty(\R))$. We have finally proved that $\rho_n$ converges almost everywhere to $\rho$ on $]0,+\infty[\times\R$. From (\ref{normeLin}) we deduce that for any $t>0$ we have for $C>0$:
\begin{equation}
\|\rho(t,\cdot)\|_{L^\infty(\R)}\leq C.
\label{normLinfin}
\end{equation}
In addition it implies using the dominated convergence and (\ref{normeLin}) that $(\rho_n^\gamma)_{n\in\mathbb{N}}$ converges up to a subsequence to $\rho^\gamma$ in $L^1_{loc}(\R^+\times\R)$.
\\
\\
\textbf{Convergence almost everywhere of $(\rho_n^\alpha u_n)_{n\in\mathbb{N}}$ when $\alpha\geq 1$}
\\
\\
We set $m_n=\rho_n^\alpha u_n$, we deduce from (\ref{superimpod1})  and (\ref{normeLin})%and for any $\varphi\in C^\infty_0(\R)$ (with $\mbox{supp}\varphi\subset[-M;M]$) , $t\geq t_1$ 
that there exists $C(t_1)>0$ such that for $t_1\in]0,T_\beta]$ and $t>t_1$:%$C(t,t_1,M)>0$:
%\begin{equation}
%\begin{aligned}
%&\|\p_x (\rho_n^\alpha u_n\varphi)\|_{L^2([t_1,t],L^1(\R))}\leq C\big( \|\rho_n^{\frac{\alpha}{2}}\p_x u_n\|_{L^2([t_1,t],L^2(\R))}\|\varphi \rho_n^{\frac{\alpha}{2}}\|_{L^2([t_1,t],L^2(\R))}\\
%&+\|\p_x(\rho_n^{\alpha-\frac{1}{2}})\|_{L^\infty([t_1,t],L^2(\R))}\|\sqrt{\rho_n}u_n\|_{L^\infty([t_1,t],L^2(\R))}(t-t_1)\|\varphi\|_{L^\infty(\R)}\\
%&+\|\sqrt{\rho_n}u_n\|_{L^\infty([t_1,t],L^2(\R))}\|\rho_n\|_{L^\infty(\R^+,L^\infty(\R)}^{\alpha-\frac{1}{2}}\|\p_x \varphi\|_{L^\infty(\R)}\\
%&\leq C(t,t_1, M) ( \|\varphi\|_{L^\infty(\R)}+\|\p_x \varphi\|_{L^\infty(\R)}) .
%\end{aligned}
%\label{Aubaa}
%\end{equation}
\begin{equation}
\begin{aligned}
&\|\p_x (\rho_n^\alpha u_n)\|_{L^2([t_1,t],L^2(\R))+L^\infty([t_1,t],L^1(\R))}\leq C\big( \|\rho_n^{\frac{\alpha}{2}}\p_x u_n\|_{L^2([t_1,t],L^2(\R))}\| \rho_n^{\frac{\alpha}{2}}\|_{L^\infty([t_1,t],L^\infty(\R))}\\
&+\|\p_x(\rho_n^{\alpha-\frac{1}{2}})\|_{L^\infty([t_1,t],L^2(\R))}\|\sqrt{\rho_n}u_n\|_{L^\infty([t_1,t],L^2(\R))}\big)\leq C (t_1),
\end{aligned}
\label{Aubaa}
\end{equation}
and:
\begin{equation}
\begin{aligned}
\|\rho_n^\alpha u_n\|_{L^\infty([t_1,t],L^2(\R))}&\leq \|
\sqrt{\rho_n}u_n\|_{L^\infty([t_1,t],L^2(\R))}\|\rho_n^{\alpha-\frac{1}{2}}\|_{L^\infty([t_1,t],L^\infty(\R))}\leq C (t_1).
\end{aligned}
\label{Aubaa1}
\end{equation}
Next we have:
$$
\begin{aligned}
&\p_t(\rho_n^\alpha u_n)=\rho_n^{\alpha-1}\p_t(\rho_n u_n)+\rho_n u_n\p_t\rho_n^{\alpha-1}\\
&=\p_x(\mu_n(\rho_n)\rho_n^{\alpha-1}\p_x u_n)-\p_x(\rho_n^{\alpha-1})\mu_n(\rho_n)\p_x u_n-\p_x(\rho_n^\alpha u_n^2)+\p_x(\rho_n^{\alpha-1})\rho_n u_n^2\\
&-\frac{a \gamma}{\gamma+\alpha-1}\p_x\rho_n^{\gamma+\alpha-1}-(\alpha-1)\rho_n^\alpha u_n\p_xu_n-\rho_n u_n^2 \p_x(\rho_n^{\alpha-1})\\
&=\p_x(\mu_n(\rho_n)\rho_n^{\frac{\alpha}{2}-1}\rho_n^{\frac{\alpha}{2}}\p_x u_n)-(\alpha-1)\frac{\mu_n(\rho_n)}{\rho_n^{\frac{3}{2}}}\p_x(\rho_n)\sqrt{\mu_n(\rho_n)}\p_x u_n \frac{\rho_n^{\alpha-\frac{1}{2}}}{\sqrt{\mu_n(\rho_n)}}-\p_x(\rho_n u_n^2 \rho_n^{\alpha-1})\\
&-(\alpha-1)\rho_n^{\frac{\alpha}{2}} u_n \rho_n^{\frac{\alpha}{2}} \p_xu_n-\frac{a \gamma}{\gamma+\alpha-1}\p_x\rho_n^{\gamma+\alpha-1}.
\end{aligned}
$$
Here since $\alpha\geq 1$ it is important to point out that from  (\ref{normeLin}), $\frac{\rho_n^{\alpha-\frac{1}{2}}}{\sqrt{\mu_n(\rho_n)}}$ is uniformly bounded in $L^\infty(\R^+,L^\infty(\R))$. 
From (\ref{superimpod1})  and (\ref{normeLin}) it yields that for $C>0, C(t_1)>0$ independent on $n$:
\begin{equation}
\begin{aligned}
&\|\p_t(\rho_n^\alpha u_n)\|_{L^2([t_1,t],H^{-1}(\R))+L^\infty([t_1,t],W^{-1,1}(\R))+L^2([t_1,t],L^{1}(\R)+L^\infty([t_1,t],W^{-1,\infty}(\R)) }\\
&\leq C\big( \|\mu_n(\rho_n)\rho_n^{\frac{\alpha}{2}-1}\|_{L^\infty ([t_1,t],L^{\infty}(\R))}\|\rho_n^{\frac{\alpha}{2}}\p_x u_n\|_{L^2([t_1,t],L^{2}(\R))}\\
&+\|\frac{\mu_n(\rho_n)}{\rho_n^{\frac{3}{2}}}\p_x(\rho_n)\|_{L^\infty([t_1,t],L^{2}(\R))} \|\sqrt{\mu_n(\rho_n)}\p_x u_n \|_{L^2([t_1,t],L^{2}(\R))} \|\frac{\rho_n^{\alpha-\frac{1}{2}}}{\sqrt{\mu_n(\rho_n)}}\|_{L^\infty ([t_1,t],L^{\infty}(\R))}\\
&+\|\rho_n u_n^2 \|_{L^\infty([t_1,t],L^{1}(\R))}\| \rho_n^{\alpha-1}\|_{L^\infty ([t_1,t],L^{\infty}(\R))}+\|\rho_n^{\gamma+\alpha-1}\|_{L^\infty ([t_1,t],L^{\infty}(\R))}\\
&+\|\rho_n^{\frac{1}{2}} u_n\|_{L^\infty ([t_1,t],L^{2}(\R))} \|\rho_n^{\frac{1}{2}(\alpha-1)}\|_{L^\infty ([t_1,t],L^{\infty}(\R))}  \| \rho_n^{\frac{\alpha}{2}} \p_x u_n\|_{L^2([t_1,t],L^{2}(\R))} \\
&\leq C(t_1).%+\|\rho_n^{\gamma+\alpha-1}\|_{L^\infty ([t_1,t],L^{\infty}(\R))} .
\end{aligned}
\label{Aubba1}
\end{equation}
We have seen that $(\rho_n^\alpha u_n)_{n\in\mathbb{N}^*}$ is uniformly bounded in $L^2_{loc}([t_1,+\infty[,W^{1,1}_{loc}(\R))$ and $(\p_t(\rho_n^\alpha u_n))_{n\in\mathbb{N}^*}$ is uniformly bounded in $L^2_{loc}([t_1,+\infty[,W^{-s,\infty}_{loc}(\R))$ for $s>0$ large enough by Sobolev embedding.\\
From (\ref{Aubaa}), (\ref{Aubaa1}), (\ref{Aubba1}), the diagonal process and Lemma \ref{Aubin}, we deduce that $(\rho_n^\alpha u_n)_{n\in\mathbb{N}}$ converges strongly up to a subsequence to $m_\alpha$
in $L^2_{loc}([t_1,+\infty[\times\R)$. In addition up to a subsequence, it implies that $\rho_n^\alpha u_n$ converges almost everywhere on $[t_1,+\infty[\times\R$ to $m_\alpha$.\\
\\
Note that we can already define the velocity $u$ with $u(t,x)=\frac{m_\alpha(t,x)}{\rho^\alpha(t,x)}$
on the set $\{(t,x)\in\R^+\times\R;\,\,\rho(t,x)>0\}$.
Let us verify now that $m_\alpha=0$ almost everywhere on the set $\{(t,x)\in\R^+\times\R;\,\,\rho(t,x)=0\}$. Indeed we have for any $t>0$ and using the Fatou Lemma and (\ref{superimpod1}):
$$
\begin{aligned}
%\int_{\{\rho(t,x)\leq 1\}}\frac{m^2_\alpha(t,x)}{\rho^{2\alpha-1}(t,x)} dx&=
\int_{\{\rho(t,x)\leq 1\}}\lim \inf_{n\rightarrow+\infty}\frac{(\rho_n^\alpha u_n(t,x))^2}{\rho_n^{2\alpha-1}(t,x)} dx
&\leq \lim \inf _{n\rightarrow+\infty}\int_{\R}\rho_n |u_n|^2(t,x) dx\\
&\leq C(1+\frac{1}{t}),
\end{aligned}
$$
with $C>0$ independent on $n$ (we mention that for any $n$ we have always $\rho_n(t,\cdot)>0$ for any $t\in\R^+$). It implies since $2\alpha-1\geq 0$ that $m_\alpha=0$ almost everywhere on the set $\{(t,x)\in\R^+\times\R;\,\,\rho(t,x)=0\}$. We have then:
\begin{itemize}
\item $m_\alpha(t,x)=0$ almost everywhere on $\{(t,x)\in\R^+\times\R;\,\,\rho(t,x)=0\}$.
\item $u(t,x)=\frac{m_\alpha(t,x)}{\rho^\alpha(t,x)}$ on $\{(t,x)\in\R^+\times\R;\,\,\rho(t,x)>0\}$ and $u(t,x)=0$ on $\{(t,x)\in\R^+\times\R;\,\,\rho(t,x)=0\}$.
\end{itemize}
We can observe that $u$ is not uniquely defined on $\{(t,x)\in\R^+\times\R;\,\,\rho(t,x)=0\}$. Now from (\ref{gainfin1}) and the Fatou lemma we deduce that there exists $C(t,t_1)>0$ such that for any $t\geq t_1$ we have:
\begin{equation}
\|\rho^{\frac{1}{4}}u\|_{L^\infty([t_1,t],L^4(\R))}\leq C(t_1,t).
\label{gainfin2}
\end{equation}
Similarly from (\ref{gainL1infini}) we have for any $t>0$:
\begin{equation}
\|\rho u(t,\cdot)\|_{L^1(\R)}\leq C(t).
\label{infoL1}
\end{equation}
We are going now to prove the strong $L^2_{loc}$ convergence of $(\sqrt{\rho_n}u_n)_{n\in\mathbb{N}}$.
\begin{lemme}
We have for any $T>0$ and any compact $K$ of $\R$:
\begin{equation}
\int^T_0\int_{K}|\sqrt{\rho_n}(t,x)u_n(t,x)-\sqrt{\rho}(t,x)u(t,x)|^2 dt dx\rightarrow _{n\rightarrow+\infty}0.
\label{gainucru}
\end{equation}
This limit is true up to a subsequence.
\end{lemme}
{\bf Proof:} Following the proof of the case $\alpha<\frac{1}{2}$ (since $(\rho_n,\frac{1}{\rho_n})_{n\in\mathbb{N}}$ is uniformly bounded on $[0,T_\beta]$), we know that up to a subsequence we have for any compact $K$ (we refer in particular to (\ref{impbdu})):
\begin{equation}
\int^{T_{\beta}}_0\int_{K}|\sqrt{\rho_n}(t,x)u_n(t,x)-\sqrt{\rho}(t,x)u(t,x)|^2 dt dx\rightarrow _{n\rightarrow+\infty}0.
\label{gainucru1}
\end{equation}
Let us deal now with the case $T>T_{\beta}$. We have seen that $(\sqrt{\rho_n}u_n)_{n\in\mathbb{N}}$ converges almost everywhere to $\sqrt{\rho}u$ on the set $\{(t,x)\in [t_1,+\infty[\times\R;\,\,\rho(t,x)>0\}$ with $t_1\in]0,T_\beta]$. It implies that we have for any $T>0$ using dominated convergence, (\ref{normeLin}) and taking $M>0$ with $K_1=K\cap \{(t,x)\in [t_1,T]\times\R;\,\,\rho(t,x)>0\}$:
\begin{equation}
\begin{aligned}
&\int^{T}_{t_1}\int_{K_1}|\sqrt{\rho_n}(t,x)u_n(t,x)1_{\{|u_n|(t,x)\leq M\}}-\sqrt{\rho}(t,x)u(t,x)1_{\{|u|(t,x)\leq M\}}|^2 dt dx \rightarrow _{n\rightarrow+\infty}0.
\end{aligned}
\label{estiu1}
\end{equation}
We have now using (\ref{gainfin1}), (\ref{gainfin2}) and (\ref{normeLin}) that there exists $C>0$ independent on $n$ such that:
\begin{equation}
\begin{aligned}
&\int^{T}_{t_1}\int_{K_1}%\cap \{(t,x)\in [t_1,T]\times\R;\,\,\rho(t,x)>0\}}
|\sqrt{\rho_n}(t,x)u_n(t,x)|^2 1_{\{|u_n|(t,x)> M\}} dt dx\\
&\leq \frac{C}{M^2}\int^{T}_{t_1}\int_{K_1%\cap \{(t,x)\in [t_1,T]\times\R;\,\,\rho(t,x)>0\}
}\rho_n |u_n|^4(t,x)dt dx\leq  \frac{C(T,t_1)}{M^2}\rightarrow_{M\rightarrow+\infty}0,
\end{aligned}
\label{estiu2}
\end{equation}
similarly we have:
\begin{equation}
\begin{aligned}
&\int^{T}_{t_1}\int_{K_1%\cap \{(t,x)\in [t_1,T]\times\R;\,\,\rho(t,x)>0\}
}|\sqrt{\rho}(t,x)u(t,x)|^2 1_{\{|u|(t,x)> M\}} dt dx\leq  \frac{C(T,t_1)}{M^2}\rightarrow_{M\rightarrow+\infty}0.
\end{aligned}
\label{estiu3}
\end{equation}
On the set $\{(t,x)\in [t_1,+\infty[\times\R;\,\,\rho(t,x)=0\}$, we have for any $M>0$:
$$\sqrt{\rho_n}|u_n(t,x)|1_{\{|u_n|(t,x)\leq M\}}\leq M\sqrt{\rho_n}(t,x)\rightarrow_{n\rightarrow+\infty}0=\sqrt{\rho}|u|(t,x).$$
From convergence dominated, (\ref{normeLin}),  (\ref{gainfin1}), (\ref{gainfin2}) and Tchebytchev lemma we deduce that:
\begin{equation}
\begin{aligned}
&\int^{T}_{t_1}\int_{K\cap \{(t,x)\in [t_1,T]\times\R;\,\,\rho(t,x)=0\}}|\sqrt{\rho_n}(t,x)u_n(t,x)|^2 dt dx\rightarrow_{n\rightarrow+\infty}0.
\end{aligned}
\label{estiu4}
\end{equation}
From (\ref{gainucru1}), (\ref{estiu1}), (\ref{estiu2}), (\ref{estiu3}), (\ref{estiu4}) and the fact that $\sqrt{\rho}u=0$ on $\{(t,x)\in[t_1,+\infty[\times\R,\;\rho(t,x)=0\}$ we deduce (\ref{gainucru}).  {\hfill $\Box$}\\
\\
Let us prove now that the sequel $(\rho_nu_n)_{n\in\mathbb{N}}$ converges in the sense of the distribution to $\rho u$. Indeed we have for any $T>0$, any $K$ a compact of $\R$ and $C(T,K)>0$:
\begin{equation}
\begin{aligned}
&\|\rho_n u_n-\rho u\|_{L^1([0,T],L^1(K))}\leq C(T,K)\big(\|\sqrt{\rho_n}\|_{L^\infty([0,T],L^\infty(K))}\|\sqrt{\rho_n} u_n-\sqrt{\rho} u\|_{L^2([0,T],L^2(K))}\\
&\hspace{4cm}+\|\sqrt{\rho_n}-\sqrt{\rho}\|_{L^2([0,T],L^2(K))}\|\sqrt{\rho}u\|_{L^2([0,T],L^2(K))}\big).
\end{aligned}
\end{equation}
From (\ref{gainucru}), (\ref{normeLin}), convergence dominated and the fact that $\sqrt{\rho_n}$ converges almost everywhere to $\sqrt{\rho}$ we deduce that 
$\rho_n u_n$ converges to $\rho u$ in $L^1_{loc}(\R^+\times\R)$. Now from (\ref{normeLin}), (\ref{gainfin1}), (\ref{gainfin2}), (\ref{impbdu1}) and by interpolation we deduce that $(\rho_n u_n)_{n\in\mathbb{N}}$ converges to $\rho u$ in $L^2_{loc}(\R^+\times\R)$.
\\
Let us deal now with the term of viscosity, we have then for any $\varphi\in C^\infty_0(\R^+\times\R)$:
$$
\begin{aligned}
&\int_0^\infty \int_{\R}\mu_n(\rho_n)\p_x u_n(t,x) \p_x \varphi(t,x) dx dt=-\mu \int_0^\infty\int_{\R}\rho_n^{\alpha-\frac{1}{2}} \sqrt{\rho_n}u_n(t,x) \p_{xx} \varphi(t,x) dx dt\\
&-\frac{1}{n}\int_0^\infty\int_{\R}\rho_n^{\theta}u_n(t,x) \p_{xx} \varphi(t,x) dx dt- \int_0^\infty\int_{\R}\p_x( \mu_n(\rho_n))(t,x) u_n(t,x) \p_x \varphi(t,x) dx dt.
\end{aligned}
$$
We have only to consider the case of the interval $[t_1,+\infty[$ with $t_1\in]0,T_\beta]$ since the procedure follows the proof of the case $\alpha<\frac{1}{2}$ on $[0,T_\beta]$. We know from (\ref{superimpod1}) and (\ref{gainfin1}) that we have for any $t\in [t_1,+\infty[$:
\begin{equation}
\begin{cases}
\begin{aligned}
&\|\p_x \rho_n^{\alpha-\frac{1}{2}}(t,\cdot)\|_{L^2(\R)}+\frac{1}{n}\|\rho_n^{\theta-\frac{3}{2}}\p_x\rho_n(t,\cdot)\|_{L^2(\R)}+\frac{1}{\sqrt{n}}\|\rho_n^{\frac{\alpha+\theta-3}{2}}\p_x\rho_n(t,\cdot)\|_{L^2(\R)}\leq C(t_1)\\
&\|\rho_n^{\frac{1}{4}}u_n(t,\cdot)\|_{L^4(\R)}\leq C(t,t_1).
\end{aligned}
\end{cases}
\label{309}
\end{equation}
We have now:
$$u_n\p_x(\mu_n(\rho_n))=\frac{\mu\alpha}{\alpha-\frac{1}{2}}\sqrt{\rho_n}u_n \p_x\rho_n^{\alpha-\frac{1}{2}}+\frac{\theta}{n}\rho_n^{\frac{1}{4}}u_n\,\rho_n^{\theta-\frac{5}{4}}\p_x\rho_n.$$
From (\ref{309}), we have setting $\theta_1\in (0,1)$ such that $\theta-\frac{5}{4}=\theta_1(\alpha-\frac{3}{2})+(1-\theta_1)(\theta-\frac{3}{2})$ (this is true because we have $\alpha\geq 1$ and $\theta\leq\frac{1}{2}$) and for any $t\geq t_1$:
\begin{equation}
\begin{aligned}
&\int_{\R}\rho_n^{2\theta-\frac{5}{2}}(t,x)|\p_x\rho_n(t,x)|^2 dx\\
&\leq (\int_{\R}\rho_n^{2(\theta-\frac{3}{2})}(t,x)|\p_x\rho_n(t,x)|^2 dx)^{(1-\theta_1)}(\int_{\R}\rho_n^{2(\alpha-\frac{3}{2})}(t,x)|\p_x\rho_n(t,x)|^2 dx)^{\theta_1}\\
&\leq n^{2(1-\theta_1)}C(t_1).
\end{aligned}
\label{utiler}
\end{equation}
Using H\"older's inequality, (\ref{309}) and (\ref{utiler}), we deduce that we have:
\begin{equation}
\begin{aligned}
&|\frac{\theta}{n}\int_{t_1}^{+\infty}\int_{\R} \rho_n^{\frac{1}{4}}u_n\,\rho_n^{\theta-\frac{5}{4}}\p_x\rho_n(t,x) \p_x\varphi(t,x) dt dx|\\
&\leq \frac{\theta}{n}\|\rho_n^{\frac{1}{4}}u_n |\p_x\varphi|^{\frac{1}{2}}\|_{L^\infty([t_1,+\infty[,L^4(\R))}\|\rho_n^{\theta-\frac{5}{4}}\p_x\rho_n |\p_x\varphi|^{\frac{1}{2}}\|_{L^1([t_1,+\infty[,L^{\frac{4}{3}}(\R))}\\
&\leq \frac{\theta}{n}\|\rho_n^{\frac{1}{4}}u_n\|_{L^\infty([t_1,+\infty[,L^4(\R))}
 \|\p_x\varphi\|_{L^\infty([t_1,+\infty[,L^\infty(\R))}^{\frac{1}{2}}\||\p_x\varphi|^{\frac{1}{2}}\|_{L^1([t_1,+\infty[,L^{4}(\R))}\\
 &\hspace{7cm}\times
\|\rho_n^{\theta-\frac{5}{4}}\p_x\rho_n \|_{L^\infty([t_1,+\infty[,L^{2}(\R))}\\
&\leq \frac{\theta}{n^{\theta_1}}C(t_1)
 \|\p_x\varphi\|_{L^\infty([t_1,+\infty[,L^\infty(\R))}^{\frac{1}{2}}\||\p_x\varphi|^{\frac{1}{2}}\|_{L^1([t_1,+\infty[,L^{4}(\R))}.
%&\leq \frac{\theta}{n}\|\rho_n^{\frac{1}{4}}u_n |\p_x\varphi|^{\frac{1}{2}}\|_{L^\infty([t_1,+\infty[,L^4(\R))}\|\rho_n^{\theta-\frac{5}{4}}\p_x\rho_n |\p_x\varphi|^{\frac{1}{2}}\|_{L^1([t_1,+\infty[,L^{\frac{4}{3}}(\R))}\\
%&\leq C(\mbox{supp}\varphi)(\|\p_x\varphi\|_{L^\infty(\R^+\times\R)}+1)\frac{\theta}{n}\|\rho_n^{\frac{1}{4}}u_n \|_{L^\infty([t_1,+\infty[,L^4(\R))}\\
%&\hspace{6cm}\times\|\rho_n^{\theta-\frac{5}{4}}\p_x\rho_n\|_{L^\infty ([t_1,+\infty[,L^{2}(\R))}\\
%&\leq C(\mbox{supp}\varphi)C(t_1)(\|\p_x\varphi\|_{L^\infty(\R^+\times\R)}+1)\frac{\theta}{n^{\theta_1}}.
\end{aligned}
\end{equation}
We deduce then that:
\begin{equation}
\frac{\theta}{n}\int_{t_1}^{+\infty}\int_{\R} \rho_n^{\frac{1}{4}}u_n\,\rho_n^{\theta-\frac{5}{4}}\p_x\rho_n(t,x) \p_x\varphi(t,x) dt dx\rightarrow_{n\rightarrow 0}0.
\label{techvisco}
\end{equation}
It is easy now to verify that:
\begin{equation}
\frac{\theta}{n}\int_{0}^{t_1}\int_{\R} \rho_n^{\frac{1}{4}}u_n\,\rho_n^{\theta-\frac{5}{4}}\p_x\rho_n(t,x) \p_x\varphi(t,x) dt dx\rightarrow_{n\rightarrow 0}0.
\label{techvisco1}
\end{equation}
We have finally proved that:
\begin{equation}
\frac{\theta}{n}\int_{0}^{+\infty}\int_{\R} \rho_n^{\frac{1}{4}}u_n\,\rho_n^{\theta-\frac{5}{4}}\p_x\rho_n(t,x) \p_x\varphi(t,x) dt dx\rightarrow_{n\rightarrow 0}0.
\label{techvisco2}
\end{equation}
Similarly if we fix now $\theta=\frac{1}{4}$, we obtain from (\ref{309}):
\begin{equation}
\frac{1}{n}\int_0^\infty\int_{\R}\rho_n^{\theta}u_n(t,x) \p_{xx} \varphi(t,x) dx dt\rightarrow_{n\rightarrow 0}0.
\label{techvisco1}
\end{equation}
We can prove easily now that up to a subsequence:
\begin{equation}
\begin{aligned}
&-\mu \int_0^\infty\int_{\R}\rho_n^{\alpha-\frac{1}{2}} \sqrt{\rho_n}u_n(t,x) \p_{xx} \varphi(t,x) dx dt
- \int_0^\infty\int_{\R}\frac{\mu\alpha}{\alpha-\frac{1}{2}}\sqrt{\rho_n}u_n \p_x\rho_n^{\alpha-\frac{1}{2}}\p_{x} \varphi(t,x) dx dt\\
&\rightarrow _{n\rightarrow+\infty}-\mu \int_0^\infty\int_{\R}\rho^{\alpha-\frac{1}{2}} \sqrt{\rho}u(t,x) \p_{xx} \varphi(t,x) dx dt\\
&\hspace{6cm}- \int_0^\infty\int_{\R}\frac{\mu\alpha}{\alpha-\frac{1}{2}}\sqrt{\rho}u \p_x\rho^{\alpha-\frac{1}{2}}\p_{x} \varphi(t,x) dx dt\\
\end{aligned}
\label{314}
\end{equation}
Indeed since $(\p_x\rho_n^{\alpha-\frac{1}{2}}\p_x\va)_{n\in\mathbb{N}}$ is uniformly bounded in $L^2([t_1,+\infty[\times\R)$, it implies that up to a subsequence $( \p_x\rho_n^{\alpha-\frac{1}{2}}\p_x\va)_{n\in\mathbb{N}}$ converges weakly in $L^2([t_1,+\infty[\times\R)$ to $w_1$. Since $(\rho_n^{\alpha-\frac{1}{2}})_{n\in\mathbb{N}}$ converges to $\rho^{\alpha-\frac{1}{2}}$ in $L^2_{loc}([t_1,+\infty[\times\R)$ via the dominated convergence, we deduce that $(\p_x\rho_n^{\alpha-\frac{1}{2}}\p_x\va)_{n\in\mathbb{N}}$ converges in the sense of the distribution to $w_1=\p_x\rho^{\alpha-\frac{1}{2}}\p_x\va $. Now using the fact that $(\sqrt{\rho_n}u_n)_{n\in\mathbb{N}}$ converges strongly to $\sqrt{\rho}u$  in $L^2_{loc}([t_1,+\infty[\times\R)$, we now obtain the convergence of the second term on the left hand side of (\ref{314}) (indeed on the time interval $[0,t_1]$ we proceed as in the previous section with $0<t_1<T_\beta$). We proceed similarly for the other term.\\
We have proved finally that $(\rho_n,u_n)_{n\in\mathbb{N}}$ converges in the sense of the distribution to a global weak solution $(\rho,u)$ of the system (\ref{1}). Finally we can prove the estimates (\ref{sestim1}) by using Fatou lemma in the framework of the weak convergence.
%Taking $s=\frac{1-\e}{2(2+\e)}$ with $\e>0$ sufficiently sm$all, we get that $u_n$ is uniformly bounded in $L^{2}_{loc}(H^{\frac{1-\e}{2(2+\e)}}(\R))$.
%Let us show the compactness result in the general case without restriction on $\alpha$, we simply assume that $\alpha>0$. A LA MELLET VASSEUR
\subsubsection*{Compactness when $\frac{1}{2}<\alpha<1$}
We proceed as in the previous section excepted that it suffices to use the Lemma \ref{Aubin} for the unknown $m_n=\rho_n u_n$. Indeed we observe that we have from (\ref{superimpod1})  and (\ref{normeLin}) with $t_1\in ]0,T_\beta]$:
\begin{equation}
\|\rho_n u_n\|_{L^\infty([t_1,+\infty[,L^2(\R))}\leq C(t_1),
\label{premiesti}
\end{equation}
and since:
$$\p_t(\rho_n u_n)=-\p_x(\rho_n u_n^2)+\p_x(\mu_n(\rho_n)\p_x u_n)-\p_x(a\rho_n^\gamma),$$
we have then:
\begin{equation}
\|\rho_n u_n\|_{L^\infty([t_1,+\infty[,W^{-1,1}(\R))+L^2([t_1,+\infty[,H^{-1}(\R))+L^\infty([t_1,+\infty[,W^{-1,+\infty}(\R))}\leq C(t_1).
\label{premiesti}
\end{equation}
We can then apply the Lemma \ref{Aubin} and conclude as in the previous section.
\subsubsection*{Continuity in time of the momentum $\rho u$ and $m^1$}
\begin{proposition}
For any $T>0$, $m$ and $\rho v$ belongs to $C([0,T],{\cal M}(\R)^{*})$.
\end{proposition}
{\bf Proof:} %VERIFIER CETTE PARTIE! \\
From (\ref{cru1}), the sequence $(\rho_n u_n)_{n\in\mathbb{N}}$  is uniformly bounded in $n$ on $[0,T]$ in ${\cal M}(\R)$ when $\|\rho_0^n u_0^n\|_{L^1}\leq M$, indeed there exists $M(T)>0$ such that for any $t\in[0,T]$ we have:
$$\|\rho_n u_n(t,\cdot)\|_{{\cal M}(\R)}\leq M(T).$$
We are now going to prove that for any $\varphi\in C_0(\R)$, $(\rho_n u_n(t,\cdot),\varphi)$ is uniformly continuous in $n$ on $[0,T]$. Indeed we have using the momentum equation of (\ref{1}), the Fubini theorem when $t>s\geq 0$ and (\ref{superimpod1}), (\ref{normeLin}), (\ref{impbdu}) and (\ref{gainfin1}) that for any $\va\in C_{0}^\infty(\R)$:
$$
\begin{aligned}
&|(\rho_n u_n(t,\cdot),\varphi)-(\rho_n u_n(s,\cdot),\varphi)|\\
&=|\int_{\R}\int^t_s [-\p_x(\rho_n u_n^2)(\tau,x)+\p_x(\mu_n(\rho_n)\p_x u_n)(\tau,x)-(\p_x P(\rho_n))(\tau,x)]\va(x)\,d\tau\, dx|\\
&=|\int^t_s\int_{\R} [\rho_n u_n^2(\tau,x)+\p_x(\mu_n(\rho_n)) u_n(\tau,x)+P(\rho_n)(\tau,x)]\p_x \va(x)\,d\tau\, dx\\
&\hspace{4cm}+ \int^t_s\int_{\R}\mu_n(\rho_n(\tau,x)) u_n(\tau,x)\p_{xx}\varphi(x) d\tau dx|
\end{aligned}
$$
$$
\begin{aligned}
&|(\rho_n u_n(t,\cdot),\varphi)-(\rho_n u_n(s,\cdot),\varphi)|\\
&\leq\int^t_s\big(\|\rho_n\|_{L^\infty(\R^+,L^\infty(\R)}\|1_{[0,T_\beta]}(\tau) u_n(\tau,\cdot)\|_{L^2(\R)}^2 +\|1_{t\geq T_{\beta}}(\tau) \sqrt{\rho_n}u_n(\tau,\cdot)\|_{L^2(\R)}^2 \big)\|\p_x\varphi\|_{L^\infty(\R)} d\tau\\
&+ \int^t_s (1_{[0,T_\beta]}(\tau)\tau  \frac{1}{\sqrt{\tau}}\|\sqrt{\rho_n}\p_x\varphi_n(\rho_n)(\tau,\cdot)\|_{L^2(\R)}\|\sqrt{\rho_n}\|_{L^\infty(\R^+,L^\infty(\R))}\frac{1}{\sqrt{\tau}}\| u_n(\tau,\cdot)\|_{L^2(\R)}\|\p_x\varphi\|_{L^\infty(\R)} \\
&+1_{\tau\geq T_\beta}(\tau)  \|\sqrt{\rho_n}u_n\|_{L^2(\R))}\|\sqrt{\rho_n}\p_x\varphi_n(\rho_n)\|_{L^2(\R)}\|\p_x\varphi\|_{L^\infty(\R)} \big) d\tau\\
&+\|P(\rho_n)\|_{L^\infty(\R^+,L^\infty(\R))}\|\p_x\varphi\|_{L^1(\R)}(t-s)\\
&+ \int^t_s \big( 1_{[0,T_\beta]}(\tau) \|\mu(\rho_n)\|_{L^\infty(\R^+,L^\infty(\R))}\|u_n(\tau,\cdot)\|_{L^2(\R)}\|\p_{xx}\varphi\|_{L^2(\R)} \\
&+1_{\tau\geq T_\beta}(\tau)\|\rho_n^{\frac{1}{4}}u_n(\tau,\cdot)\|_{L^4(\R)}\|\p_{xx}\varphi\|_{L^{\frac{4}{3}}(\R)}\|\frac{\mu_n(\rho_n)}{\rho_n^{\frac{1}{4}}}\|_{L^\infty(\R^+,L^\infty(\R))}\big)d\tau\\
&\leq C(|t-s|+|t-s|^{\alpha_1})(1+\|\p_x\varphi\|_{L^\infty(\R)}+\|\p_x\varphi\|_{L^1(\R)}+\|\p_{xx}\varphi\|_{L^\infty(\R)}+\|\p_{xx}\varphi\|_{L^1(\R)}),
%&=|\int_{\R^2}\big(\int_{\R^2}p_1^n(t,x,s,y)v^n(s,y)dy \big)\varphi(x) dx-\int_{\R^2}\big(\int_{\R^2}p_1^n(t,x,y,s)dx \big)v^n(s,y) \varphi(y) dy|\\
%&\leq \int_{\R^2}\int_{\R^2}p_1^n(t,x,s,y)|v^n(s,y)|\,| \varphi(x) - \varphi(y)| dx dy.
\end{aligned}
$$
with $\alpha_1>0$. We note that we have used the fact that $(u_n)_{n\in\mathbb{N}}$ is uniformly bounded in $L^{3-\e}([0,T_\beta],L^2(\R))$ with $\e>0$ small enough and the fact that we have choose $\theta=\frac{1}{4}$. It is also important to recall that $(\frac{1}{\rho_n})_{n\in\mathbb{N}}$ is uniformly bounded in $L^\infty(\R^+,L^\infty(\R))$ when $0<\alpha<\frac{1}{2}$.\\
By a density argument and the fact that $(\rho_nu_n)_{n\in\mathbb{N}} $ is uniformly bounded in $L^\infty([0,T],{\cal M}(\R))$ it implies that for $\varphi\in C_0(\R)$, $(\rho_n u_n(t,\cdot),\varphi)_{n\in\mathbb{N}}$ is uniformly continuous in $n$ on $[0,T]$. In other words the sequence $(\rho_n u_n)_{n\in\mathbb{N}}$ is equicontinuous on $[0,T]$ for the weak * topology on ${\cal M}(\R)$ when $\|\rho_0^n u_0^n\|_{L^1}\leq M$. Using the Ascoli theorem on the metric space $B(0,M(T))_{{\cal M}(\R)}$
endowed with the weak * topology (this is true because $C_0(\R)$ is separable), we deduce that up to a subsequence $\rho_n u_n (t,\cdot)$ converges uniformly to $m(t,\cdot)$ on $[0,T]$ in the weak * topology of ${\cal M}(\R)$. 
In particular it implies that $m$ belongs to $C([0,T],{\cal M}(\R)*)$ for any $T>0$. The proof is similar for $m^2=\rho v$.
{\hfill $\Box$}
%\subsubsection*{Particular case $\gamma=\alpha$}
%In this case since $\gamma>1$, we are in the case $\alpha>\frac{1}{2}$. In addition from (\ref{cru1}) we deduce that for any $t>0$ we have for any $n\in\mathbb{N}^*$ and $C>0$ independent on $n$:
%\begin{equation}
%\|\rho_n u_n(t,\cdot)\|_{L^1(\R)}+ \|\rho_n v_n(t,\cdot)\|_{L^1(\R)}\leq C (2\|\rho_0 v_0\|_{{\cal M}(\R)}+\|\va_1(\rho_0)\|_{BV(\R)})e^{Ct}.
%\label{L1}
%\end{equation}
%It implies in particular that $\rho_n$ is uniformly bounded in $L^\infty_{loc}(L^\infty(\R))$. From (\ref{gainu}), (\ref{BD}) and (\ref{dgainu}) we deduce that there exists $C(t)>0$ independent on $n$ such that for any $t>0$ we have:
%\begin{equation}
%\|\sqrt{\rho_n}u_n\|_{L^2([0,t],L^2(\R))}\leq C(t).
%\label{fingtech}
%\end{equation}
%Using (\ref{L1}) and (\ref{fingtech}) we can conclude as in the case $\alpha>\frac{1}{2}$ by using the same regularizing effects on $u_n$ and the same compactness arguments.
\section{Proof of the Corollary \ref{corbis}}%and \ref{cor1}}
\label{section3}
%\subsection{Proof of the Corollaries \ref{corbis}}
We proceed as in the proof of the theorem \ref{theo1}, we use in particular the same regularizing process. We have then a sequence $(\rho_n,u_n)_{n\in\mathbb{N}}$ solution of the system (\ref{1bis}) and we want to verify that up to a subsequence this sequence converges in the sense of the distributions to a weak solution of the system (\ref{1}). Compared with the theorem \ref{theo1}, we have no estimate as (\ref{BD}) since $v_0$ does not belong to $L^2(\R)$. However we have the classical energy estimate, indeed this is due to the fact that $u_n(0,\cdot)$ is uniformly bounded in $L^2(\R)$. Multiplying the momentum equation of (\ref{1bis}) and integrate over $(0,T)\times\R$ we have then:
\begin{equation}
\begin{aligned}
&\frac{1}{2}\int_{\R}\big(\rho_n(T,x)|u_n(T,x)|^2+(\Pi(\rho_n(T,x))-\Pi(\bar{\rho})) \big) dx\\
&+\int^T_0\int_{\R}\mu_n(\rho_n(s,x))(\p_x u_n(s,x))^2 dx ds\\
&\hspace{4cm}\leq \frac{1}{2}\int_{\R}\big(\rho_n(0,x)|u_n(0,x)|^2+(\Pi(\rho_n(0,x))-\Pi(\bar{\rho})) \big) dx.
\end{aligned}
\label{energiea}
\end{equation}
Now proceeding as in the proof of the theorem \ref{theo1}, we establish that there exists $T_\beta>0$, $C>0$ independent on $n$ such that:
\begin{equation}
\|(\frac{1}{\rho_n},\rho_n)\|_{L^\infty([0,T_\beta],L^\infty(\R))}\leq C.
\label{controlunif}
\end{equation}
From (\ref{cru1}), we deduce that it exists $C_1>0$ independent on $n$ such that:
\begin{equation}
\begin{aligned}
&\|\rho_nu_n\|_{L^\infty([0,T_\beta],L^1(\R))}+\|\rho_n v_n\|_{L^\infty([0,T_\beta],L^1(\R))} \\
&\hspace{1cm}\leq C(\|\rho_nu_n(0,\cdot)\|_{L^1(\R)}+ \|\rho_n v_n(0,\cdot)\|_{L^1(\R)})e^{C_1T_\beta^{\gamma-\alpha}}\leq C_1 \e_0 e^{C_1 T_\beta^{\gamma-\alpha}}.
\end{aligned}
\label{controlunif1}
\end{equation}
It implies in particular from (\ref{controlunif}) and (\ref{controlunif1}) that for any $t\in[0,T_\beta]$, there exists $C>0$ independent on $n$ such that:
\begin{equation}
\begin{aligned}
&\|\rho_n^\alpha(t,\cdot)\|_{BV(\R)}\leq C+C_1\e_0 e^{C_1 T_\beta^{\gamma-\alpha}}.
\end{aligned}
\label{controlunif2}
\end{equation}
From (\ref{controlunif}), (\ref{controlunif2}), composition theorem we deduce that for $C>0$ independent on $n$ and any $t\in[0,T_\beta]$:
\begin{equation}
\begin{aligned}
&\|\rho_n(t,\cdot)\|_{BV(\R)}\leq C(1+\e_0 e^{C T_\beta^{\gamma-\alpha}}).
\end{aligned}
\label{controlunif3}
\end{equation}
We wish now to prove that $\rho_n$ converges almost everywhere on $]0,T_\beta]\times\R$ up to a subsequence. We are going to use the lemma \ref{compac}. From (\ref{controlunif3}), it suffices only to prove the third assumption of the lemma \ref{compac}.
%To do this, we must prove each assumption of the lemma \ref{compac}. First we have since $\rho_n$ is regular on $[0,T_\beta]\times\R$ and for any increasing $-\infty<x_0<\cdots<x_N<+\infty$ and any $t\in[0,T_\beta]$ we have from (\ref{controlunif1}) and (\ref{controlunif}) for $C>0$:
%\begin{equation}
%\begin{aligned}
%\sum_{i=0}^{N-1} |\rho_n(t,x_{i+1})-\rho_n(t,x_{i})|&\leq \sum_{i=0}^{N-1}|\int_{x_i}^{x_{i+1}}\p_x\rho_n(t,x) dx|\\
%&\leq C \e_0 e^{C T_\beta^{\gamma-\alpha}}.
%\end{aligned}
%\label{BV1}
%\end{equation}
%From (\ref{controlunif}) , we have for any $t\in[0,T_\beta]$:
%\begin{equation}
%|\rho_n(t,0)|\leq C
%\label{BV2}
%\end{equation}
From the mass equation, we deduce now that for any compact $K\in\R$ we have from (\ref{energiea}), (\ref{controlunif}), Sobolev embedding that for any $(s,t)\in[0,T_\beta]^2$ and for $C>0$ independent on $n$ and depending on $T_\beta$:
\begin{equation}
\begin{aligned}
&\|\rho_n(t,\cdot)-\rho_n(s,\cdot)\|_{L^1(K)}\leq \int^t_s \|\p_x(\rho_n u_n)(s',\cdot)\|_{L^1(K)} ds'\\
&\leq \int^t_s\big(\sqrt{|K|}\,\|\frac{\rho_n}{\sqrt{\mu_n(\rho_n)}}(t,\cdot)\|_{L^\infty(K)}\|\sqrt{\mu_n(\rho_n)}\p_x u_n(t,\cdot)\|_{L^2(K)}\\
&+\|\p_x\va_{1,n}(\rho_n)(t,\cdot)\|_{L^1(K)}\|\frac{\rho_n}{\mu_n(\rho_n)}(t,\cdot)\|_{L^\infty(K)}\| u_n(t,\cdot)\|_{H^1(\R)}\big) ds\\
&\leq C\sqrt{t-s}(1+\sqrt{|K|}).
\end{aligned}
\end{equation}
From the Lemma \ref{compac} by localizing the argument (it suffices to consider the sequel $(\varphi\rho_n)_{n\in\mathbb{N}}$ with $\varphi\in C^\infty_0(\R)$ and to use the diagonal process), we deduce that up to a subsequence $(\rho_n)_{n\in\mathbb{N}}$ converges strongly in $L^1_{loc}([0,T_\beta]\times\R)$. It implies in particular that up to a subsequence  $(\rho_n)_{n\in\mathbb{N}}$ converges almost everywhere to $\rho$ on $[0,T_\beta]\times\R$. From (\ref{controlunif}) and the dominated convergence theorem we deduce that:
\begin{itemize}
\item $\rho_n^\gamma$, $\rho_n^\alpha$ converges strongly in $L^p_{loc}([0,T_\beta]\times\R)$ for any $1\leq p<+\infty$ to respectively $\rho^\gamma$ and $\rho^\alpha$.
\end{itemize}
 From (\ref{energiea}) and (\ref{controlunif}) we observe that  $u_n$ is uniformly bounded in $L^2([0,T_\beta],H^1(\R))$. Up to a subsequence $(u_n)_{n\in\mathbb{N}}$ converges weakly to $u$ in $L^2([0,T_\beta],H^1(\R))$. Since $\rho_n^\alpha$ converges strongly to $\rho^\alpha$ in $L^2_{loc}([0,T_\beta]\times\R)$, $\mu_n(\rho_n)\p_x u_n$ converges to $\mu(\rho)\p_x u$ in the sense of the distribution on $[0,T_\beta]\times\R$.\\
 Easily we can prove now that:
\begin{itemize}
\item $(\rho_n u_n)_{n\in\mathbb{N}}$ converges up to a subsequence to $\rho u$ in ${\cal D}'([0,T_{\beta}]\times\R)$.
\item $\rho_n u^2_n$ converges up to a subsequence to $\rho u^2$ in ${\cal D}'([0,T_{\beta}]\times\R)$.
\end{itemize}
Indeed it suffices to use the arguments developed in \cite{Lio98}. To finish since $\p_x\rho_n$ is uniformly bounded in $L^\infty([0,T_\beta],L^1(\R))$ we deduce that up to a subsequence 
$\p_x\rho_n$ converges * weakly to $\p_x\rho$ in $L^\infty([0,T_\beta],{\cal M}(\R))$.\\
In addition we have the following estimates for $C>0$ large enough by using the Fatou lemma for weak convergence:
\begin{equation}
\begin{cases}
\begin{aligned}
&\|u\|_{L^\infty([0,T_\beta],L^2(\R))}\leq C\\
&\|\p_x u\|_{L^2([0,T_\beta],L^2(\R))}\leq C\\
&\|(\rho,\frac{1}{\rho})\|_{ L^\infty([0,T_\beta],L^\infty(\R))}\leq C\\
&\| \p_x\rho\|_{ L^\infty([0,T_\beta],{\cal M}(\R))}\leq C\\
&\|(\rho-\bar{\rho})\|_{L^\infty([0,T_\beta],L^\gamma_2(\R))}\leq C.
%&\frac{1}{2}\int_{\R}\big(\rho(T,x)|u(T,x)|^2+(\Pi(\rho(T,x))-\Pi(\bar{\rho})) \big) dx\\
%&+\int^T_0\int_{\R}\mu_n(\rho(s,x))(\p_x u(s,x))^2 dx\leq \frac{1}{2}\int_{\R}\big(\rho_n(0,x)|u_n(0,x)|^2+(\Pi(\rho_n(0,x))-\Pi(\bar{\rho})) \big) dx.
\end{aligned}
\end{cases}
\end{equation}
In addition from (\ref{energiea}), (\ref{controlunif}), (\ref{controlunif1}) and Sobolev embedding, we deduce that $(\rho_n-\bar{\rho})_{n\in\mathbb{N}}$ is uniformly bounded in $L^\infty([0,T_\beta], H^{\frac{1}{2}}(\R))$.
From the mass equation $\p_t\rho_n$ is uniformly in $L^\infty([0,T_\beta],W^{-1,1}(\R))$. Using the Lemma \ref{Aubin}, we deduce that $(\rho_n)_{n\in\mathbb{N}}$ converges strongly to $\rho$ in $C([0,T_\beta], L^p_{loc}(\R))$ for $1\leq p<+\infty$.\\
Following the same argument than in the previous section, we prove also that $(\rho_n u_n)_{n\in\mathbb{N}}$ converges to $\rho u$ in $C([0,T_\beta],{\cal M}(\R)^*)$.
%\begin{equation}
%\begin{aligned}
%&u\in L^\infty([0,T_\beta],L^2(\R))\\
%&\n u\in L^2([0,T_\beta],L^2(\R))\\
%&\rho,\frac{1}{\rho}\in L^\infty([0,T_\beta],L^\infty(\R))\\
%& \p_x\rho\in L^\infty([0,T_\beta],{\cal M}(\R))\\
%&(\rho-\bar{\rho})\in L^\infty([0,T_\beta],L^\gamma_2(\R))
%&\frac{1}{2}\int_{\R}\big(\rho(T,x)|u(T,x)|^2+(\Pi(\rho(T,x))-\Pi(\bar{\rho})) \big) dx\\
%&+\int^T_0\int_{\R}\mu_n(\rho(s,x))(\p_x u(s,x))^2 dx\leq \frac{1}{2}\int_{\R}\big(\rho_n(0,x)|u_n(0,x)|^2+(\Pi(\rho_n(0,x))-\Pi(\bar{\rho})) \big) dx.
%\end{aligned}
%\end{equation}
\section{Proof of the theorem \ref{theo2}}
\label{section4}
%POURQUOI VHANGER LA REGULARISATION SUR DONNEES INITIALES, FAIRE!\\
Let $(\rho_0,u_0,v_0)$ verifying the assumptions of the Theorem \ref{theo2}, we define
now the regularizing initial data
$(\rho_0^n,m_0^n)_{n\in\mathbb{N}}$ as in the previous section. From \cite{Ka,MV} since $(\rho_0^n-\bar{\rho},u_0^n)$ belongs to $H^1(\R)\times H^1(\R)$ and we have $0\leq c\leq \rho_0^n\leq M<+\infty$, we know that it exists a sequence of global strong solution $(\rho_n, u_n)_{n\in\mathbb{N}}$ for the system (\ref{1}) with constant viscosity coefficient $\mu$. In addition since the regularity $H^s$ with $s\geq 1$ is preserved the sequence $(\rho_n,u_n)_{n\in\mathbb{N}}$ is in $C^\infty(\R^+\times\R)$. We are going now to prove uniform estimates in $n$ on the solution $(\rho_n,u_n)$. When we will have sufficiently uniform estimates, we will verify that $(\rho_n,u_n)_{n\in\mathbb{N}}$ converges up to a subsequence in a suitable functional space to a global weak solution $(\rho,u)$ of (\ref{1}).\\
Similarly we have the following entropy and for $C>0$ large enough:
\begin{equation}
\begin{aligned}
&\frac{1}{2}\int_{\R}(\rho_n |v_n|^2(T,x)+(\Pi(\rho_n)-\Pi(\bar{\rho}))(T,x) dx+\frac{4a\gamma\mu}{(\gamma+\alpha-1)^2} \int^T_0\int_{\R}|\p_x\rho_n^{\frac{1}{2}(\gamma-1})|^2 (s,x) ds dx\\
&\hspace{4cm}\leq \frac{1}{2}\int_{\R}(\rho_n |v_n|^2(0,x)+(\Pi(\rho_n)-\Pi(\bar{\rho}))(0,x) dx\leq C.
\end{aligned}
\label{BD1ac}
\end{equation}
\subsubsection*{$L^\infty$ estimate of $\frac{1}{\rho}$ and $\rho$ in $L^\infty$ norm on a finite time intervall}
In this case the only difference is the way that we are going to use to prove  $L^\infty$ estimate on $\rho_n$ and $\frac{1}{\rho_n}$. We can show the estimate (\ref{cru1}) in a similar way as the proof of the Theorem \ref{theo1}. We deduce as (\ref{Linf}) that we have for $C>0$ large enough independent on $n$ and any $t>0$:
\begin{equation}
\begin{aligned}
&\|\rho_n u_n(t,\cdot)\|_{L^1(\R)}+\|\rho_n v_n(t,\cdot)\|_{L^1(\R)}\\
&\hspace{2cm}\leq C(  \|\rho_n v_n(0,\cdot)\|_{L^1(\R)}+\|\rho_n u_n(0,\cdot)\|_{L^1(\R)} ) e^{\frac{C}{\mu}\int^t_0\|P'(\rho_n)\rho_n (s,\cdot)\|_{L^\infty}ds}.
\end{aligned}
\label{pas1}
\end{equation}
Since we have $\rho_n v_n=\rho_nu_n+\mu\p_x\ln\rho_n$, we deduce in particular that for any $t>0$ we have using (\ref{BD1ac}):
\begin{equation}
\begin{aligned}
\|\ln\rho_n(t,\cdot)\|_{L^\infty(\R)}&\leq|\ln\bar{\rho}|+\|\p_x\ln(\rho_n(t,\cdot))\|_{L^1(\R)}\\
&\leq C\big(|\ln\bar{\rho}|+\big(\|\rho_0v_0\|_{{\cal M}(\R)}+\|\rho_0 u_0\|_{{\cal M}(\R)}\big)e^{\frac{C}{\mu}\int^t_0\|P'(\rho_n)\rho_n(s,\cdot)\|_{L^\infty(\R)}ds}\big).
\end{aligned}
\label{5.123}
\end{equation}
An it yields that for $C'>0$:
\begin{equation}
\begin{aligned}
\|\ln\rho_n(t,\cdot)\|_{L^\infty(\R)}
&\leq C'\big(|\ln\bar{\rho}|+\big(\|\rho_0v_0\|_{{\cal M}(\R)}+\|\rho_0 u_0\|_{{\cal M}(\R)}\big)e^{C' t e^{\gamma\|\ln\rho_n\|_{L^\infty_t(L^\infty(\R))}} }\big).
\end{aligned}
\label{Linf1ad}
\end{equation}
Let us prove now that $(\rho_n)_{n\in\mathbb{N}}$ and $(\frac{1}{\rho})_{n\in\mathbb{N}}$ are uniformly bounded in $L^\infty$ in $n$ and on an time interval $[0,T^*]$ with $T^*$ independent on $n$. More precisely we define by:
$$
\begin{aligned}
&T_n=\sup\{ t\in(0,+\infty),%\|\rho^n(t,\cdot)\|_{L^\infty}\leq \sup (2\|\rho_0\|_{L^\infty}, C'\bar{\rho}+2C'(\|\rho_0 v_0\|_{L^1(\R)}+\|\rho_0 u_0\|_{L^1(\R)} )^{\frac{1}{\alpha}})
\|\ln\rho^n(t,\cdot)\|_{L^\infty(\R)}
\leq M_1.\},
\end{aligned}
$$
with $M_1= \sup(2\|\ln\rho_0\|_{L^\infty}, C'\big(|\ln\bar{\rho}|+2(\|\rho_0 v_0\|_{L^1(\R)}+\|\rho_0 u_0\|_{L^1(\R)} )\big))$.
We observe that $T_n>0$ since $\ln \rho_n$ belongs to $C([0,+\infty[, L^\infty(\R))$ for any $n\in\mathbb{N}$ (indeed this is due to the regularity of the solution $(\rho_n,u_n)$) and since $\|\ln \rho_0^n\|_{L^\infty}<M_1$). Let us define now $T_2$ such that:
$$e^{C' T_2 e^{\gamma M_1}}=\frac{3}{2}\;\;\;\mbox{and}\;\;\; T_2=\frac{\ln (\frac{3}{2})}{C' e^{\gamma M_1}}.$$
From the definition of $T_n$ and from (\ref{Linf1ad}) we have for any $n\in\mathbb{N}$:
$$T_n\geq T_2>0.$$
We have then proved that $(\ln \rho_n)_{n\in\mathbb{N}}$ is uniformly bounded in $n$ in $L^\infty((0,T_2),L^\infty(\R))$ with $T_2$ independent on $n$. In particular this gives for any $n$:
\begin{equation}
\|(\rho_n,\frac{1}{\rho_n})\|_{L^\infty([0,T_2],L^\infty(\R))}\leq e^{M_1}.
\label{linfad}
\end{equation}
Now proceeding as in the previous section we can show additional regularity information on the velocity $u_n$ and $\rho_n$, we have in particular for any $t>0$ there exists $C>0$ independent on $n$ such that:
\begin{equation}
\begin{aligned}
&\| \sqrt{\rho_n}u_n(t,\cdot)\|_{L^2(\R)}\leq C(\frac{1}{\sqrt{t}}+1)\\
&\|(\sqrt{s}1_{\{s\leq 1\}}+1_{\{s\geq 1\}})\p_x u_n\|_{L^2_t(L^2(\R))}\leq C.
\end{aligned}
\label{cru2ad}
\end{equation}
We recall in addition that from (\ref{BD1ac}) we have also for any $t>0$ and $C>0$ independent on $n$:
\begin{equation}
\begin{cases}
\begin{aligned}
&\|\sqrt{\rho_n}v_n(t,\cdot)\|_{L^2}\leq C\\
&\|\Pi(\rho_n(t,\cdot))-\Pi(\bar{\rho})\|_{L^1(\R)}\leq C\\
&\|\p_x(\rho_n^{\frac{1}{2}(\gamma-1)})\|_{L^2_t(L^2(\R))}\leq C.
\end{aligned}
\end{cases}
\label{superimpad}
\end{equation}
From (\ref{cru2ad}), (\ref{superimpad}) and the definition of the effective velocity $v_n$, we deduce that for any $t>0$ and $C>0$ large enough:
\begin{equation}
\|\p_x(\frac{1}{\sqrt{\rho_n}(t,\cdot)})\|_{L^2(\R)}\leq C(\frac{1}{\sqrt{t}}+1).
\label{superimpadd}
\end{equation}
\subsubsection*{Control of the $L^\infty$ norm for long time of $(\rho_n)_{n\in\mathbb{N}} $ and $(\frac{1}{\rho_n})_{n\in\mathbb{N}}$}
We are going to follow the method develop by Hoff in \cite{Hof98}. We fix now $t_1\in]0,T_2]$. We are going now to prove that for every $T>t_1$, it exists $C(T)>0$ such that for any $x_0\in \R$ and $t_2\in[t_1,T]$ it exists necessary a element $x_1$ such that:
%First for any $t_1<t_2\leq T$ with $t_1\in]0,T_2]$ and $x_0\in\R$, then there exists $x_1$ and a constant $C(T)>0$ depending only on $T$ such that we have:
\begin{equation}
\begin{cases}
\begin{aligned}
&|x_0-x_1|\leq C(T)\\
&C(T)^{-1}\leq \rho_n(t_2,x_1)\leq C(T).
\end{aligned}
\end{cases}
\label{Hoftech}
\end{equation}
It is important to point out that $C(T)$ is independent on $n$ and depends only on $T$.
We observe that there exists $C>0$ such that:
$$\lim_{\rho\rightarrow 0}\inf(\Pi(\rho)-\Pi(\bar{\rho}))\geq C^{-1}.$$
We deduce in particular that there exists $\delta>0$ such that for all $0<\rho<\delta$ we have:
$$\Pi(\rho)-\Pi(\bar{\rho})\geq\frac{C^{-1}}{2}.$$
Then if for any $R>0$ we have:
$$\sup_{x\in[x_0-R,x_0+R]}\rho_n(t_2,x)<\delta,$$
then we have from the previous estimate and (\ref{BD1ac}) (with $C_1>0$ depends on the initial data $(\rho_0,u_0)$):
\begin{equation}
C^{-1}R\leq \int_{\R}(\Pi(\rho_n(t_2,x)-\Pi(\bar{\rho}))dx\leq C_1.
\label{vgtech}
\end{equation}
In particular $R$ can not goes to $0$. In particular it implies that for $R_1=2 C_1 C$, it exists necessary $x_2\in [x_0-R_1,x_0+R_1]$ such that $\rho_n(t_2,x_2)>\delta$. Let us define now $C_{n,\delta}=\{x\in [x_0-R_1,x_0+R_1],\,\rho_n(t,x)>\delta\}$, it is now obvious that the measure of $C_{n,\delta}$ verifies $|C_{n,\delta}|> \frac{ R_1}{2}$. Indeed we have from (\ref{BD1ac}) for $C_1>0$ large enough:
\begin{equation}
\int_{[x_0-R_1,x_0+R_1]/ C_{n,\delta}}(\Pi(\rho_n(t_2,x)-\Pi(\bar{\rho}))dx+\int_{C_{n,\delta}}(\Pi(\rho_n(t_2,x)-\Pi(\bar{\rho}))dx \leq C_1.
\label{vgtecha}
\end{equation}
If $|C_{n,\delta}|<\frac{R_1}{2}$ then we have:
$$\frac{C^{-1}}{2}\big|[x_0-R_1,x_0+R_1]/ C_{n,\delta}\big|\leq C_1.$$
This is absurd since  $R_1=2 C_1 C$ then we have $|C_{n,\delta}|>\frac{R_1}{2}$. We can verify now that for $C_2>0$ large enough we have:
$$\rho+P(\rho)\leq C_2[1+(\Pi(\rho)-\Pi(\bar{\rho}))].$$
It implies from (\ref{vgtech}) that it exists $x_1\in C_{n,\delta}$ such that $\rho_n(t_2,x_1)\leq(4C_2+\frac{2C_2C_1}{R_1})$. This concludes the proof of (\ref{Hoftech}) since $R_1$ depends only on $T$.
\\
We use now the characteristic method, we set for  $t\in[t_1,t_2]$ with $t_1<t_2\leq T$:
$$
\begin{cases}
\begin{aligned}
&\frac{d}{dt}X^n_j(s)=u_n(t,X^n_j(s))\\
&X^n_j(t_2)=x_j\in\R,\;\;j=0,1.
\end{aligned}
\end{cases}
$$
We set $\D L_n(t)=\log\rho_n (t,X_1(t))-\log\rho_n(t,X_0(t))$. From (\ref{1}) we deduce that:
$$
\begin{aligned}
&\frac{d}{dt} \D L_n(t)=-\int^{X^n_1(t)}_{X^n_0(t)}\p_{xx}u_n(t,x) dx\\
&=-\frac{1}{\mu}\int^{X^n_1(t)}_{X^n_0(t)}(\p_t(\rho_n u_n)(t,x)+\p_x(\rho_n u_n^2)(t,x)+\p_x P(\rho_n)(t,x) dx=-\frac{d}{dt} I_n(t)-\D P_n(t),% (we mention that $\alpha$ is well defined since we considered approximate solution $\rho^n$ which does not vanish). 
\end{aligned}
$$
where $I_n(t)=\int^{X^n_1(t)}_{X^n_0(t)} \rho_n u_n(t,x) dx$ and $\D P_n(t)=P(\rho_n(t,X^n_1(t)))-P(\rho_n(t,X^n_0(t)))$. We observe now as in \cite{Hof98} that $\alpha_n(t)=\frac{\D L_n(t)}{\D P_n(t)}$ is positive. We deduce then that we have the following ordinary differential equation:
$$
\frac{d}{dt}\D L_n(t)+\alpha_n(t)\D L_n(t)=-\frac{d}{dt}I_n(t).
$$
We deduce that we have for any $t\in [t_1,t_2]$:% $t\geq t_1$ with $t_1\in (0,T_2)$:
\begin{equation}
\begin{aligned}
&|\D L_n(t)|\leq |\D L_n(t_1)|+|I_n(t_1)|+|I(t)|+\int^t_{t_1} e^{-\int^t_s \alpha_n(\tau) d\tau}\alpha_n(s)|I_n(s)| ds.
\end{aligned}
\label{2.8}
\end{equation}
For the latter, we take for example $X^n_1>X^n_0$. We observe easily that for any $t\in[t_1,t_2]$ we have from (\ref{cru2ad}), (\ref{BD1ac}) and for $C>0$ large enough:
\begin{equation}
\begin{aligned}
&I_n(t)\leq \|\sqrt{\rho_n}u_n(t,\cdot)\|_{L^2(\R)}\|\sqrt{\rho_n}\|_{L^2([X^n_0(t),X^n_1(t)])}\\
&\leq C(1+\frac{1}{\sqrt{t}})\big(\|\rho_n(t,\cdot)-\bar{\rho}\|_{L^\gamma_2(\R)}(|X^n_0(t)-X^n_1(t)|^{\frac{1}{2}}+|X^n_0(t)-X^n_1(t)|^{1-\frac{1}{\gamma}})\\
&\hspace{9cm}+|X^n_0(t)-X^n_1(t)|\big)^{\frac{1}{2}}\\
&\leq C(1+\frac{1}{\sqrt{t}})\big(|X^n_0(t)-X^n_1(t)|^{\frac{1}{2}}+|X^n_0(t)-X^n_1(t)|^{1-\frac{1}{\gamma}})+|X^n_0(t)-X^n_1(t)|\big)^{\frac{1}{2}}.
%(1+\|\rho(t,\cdot)-\bar{\rho}\|_{L^\gamma_2})\\
%&\hspace{2cm}\times((|X_0(t)-X_1(t)|^{\frac{1}{2}}+|X_0(t)-X_1(t)|^{1-\frac{1}{\gamma}}+|X_0(t)-X_1(t)|\big).
\end{aligned}
\label{adimp}
\end{equation}
It remains now only to estimate $|X^n_0(t)-X^n_1(t)|$, next we have:
$$
\begin{aligned}
&\frac{d}{dt}(X^n_1-X^n_0)(t)=\int^{X^n_1(t)}_{X^n_0(t)}\p_x u_n (t,x) dx\geq -\frac{1}{2}(X^n_1(t)-X^n_0(t))-\frac{1}{2}\int^{X^n_1(t)}_{X^n_0(t)}(\p_x u_n)^2 (t,x) dx.
\end{aligned}
$$
Using (\ref{cru2ad}) and the previous inequality, we get for $t_1\leq t\leq t_2$ and $C>0$ large enough:
$$
\begin{aligned}
&X^n_1(t_2)-X^n_0(t_2)\geq X^n_1(t)-X^n_0(t)- \int^{t_2}_{t}(X_1(s)-X_0(s))ds-\int^{t_2}_{t} \int^{X^n_1(s)}_{X^n_0(s)}(\p_x u_n)^2 (s,x) dx ds\\
&X^n_1(t)-X^n_0(t)\leq  |x_1-x_0|+ \int^{t_2}_{t}(X^n_1(s)-X^n_0(s))ds+\int^{t_2}_{t} \int_{\R} (\p_x u_n)^2 (s,x) dx\\
&\leq  |x_1-x_0|+ \int^{t_2}_{t}(X^n_1(s)-X^n_0(s))ds+\frac{C}{t_1} .
\end{aligned}
$$
Setting $F_n(t)= \int^{t_2}_{t}(X^n_1(s)-X^n_0(s))ds$, we have then:
$$F_n(t)+F'_n(t)\geq -|x_0-x_1|-\frac{C}{t_1},$$
and:
$$F_n(t)\leq\int^{t_2}_t(|x_0-x_1|+\frac{C}{t_1})e^{s-t} ds.$$
%From Gronwall lemma and (\ref{Hoftech}) 
We deduce that for $t_1\leq t\leq t_2\leq T$, we have:
\begin{equation}
|X^n_1(t)-X^n_0(t)|\leq  |x_1-x_0|+\frac{C}{t_1} +\int^{T}_t(|x_0-x_1|+\frac{C}{t_1})e^{s-t} ds.
% C_1(T,t_1),
\label{ad1}
\end{equation}
Combining (\ref{adimp}) and (\ref{ad1}), it exists a continuous fonction $C$ independent on $n$ such that for $t_1\leq t\leq t_2\leq T$ we have:
\begin{equation}
I_n(t)\leq C(T,t_1).
\label{Linfimp}
\end{equation}
Now using (\ref{2.8}) and (\ref{Linfimp}) we get for any $t_1\leq t\leq t_2\leq T$:
\begin{equation}
\begin{aligned}
|\D L_n(t)|&\leq |\D L_n(t_1)|+2 C(T,t_1)+C(T,t_1)\int^t_{t_1} e^{-\int^t_s \alpha_n(\tau) d\tau}\alpha_n(s)ds\\
&\leq |\D L_n(t_1)|+4 C(T,t_1).
\end{aligned}
\label{2.9}
\end{equation}
Now since $(\frac{1}{\rho_n}(t_1,\cdot),\rho_n(t_1,\cdot))_{n\in\mathbb{N}}$ is uniformly bounded in $n$ in $L^\infty(\R)$, we deduce that:
\begin{equation}
|\ln\rho_n(t_2,x_0)-\ln\rho_n(t_2,x_1)|\leq C(T,t_1).
\label{Hoftech1}
\end{equation}
From (\ref{Hoftech}) and (\ref{Hoftech1}), we deduce that for any $t_2$ with $0<t_1\leq t_2\leq T$ and any $x_0\in\R$ there exists $C_1(T,t_1)$ independent on $n$ such that:
$$C(t_1,T)^{-1}\leq\rho_n(t_2,x_0)\leq C(T,t_1).$$
It implies in particular that we have for $C_1$ an increasing continuous function and any $t\geq t_1$
\begin{equation}
\begin{aligned}
&\|(\rho_n(t,\cdot),\frac{1}{\rho_n}(t,\cdot)\|_{L^\infty}\leq C_1(t,t_1).
\end{aligned}
\label{2.10}
\end{equation}
If we summarize the previous estimates (\ref{pas1}), (\ref{cru2ad}), (\ref{superimpad}), (\ref{superimpadd}) and (\ref{2.10}), we have proved that for any $t>0$ there exists $C>0$ and $C_1$ a continuous increasing function  not depending on $n$ such that:
\begin{equation}
\begin{cases}
\begin{aligned}
&\|(\rho_n u_n(t,\cdot),\rho_n v_n(t,\cdot))\|_{L^1(\R)}\leq C_1(t)\\
&\|(\rho_n(t,\cdot),\frac{1}{\rho_n}(t,\cdot)\|_{L^\infty}\leq C_1(t)\\
&\| \sqrt{\rho_n}u_n(t,\cdot)\|_{L^2(\R)}\leq C(\frac{1}{\sqrt{t}}+1)\\
&\|(\sqrt{s}1_{\{s\leq 1\}}+1_{\{s\geq 1\}})\p_x u_n\|_{L^2_t(L^2(\R))}\leq C\\
&\|\sqrt{\rho_n}v_n(t,\cdot)\|_{L^2}\leq C\\
&\|\rho_n(t,\cdot)-\bar{\rho}\|_{L^\gamma_2(\R)}\leq C\\
&\|\p_x\rho_n\|_{L^2_t(L^2(\R))}\leq C\\
%&\|\p_x(\rho^{\frac{1}{2}(\gamma-1)})\|_{L^2_t(L^2(\R))}\leq C\\
&\|\p_x\rho_n(t,\cdot)\|_{L^2(\R)}\leq C(\frac{1}{\sqrt{t}}+1).
\end{aligned}
\end{cases}
\label{resumad}
\end{equation}
\subsubsection*{Compactness arguments}
From (\ref{resumad}) we deduce %that these estimates are uniform in $n$ and it implies in particular since $\|\rho_n-\bar{\rho}\|_{L^\infty_{loc}(L^2(\R))}$ and $\|\rho_n\|_{L^\infty_{loc}(L^\infty(\R))}$ is uniformly bounded in $n$ that 
%$\rho_n-\bar{\rho}$ is uniformly bounded in $L^\infty_{loc}(L^2(\R))$. Combining this information with (\ref{resumad}), we obtain 
that $(\rho_n)_{n\in\mathbb{N}}$ is uniformly bounded in
$L^2_{loc}(H^1(\R))$. In addition from (\ref{resumad}) and from the fact that $\p_t\rho_n=-\p_x(\rho_n u_n)$ we deduce that $(\p_t\rho_n)_{n\in\mathbb{N}}$ is uniformly bounded in $L^{\infty}_{loc}(W^{-1,1}(\R))$ and is bounded uniformy in $n$ in $L^{2-\e}_{loc}(W^{-\frac{3}{2},2}(\R))$. Using the Aubin-Lions lemma \ref{Aubin}, we deduce that the sequence $(\rho_n)_{n\in\mathbb{N}}$ converges up to a subsequence in $L^2_{loc}(H^s_{loc}(\R))$ with $0\leq s<1$ to $\rho$.
It implies in particular that up to a subsequence $(\rho_n)_{n\in\mathbb{N}}$ converges almost everywhere to $\rho$. From convergence dominated and (\ref{resumad}) we show that  $\rho_n^\gamma$  converges to $\rho^\gamma$ in $L^p_{loc}(\R\times\R^+)$ for any $1\leq p<+\infty$. In addition from (\ref{resumad}) for any $t>0$, there exists $C(t)>0$ such that we have:
\begin{equation}
\begin{aligned}
&\|(\frac{1}{\rho}(t,\cdot),\rho(t,\cdot))\|_{L^\infty(\R)}\leq C(t).
\end{aligned}
\end{equation}
% by Sobolev embedding that $\rho_n$ converges uniformly to $\rho$ on every compact $[0,T]\times [-M,M]$ with $M>0$. We deduce that $\rho_n^\gamma$ converges uniformy to $\rho^\gamma$  on every compact $[0,T]\times [-M,M]$ with $M>0$.\\
Proceeding as for (\ref{cru3}), we deduce that $(u_n)_{n\in\mathbb{N}}$ is uniformly bounded in $n$ in $L^{p(s)-\e}_{loc}(H^s(\R))$ with $0\leq s<1$ and $p(s)=\frac{3}{1+2s}$ with $\e>0$ sufficiently small such that
$p(s)-\e\geq 1$. Now using paraproduct law we deduce that we have for any $t>0$ and $0\leq s\leq\frac{1}{2}$:
\begin{equation}
\begin{aligned}
&\|(\rho_n-\bar{\rho}) u_n(t,\cdot)\|_{H^s(\R)}\leq\\
& \|u_n(t,\cdot)\|_{H^s(\R)}\|\rho_n(t,\cdot)-\bar{\rho}\|_{L^\infty}+\|\rho_n(t,\cdot)-\bar{\rho}\|_{H^1(\R)}\|u_n(t,\cdot)\|_{H^{-\frac{1}{2}}(\R)}
\end{aligned}
\end{equation}
Combining the fact that $(u_n)_{n\in\mathbb{N}}$ is uniformly bounded in $n$ in $L^{p(s)-\e}_{loc}(H^s(\R))$ and (\ref{resumad}) we deduce that $(\rho_n u_n)_{n\in\mathbb{N}}$ is uniformly bounded in $L^{q(s)}_{loc}(H^s(\R))$
with $q(s)=\min (2,p(s))$ and $s\in[0,\frac{1}{2}]$.\\
Now we recall that we have $\p_t(\rho_n u_n)=-\p_x(\rho_n u^2_n)+\mu\p_{xx}u_n-\p_x P(\rho_n)$, taking in the previous estimate $s=0$ we deduce that $(u_n^2)_{n\in\mathbb{N}}$ is bounded in $L^{\frac{3}{2}}_{loc}(L^1(\R))$ and $(\rho_n)_{n\in\mathbb{N}}$ is uniformly bounded in $L^\infty_{loc}(L^\infty(\R))$ then $(\rho_n u_n^2)_{n\in\mathbb{N}}$ is uniformly bounded in $L^{\frac{3}{2}}_{loc}(L^1(\R))$. Proceeding similarly we observe by Sobolev embedding that $(\p_t(\rho_nu_n))_{n\in\mathbb{N}}$ is bounded uniformly in $L^1_{loc}(H^{-s_2}(\R))$ with $s_2>0$ large enough. Using the Aubin-Lions lemma \ref{Aubin} we deduce that $(\rho_n u_n)_{n\in\mathbb{N}}$ converges strongly up to a subsequence to $m$
in in $L^{q(s)}_{loc}(H_{loc}^s(\R))$
with $q(s)=\min (2,p(s))$ and $s\in[0,\frac{1}{2}[$. By Sobolev embedding we deduce that $(\rho_n u_n)_{n\in\mathbb{N}}$ converges strongly to $m$ in $L^2_{loc}(L^2_{loc}(\R))$. Now since $(\rho_n)_{n\in\mathbb{N}}$ and $(\frac{1}{\rho_n})_{n\in\mathbb{N}}$ converges respectively to $\rho$ and $\frac{1}{\rho}$ in $L^p_{loc}(\R^+\times\R)$ for $1\leq p<+\infty$ and since $(\frac{1}{\rho_n})_{n\in\mathbb{N}}$ is uniformly bounded in $L^\infty_{loc}(L^\infty(\R))$ we get that $(u_n)_{n\in\mathbb{N}}$ converges up to a subsequence strongly to $\frac{m}{\rho}=u$ in  $L^{2-\e}_{loc}(L^{2-\e}_{loc}(\R))$ with $\e>0$ and $(u_n)_{n\in\mathbb{N}}$ converges up to a subsequence almost everywhere to $u$.
\\
All these estimates are now sufficient to show that $(\rho_n,u_n)_{n\in\mathbb{N}}$ converges up to a subsequence in the sense of the distribution to a global weak solution $(\rho,u)$.\\
In addition using Fatou lemma, convergence almost everywhere and weak convergence, we deduce from (\ref{resumad}) that for any $t>0$ there exists $C>0$ and $C_1(t)>0$ such that:
\begin{equation}
\begin{cases}
\begin{aligned}
&\|\rho u(t,\cdot)\|_{L^1(\R)}\leq C_1(t)\\
%&\|\rho v(t,\cdot))\|_{{\cal M}(\R)}\leq C_1(t)\\
&\|(\rho(t,\cdot),\frac{1}{\rho}(t,\cdot)\|_{L^\infty}\leq C_1(t)\\
&\| \sqrt{\rho}u(t,\cdot)\|_{L^2(\R)}\leq C(\frac{1}{\sqrt{t}}+1)\\
%&\|(\sqrt{s}1_{\{s\leq 1\}}+1)\p_x u\|_{L^2_t(L^2(\R))}\leq C\\
%&\|\sqrt{\rho}v(t,\cdot)\|_{L^2}\leq C\\
&\|\rho(t,\cdot)-\bar{\rho}\|_{L^\gamma_2(\R)}\leq C\\
&\|\p_x\rho\|_{L^2_t(L^2(\R))}\leq C.
%&\|\p_x(\rho^{\frac{1}{2}(\gamma-1)})\|_{L^2_t(L^2(\R))}\leq C\\
%&\|\p_x\rho(t,\cdot)\|_{L^2(\R)}\leq C C_1(t) (\frac{1}{\sqrt{t}}+1).
\end{aligned}
\end{cases}
\label{resumad2d}
\end{equation}
In addition $u$ belongs to $L^{p(s)-\e}_{loc}(H^s(\R))$ with $0\leq s<1$ and $p(s)=\frac{3}{1+2s}$ with $\e>0$ sufficiently small such that
$p(s)-\e\geq 1$. As in the previous proof, we can show also that $\rho$ is in  $C([0,+\infty[,W^{-\e,p}_{loc}(\R))$ and $\rho u$ is in $C([0,T],{\cal M}(\R)^*)$. This achieves the proof of the Theorem \ref{theo2}. %DONNER L'ARGUMENT DE CONTINUITE EN TEMPS $0$.
%\|\rho u(t,\cdot)\|_{L^1(\R)}\leq (\|\rho_0 v_0\|_{L^1(\R)}+\|\rho_0 u_0\|_{L^1(\R)} ) e^{2\int^T_0\|\frac{P'(\rho)\rho^2}{\mu(\rho)}(s,\cdot)\|_{L^\infty}ds}.
%(VOIR L'ARGUMENT QUI PERMET D'ASSURER QUE L'ON EST DANS $L^1(\R)$ INSTANTANEMENT, PAS CLAIR? VOIRE POUR EULER INCOMPRESSIBLE.
\section*{Acknowledgements}
The author has been partially funded by the ANR project 
INFAMIE ANR-15-CE40-0011. This work was realized during the secondment of the author in the ANGE Inria team.

\end{document}